\documentclass[11pt,american,oneside]{book}
\usepackage[T1]{fontenc}
\usepackage{babel}
\usepackage{csquotes}
\usepackage{amsmath,amssymb,amsthm}
\usepackage{booktabs,array}
\usepackage{xcolor}
\usepackage[left=2.5cm, right=2.5cm, top=2.6cm, bottom=2.6cm]{geometry}
\usepackage{fancyhdr}
\usepackage[colorlinks=true,linkcolor=blue!50!black,citecolor=blue!50!black,urlcolor=blue!50!black]{hyperref}

\DeclareMathOperator{\cl}{cl}
\DeclareMathOperator{\cone}{cone}

\DeclareMathOperator{\rank}{rank}
\DeclareMathOperator{\coker}{coker}
\DeclareMathOperator{\ad}{ad}

\newcommand{\R}{\mathbb{R}}
\newcommand{\Z}{\mathbb{Z}}
\newcommand{\N}{\mathbb{N}}
\newcommand{\cB}{\mathcal{B}}
\newcommand{\cD}{\mathcal{D}}
\newcommand{\cH}{\mathcal{H}}
\newcommand{\cK}{\mathcal{K}}
\newcommand{\cP}{\mathcal{P}}
\newcommand{\cC}{\mathcal{C}}
\newcommand{\cN}{\mathcal{N}}
\newcommand{\Rlp}{\R[s,s^{-1}]}
\newcommand{\Rp}{\R_+[s,s^{-1}]}
\newcommand{\us}{\_\allowbreak}

\theoremstyle{plain}
\newtheorem{theorem}{Theorem}[chapter]
\newtheorem{proposition}[theorem]{Proposition}
\newtheorem{lemma}[theorem]{Lemma}
\newtheorem{corollary}[theorem]{Corollary}
\newtheorem{conjecture}[theorem]{Conjecture}
\theoremstyle{definition}
\newtheorem{definition}[theorem]{Definition}
\newtheorem{example}[theorem]{Example}
\newtheorem{question}[theorem]{Question}
\newtheorem{problem}[theorem]{Problem}
\theoremstyle{remark}
\newtheorem{remark}[theorem]{Remark}

\title{A Categorical Approach to Four Problems in Control Theory\\[6pt]
\large Diagrams, Degenerate Sets, and Universal Substitutes}
\author{Alex Borisevich}
\date{November 2026}
\begin{document}

\maketitle
\tableofcontents
\chapter{Introduction: from diagrams to four problems}\label{ch:intro}

\section{The programme and the selection}

In 2004 Blondel and Megretski published sixty-three open problems in mathematical systems and control theory, each contributed by a specialist with a statement, a motivation, and the state of the art \cite{BM04}. This volume starts from a top-down proposal rather than from four independent calculations. A control problem is first drawn as a diagram between a category of representations and a category of behaviours or geometric objects. The question is then whether a semantic fact has a syntactic witness---a lift, factorization, section, or splitting. On the regular part the required arrow is obtained by inversion. On the degenerate set it can fail, and the diagram tells us to look for the defect or for a universal substitute.

The four problems were selected as tests of that programme. Their defects are, respectively, non-fullness of a functor, failure of a smooth lift, a nonzero cokernel, and failure of exactness after passing to an associated graded object. The diagrammatic framework is developed in Chapter~\ref{ch:diag} before the four technical chapters. It produces the full algebraic mechanism for Problems~1.10 and 3.2, the closure criterion and search target for Problem~1.9, and the shape and location of a counterexample for Problem~7.3; problem-specific analysis then supplies what the general structure cannot.

The four problems and their outcomes:

\begin{itemize}
\item \textbf{Problem 1.9} (ten Dam, Nieuwenhuis): a Farkas lemma for behavioral inequalities. \emph{The conjecture is false.} Chapter~\ref{ch:fa}.
\item \textbf{Problem 7.3, first question} (Jouan): does a $C^r$ observer force a finite-order $C^r$ relation among the output derivatives? \emph{No, in the smooth category.} Chapter~\ref{ch:ob}.
\item \textbf{Problem 3.2} (Devanathan): remove all nonlinear terms of a linearly controllable single-input system by state change and feedback, to arbitrary order. \emph{Impossible for $n \ge 3$; the obstruction space is computed exactly and its kernel identified.} Chapter~\ref{ch:nf}.
\item \textbf{Problem 1.10} (Trentelman): characterize the behaviours implementable from a given plant by \emph{regular feedback}, an interconnection in which both output cardinalities and McMillan degrees add. \emph{For full interconnection, the 2010 Vinjamoor--Belur solution is recovered as a graded direct-summand criterion and a finite leading-matrix rank test. For partial interconnection, the criterion is completed by identifying the exact defect introduced by elimination.} Chapter~\ref{ch:rfi}.
\end{itemize}

The first two are new results. The third is a confirmation of known normal-form theory in the coordinates of the problem, with one small addition (the identification of the kernel as the prolonged flat-output reparametrization); it is included because it completes the account of the problem and because it is the case where the method reduces to a mechanical computation. For the fourth, Vinjamoor and Belur already solved the full-interconnection problem, constructed an RFI controller when one exists, and parametrized all such controllers \cite{VB10}. The contribution here is an equivalent filtered-module formulation and structural proof of that case, together with a necessary-and-sufficient partial-interconnection criterion that makes the elimination obstruction explicit.

Chapter~\ref{ch:diag} therefore comes first as the source of the four lines of investigation, not as a retrospective translation. It states the syntax--semantics diagram, identifies the four failure modes, and records exactly which conclusions follow from structure alone. The problem chapters then perform the required analysis and computation. Chapter~\ref{ch:formal} describes the Lean~4/Mathlib formalization: what is machine-checked, what is not, and why.

\section{The problems as stated}

We reproduce the questions in the form in which the chapters answer them. The original chapters of \cite{BM04} contain more context and are cited in each chapter's literature section.

\begin{problem}[1.9, ten Dam--Nieuwenhuis]\label{prob:farkas}
Let $\Rlp$ be the ring of real Laurent polynomials and $\sigma$ the shift on sequences $w : \Z\to\R^q$, $(\sigma w)(t) = w(t+1)$. For $R \in \R^{g\times q}[s,s^{-1}]$ put $\cB(R) = \{w : R(\sigma,\sigma^{-1})w \ge 0\}$, the inequality holding componentwise at every time. Is it true that, for all $R$ and $R'$,
\[
\cB(R) \subseteq \cB(R') \iff R' = HR \ \text{ for some } H \in \R_+^{g'\times g}[s,s^{-1}]
\]
with $\R_+[s,s^{-1}]$ the Laurent polynomials with nonnegative coefficients? The \enquote{if} direction is immediate, and the \enquote{only if} direction is proved in \cite{BM04} when $R$ has full row rank.
\end{problem}

\begin{problem}[7.3/P1, Jouan]\label{prob:jouan}
Let $X$ be a compact connected $C^\infty$ (resp.\ $C^\omega$) manifold and $\dot x = f(x)$, $y = h(x)\in\R^p$ an observable system of the same class. A $C^r$ observer is a system $\dot z = F(z,y)$, $\hat x = \theta(z)$ on an open $V\subseteq\R^N$ with $F$ of class $C^r$ and $\theta : V\to X$ continuous, such that $d(x(t),\theta(z(t)))\to0$ for every $x(0)\in X$ and $z(0)\in V$. Does the existence of a $C^r$ observer imply the existence of $n$ and a $C^r$ function $\varphi$ with
\[
L_f^nh = \varphi\bigl(h, L_fh,\dots,L_f^{n-1}h\bigr)\ ?
\]
\end{problem}

\begin{problem}[3.2, Devanathan]\label{prob:dev}
Let $\dot\xi = f(\xi) + g(\xi)\zeta$ be a single-input system on $\R^n$ with controllable linearization at the origin, brought by linear coordinates and linear feedback to $\dot x = Ax + Gu + \sum_{m\ge2}f_m + \sum_{m\ge1}g_mu$ with $A$ diagonal, nonresonant, and $G = (1,\dots,1)^T$. Find, for every $m\ge2$, a near-identity change of state $x = y + h_m(y)$ and a change of input $v = (1+\beta_{m-1})u + \alpha_m$ solving the generalized homological equations
\[
Dh_m\,Ay - Ah_m + G\alpha_m = f''_m,\qquad Dh_m\,G + G\beta_{m-1} = g''_{m-1},
\]
so that the nonlinear terms are removed up to arbitrary order. The case $m = 2$ was solved in the sense that the reduced linear system was shown to have $n(n-1)/2$ equations in $n$ unknowns and rank $n-1$.
\end{problem}

\begin{problem}[1.10, Trentelman]\label{prob:rfi-intro}
Let $\cP_{\rm full}$ be a linear differential behaviour in the variables $(w,c)$, $w\in\R^{\mathtt w}$ the manifest variable and $c\in\R^{\mathtt c}$ the interconnection variable, and let $\cC$ be a behaviour in $c$. The controlled behaviour is $\cK = \{w : \exists c,\ (w,c)\in\cP_{\rm full},\ c\in\cC\}$. Write $p(\cdot)$ for output cardinality and $n(\cdot)$ for McMillan degree. The interconnection is \emph{regular} if $p(\cP_{\rm full}\wedge_c\cC) = p(\cP_{\rm full})+p(\cC)$ and a \emph{regular feedback interconnection} if in addition $n(\cP_{\rm full}\wedge_c\cC) = n(\cP_{\rm full})+n(\cC)$. Characterize the $\cK$ for which a controller $\cC$ exists implementing $\cK$ by a regular feedback interconnection. Implementability and regular implementability had been characterized by Willems--Trentelman and Belur--Trentelman; the McMillan-degree condition was the open part of the 2004 problem.
\end{problem}

\section{The top-down route from diagrams to results}

The common diagram isolates a \emph{regular part}, where the requested witness is obtained by inverting a linear map, and a \emph{degenerate set}, where that map fails. Applying this distinction before attempting a proof gives the following four targets.

In Problem~\ref{prob:farkas} the regular part is $R$ of full row rank, where the shift operator $R(\sigma)$ is onto and every implied inequality has a unique certificate. The degenerate set is where $R$ has syzygies, and there the behaviour $\cB(R)$ can satisfy a linear \emph{equation} that is true only because of what trajectories do as $t\to\pm\infty$. A certificate $H$ is a finite computation at finitely many time instants and cannot detect that. Chapter~\ref{ch:fa} exhibits the simplest instance, $R = [\,1\,;\,-(1-s)^2\,]$, whose behaviour is the nonnegative constants, and proves that the conjecture holds for a given $R$ exactly when a certain cone of certificates is closed, which fails whenever such a hidden equality exists without a nonnegative syzygy to certify it.

In Problem~\ref{prob:jouan} the regular part is where the observability map $\Phi_n = (h, L_fh,\dots,L_f^{n-1}h)$ is an injective immersion for some $n$: then $\varphi = L_f^nh\circ\Phi_n^{-1}$ is as smooth as the data and the question is empty. The degenerate set is $\bigcap_n\{\rank D\Phi_n < \dim X\}$. In the smooth category a flat output puts a point there, and Chapter~\ref{ch:ob} shows that at such a point a smooth observer still exists, global in Jouan's sense: a Kazantzis--Kravaris/Luenberger filter whose immersion is an exponential average of $h$ along the flow (averaging removes the flat point), followed by an angle tracker on a circle factor of the observer state space, which is needed because no continuous estimate map on a contractible domain can invert the immersion on the circle. Meanwhile no relation $L_f^nh = \varphi\circ\Phi_n$ is even Lipschitz.

In Problem~\ref{prob:dev} the regular part is the invertibility of $\ad_{Ax}$ on homogeneous vector fields under nonresonance, which is Poincar\'e's theorem. The degenerate set is the cokernel of the augmented operator $[\ad_{Ax},\ad_G\mid G\alpha, G\beta]$, which the problem implicitly assumes to vanish. Chapter~\ref{ch:nf} computes it: it has dimension $(n-1)\binom{m+n-2}{n-1}-\binom{m+n-1}{n-1}+1$, positive for $n\ge3$, and the one-dimensional kernel of the operator is the infinitesimal symmetry of the linear system that reparametrizes its flat output.

In Problem~\ref{prob:rfi-intro} the regular part is regularity itself: the stacked kernel representation $[R_\cP;R_\cC]$ has full row rank, the equation modules of plant and controller are independent, and the output cardinalities add. The McMillan degree is not additive under this condition. It is the invariant of the \emph{associated graded} object, the module generated by the leading coefficient rows, and Chapter~\ref{ch:rfi} shows that regular feedback is exactly regularity after passing to the associated graded: the stacked leading coefficient matrix must have full row rank. The degenerate set is where leading coefficients cancel. In the full interconnection case this gives the complete answer, a rank condition on the leading matrix of the factor $F$ with $R_\cP = FR_\cK$, tested on finitely many thresholds, with the controller written down from a homogeneous completion of that matrix. In the partial case the same condition is applied to the canonical controller inside the interconnection-variable space, and one more identity is needed: the McMillan degree the plant carries in the eliminated variable must not be something the desired behaviour has to see.

The programme is therefore: draw the requested witness, isolate the map whose inversion supplies it on the regular part, and replace failed inversion on the degenerate set by its defect or universal substitute. In the four cases below this produces, respectively, an equality at infinity and a closure operation, a new filtered immersion, a cokernel, and an associated-graded exactness condition. The technical chapters compute those objects.

\section{What is machine-checked}

For Chapter~\ref{ch:fa}, Problem~\ref{prob:farkas} is stated in Lean~4 with Mathlib as a proposition over all Laurent polynomial matrices and its negation is proved; the same holds for the companion minimality conjecture. The general obstruction and the degree bound behind the closedness criterion are also formalized. The obstruction half of Chapter~\ref{ch:ob} is formalized in a one-dimensional model; of the observer half, the abstract convergence theorem (the circle lemma, the angle lift, and the two-stage argument, for every initial condition) is formalized, while the star-shapedness of the concrete immersion is not. The kernel theorem of Chapter~\ref{ch:nf}, the diagonal action of $\ad_{Ax}$, and the rank--nullity count are formalized; the binomial dimension formula is arithmetic on top of them. For Chapter~\ref{ch:rfi}, two layers are formalized: the module-theoretic layer of the partial interconnection case over an arbitrary commutative ring (implementation and regularity as lattice identities, the reduction to the interconnection-variable space, the recovery of the Belur--Trentelman theorem through torsion-freeness over a PID), and the filtered layer over a field, the lifting lemma \enquote{leading forms detect bases} with the full-interconnection theorem in the form \enquote{a regular feedback controller exists iff the plant's leading rows admit a homogeneous completion in $\operatorname{gr}E(\cK)$}.

\chapter[The Diagrammatic Programme]{The Diagrammatic Programme\\[4pt]\large From universal shapes to four control problems}\label{ch:diag}

\section*{Summary}
The companion notes on control architectures as commutative diagrams propose a principle: when the arrow a design asks for does not exist, study its universal substitute, and read the design's properties off the universal property. This chapter applies that principle before the four technical investigations. In each problem there is a functor from a syntactic category (representations: polynomial matrices, certificates, jets, Taylor data) to a semantic one (behaviours, trajectories, geometry), and the question is whether a semantic fact is reflected syntactically. On the \emph{regular part} it is, because a linear map is invertible there; the \emph{degenerate set} is where the map fails. The four failures are of four standard kinds: non-fullness of a functor (1.9), non-existence of a lift in a smooth category (7.3), a nonzero cokernel (3.2), and non-strictness of a filtered exact sequence (1.10). The framework gives the search target in all four cases. For 1.9 the closedness criterion and the implicit-equality obstruction follow from a Galois connection; only the identification of the Galois closure with the weak closure is analysis. For 1.10 the McMillan degree becomes a Hilbert-function invariant of the associated graded module, regular feedback becomes strict exactness, and the partial case is forced to carry the defect of elimination. For 3.2 the obstruction count is the augmentation exact sequence plus one number. For 7.3 the framework identifies a missing lift, prescribes enlargement by a different immersion, and points to the flat degenerate set, while analysis must construct the example. The rest of the volume follows these four routes.

\section{The framework, recalled}

The companion notes \cite{Bor26diag} define a category $\mathsf{Beh}$ whose objects are admissible signal sets and whose morphisms are well-posed relations, with a monoidal product for parallel placement, a trace for feedback, a dagger for inversion, an order for uncertainty and a gain enrichment for performance. Eleven control architectures are drawn as diagrams asked to commute; four shapes emerge (a section of an epimorphism, a fixed point of a trace, a factorization through a subobject, a section of a sheaf), each with a precise existence statement. The main structural result is an exact sequence for the augmented map $A = [\Lambda\mid B]$ of the continuation controller,
\begin{equation}\label{eq:augseq}
0\to\ker\Lambda\to\ker A\to\R^p\xrightarrow{\ q\circ B\ }\coker\Lambda\to\coker A\to0,
\end{equation}
whose readings are the repair criterion, the authority collapse, and the residual-authority count; it is machine-checked in \texttt{ExactSequence.lean}. The notes hold the formalism to a standard: a categorical statement earns its place if it predicts something not already in hand or connects the construction to a known result elsewhere; otherwise it is an encoding.

\section{The thesis: syntax, semantics, and the degenerate set}

Each of the four problems has the following shape. There is a category $\mathsf{Syn}$ of representations, a category $\mathsf{Sem}$ of the objects represented, and a functor $\mathcal B : \mathsf{Syn}\to\mathsf{Sem}$. A property holds in $\mathsf{Sem}$, and the problem asks whether it is witnessed in $\mathsf{Syn}$. Where a certain linear map is invertible (the regular part), the witness is obtained by inverting it. Where the map fails (the degenerate set), the witness may not exist, and the problem is to say exactly when.

\begin{center}\small
\begin{tabular}{@{}>{\raggedright\arraybackslash}p{1.6cm}>{\raggedright\arraybackslash}p{3.1cm}>{\raggedright\arraybackslash}p{3.1cm}>{\raggedright\arraybackslash}p{3.3cm}>{\raggedright\arraybackslash}p{3.6cm}@{}}
\toprule
Problem & Syntax $\to$ semantics & Regular part & Degenerate set & Failure mode \\
\midrule
1.9 & nonnegative certificates $H$ $\to$ inclusions of behaviours & $R$ of full row rank: $R(\sigma)$ onto & syzygies; implicit equalities at $t\to\pm\infty$ & functor not full \\
7.3/P1 & output derivatives $\Phi_n$ $\to$ the map $L_f^nh$ & $\Phi_n$ injective immersion & $\bigcap_n\{\rank D\Phi_n<\dim X\}$ & no $C^1$ lift \\
3.2 & $(h_m,\alpha_m,\beta_{m-1})$ $\to$ degree-$m$ terms & $\ad_{Ax}$ invertible (nonresonance) & $\coker[\Lambda_m\mid B_m]$ & cokernel $\ne0$ \\
1.10 & controller rows $\to$ interconnection & $[R_\cP;R_\cC]$ full row rank & leading-coefficient dependence & graded sequence not exact \\
\bottomrule
\end{tabular}
\end{center}

The rest of the chapter makes each row precise.

\section{Problem 1.9: fullness of the behaviour functor}

Let $\mathsf{Rep}_+$ be the category whose objects are Laurent polynomial matrices $R\in\R^{g\times q}[s,s^{-1}]$ and whose morphisms $R\to R'$ are nonnegative Laurent matrices $H$ with $R' = HR$; let $\mathsf{Sub}$ be the poset of subsets of $(\R^q)^\Z$. The behaviour functor $\cB : \mathsf{Rep}_+^{\rm op}\to\mathsf{Sub}$ sends $R$ to $\cB(R)$ and a morphism $H$ to the inclusion $\cB(R)\subseteq\cB(HR)$ (Lemma 4.5.9 of the thesis \cite{tenDam97}, the \enquote{if} direction).

\begin{proposition}\label{prop:diag-full}
Conjecture~\ref{fa:conj:farkas} is the statement that $\cB$ is \emph{full}: every inclusion $\cB(R)\subseteq\cB(R')$ is the image of a morphism. Fullness holds on the regular part ($R$ of full row rank) and fails in general; it holds for a given $R$ iff the certificate cone $\cK(R)$ is weakly closed (Proposition~\ref{fa:prop:bipolar}), and fails whenever $\cB(R)$ has an implicit equality not generated by a nonnegative syzygy (Proposition~\ref{fa:prop:obstruction}).
\end{proposition}

The contrast with the linear case is the point. For linear equality behaviours over $\R^\Z$ or $C^\infty(\R)$, the functor $R\mapsto\ker R(\sigma)$ from polynomial matrices with morphisms $R' = HR$ ($H$ arbitrary) is full: $\ker R\subseteq\ker R'$ iff $R' = HR$, for \emph{all} $R$, not only full row rank. This is the fundamental principle of Oberst's duality \cite{Obe90}, and its module-theoretic reason is that the signal space is an injective module over the operator ring, so that $\operatorname{Hom}$ into it turns the inclusion of row modules into a surjection. The positive structure breaks injectivity: the \enquote{module} of nonnegative certificates is a cone, inclusions of behaviours are detected only by the closure of the cone, and the closure is the universal substitute. In the language of the companion notes, the arrow $H$ is missing, but $\cK(R)^{\circ\circ} = \cl\cK(R)$, the double polar, exists and is the best approximation: it is the set of all implied rows, and Theorem~\ref{fa:thm:main} is a row in $\cl\cK(R)\setminus\cK(R)$.

What the structure derived here: the \emph{shape} of the answer, that the conjecture is equivalent to closedness of a cone and that closedness can be tested on finite-dimensional slices (Lemma~\ref{fa:lem:finite}); and the necessary condition through implicit equalities, which is a statement about morphisms in $\mathsf{Rep}_+$ (nonnegative syzygies). What analysis supplied: the counterexample itself, whose proof is the behaviour of concave sequences at infinity; and the same-sign-ends sufficient condition, which is a degree bound.

\subsection*{The closure is the monad of a Galois connection}

The criterion of Proposition~\ref{fa:prop:bipolar} and the obstruction of Proposition~\ref{fa:prop:obstruction} are formal consequences of an adjunction; only the identification of the abstract closure with the weak closure is analysis. Let $T = (\R^q)^\Z$, $E^* = \R^{1\times q}[s,s^{-1}]$, and write $r\vdash w$ if $(r(\sigma)w)(t)\ge0$ for all $t\in\Z$. The relation $\vdash$ induces an antitone Galois connection between subsets of $E^*$ and subsets of $T$,
\[
S^\triangleright = \{w : r\vdash w\ \ \forall r\in S\},\qquad W^\triangleleft = \{r : r\vdash w\ \ \forall w\in W\},\qquad S\subseteq W^\triangleleft\iff W\subseteq S^\triangleright ,
\]
so $S\mapsto S^{\triangleright\triangleleft}$ is a closure operator and $S^{\triangleright\triangleleft\triangleright} = S^\triangleright$. By definition $\{\text{rows of }R\}^\triangleright = \cB(R)$. Since $r\vdash w$ implies $s^kr\vdash w$ (the quantifier over $t$ absorbs the shift) and $\vdash$ is preserved by nonnegative combinations, $\cK(R)\subseteq\{\text{rows of }R\}^{\triangleright\triangleleft}$, hence $\cK(R)^\triangleright = \cB(R)$.

\begin{proposition}\label{prop:diag-galois}
$\cB(R)\subseteq\cB(R')$ iff every row of $R'$ lies in $\cK(R)^{\triangleright\triangleleft}$. Consequently Conjecture~\ref{fa:conj:farkas} holds for a given $R$ iff $\cK(R)$ is closed for the Galois closure, $\cK(R)^{\triangleright\triangleleft} = \cK(R)$; and the functor $\cB$ of Proposition~\ref{prop:diag-full} is full iff the monad $S\mapsto S^{\triangleright\triangleleft}$ is trivial on every $\cK(R)$. Moreover $\cK(R)^{\triangleright\triangleleft}\cap(-\cK(R)^{\triangleright\triangleleft})$ is the space of implicit equalities of $\cB(R)$, while $\cK(R)\cap(-\cK(R))$ is spanned by the rows $H_1R$ with $(H_1+H_2)R = 0$ for some nonnegative $H_2$; so the conjecture forces every implicit equality to come from a nonnegative syzygy (Proposition~\ref{fa:prop:obstruction}).
\end{proposition}
\begin{proof}
$\cB(R)\subseteq\{\text{rows of }R'\}^\triangleright$ iff $\{\text{rows of }R'\}\subseteq\cB(R)^\triangleleft = \cK(R)^{\triangleright\triangleleft}$. The conjecture asserts that the rows of $R'$ then lie in $\cK(R)$; letting $R'$ run over single rows gives the equivalence. An implicit equality is a $p$ with $p\vdash w$ and $-p\vdash w$ for all $w\in\cB(R)$, i.e.\ $\pm p\in\cB(R)^\triangleleft$. If $p = H_1R$ and $-p = H_2R$ with $H_1,H_2\ge0$ then $(H_1+H_2)R = 0$.
\end{proof}

The analysis enters once: $\cK(R)$ and $\cB(R)$ are shift-invariant convex cones, so $S^\triangleright$ and $W^\triangleleft$ are polars for the bilinear pairing $\langle r,w\rangle = (r(\sigma)w)(0)$, and $\cK(R)^{\triangleright\triangleleft} = \cK(R)^{\circ\circ}$ is the $\sigma(E^*,E)$-closure by the bipolar theorem, once $E^*$ and $T$ are placed in duality (this is where Kakutani--Klee is used). Everything in Proposition~\ref{prop:diag-galois} holds for any relation $\vdash$ preserved by the operations that generate $\cK(R)$; the Hahn--Banach step is what makes the closure computable. What the structure cannot do is decide for which $R$ the cone is closed; that is the counterexample and the same-sign-ends bound.

\section{Problem 7.3/P1: a lift that exists continuously but not smoothly}

For the system $\dot x = f(x)$, $y = h(x)$, the relation $L_f^nh = \varphi\circ\Phi_n$ is a lifting problem:
\[
\begin{array}{ccc}
X & \xrightarrow{\ L_f^nh\ } & \R^p \\
{\scriptstyle\Phi_n}\downarrow & \nearrow_{\varphi} & \\
\R^{pn} & &
\end{array}
\]
asked in the category $C^r$. On the regular part, where $\Phi_n$ is an injective immersion, the lift exists in $C^r$ for the same $r$ as the data: $\varphi = L_f^nh\circ\Phi_n^{-1}$ on the embedded image, extended by a tubular retraction. In $C^0$ the lift exists whenever $\Phi_n$ is injective (Tietze). Jouan's question is whether the existence of an observer forces the $C^r$ lift at the degenerate set $\bigcap_n\{\rank D\Phi_n<\dim X\}$. Theorem~\ref{ob:thm:main} says it does not: the obstruction is a Lipschitz bound that fails at a flat point.

The universal-substitute reading is instructive because the substitute is not a lift. The KKL observer replaces the object: instead of asking for $\varphi$ on $\Phi_n(X)$, it builds a different immersion $T = \lambda\int_{-\infty}^0e^{\lambda\tau}h\circ\phi_\tau\,d\tau$ into the filter state space and takes a section of it. The companion notes list three classical instances of \enquote{enlarging the object} to repair a missing arrow: dynamic extension (derivatives of the input), preview (the time axis), and the continuation clock (the state). The KKL immersion is a fourth: enlargement by a filter state $z$, with the arrow $x\mapsto T(x)$ total and injective where $\Phi_n$ is injective but not immersive. The observer is then a section shape in $\mathsf{Beh}$: a left inverse of $T$ on its image. Here the shape has a global constraint that a local reading would miss: the section must be defined on a forward-invariant open set and be continuous there, and since it restricts to a homeomorphism on the embedded circle $T(S^1)$, that set cannot be contractible. A nearest-point retraction on a tube fails (the tube is not invariant), and so does any $\theta$ on $\R^N$; the paper realizes the section dynamically, by an angle tracker on the circle factor of $\R^2\times(\R^2\setminus\{0\})$, using that $T(S^1)$ is star-shaped. The regularity of this section is the whole content of Jouan's question, and the answer is that it can be smooth while every lift $\varphi$ is non-Lipschitz, because convolution along the flow is a smoothing operation and differentiation along the flow is not.

What the structure derived: that the question is empty off the degenerate set, so that a counterexample must have $\bigcap_n\cD_n\ne\emptyset$, and that in the smooth category such points are produced by flat functions; and the identification of KKL as an enlargement of the object, which predicted that its regularity is independent of the rank of $D\Phi_n$. What analysis supplied: the star-shapedness estimate for $T_\lambda$ (Lemma~\ref{ob:lem:embedding}), the one-sided attractor $\dot\epsilon = -K(1-\cos\epsilon)$ on the circle, and the asymptotics at the flat point (Lemma~\ref{ob:lem:asymp}). What topology supplied: the non-contractibility of the observer's domain.

\section{Problem 3.2: the augmentation sequence, literally}

The degree-$m$ homological equations are $L_m(h,\alpha,\beta) = (f''_m,g''_{m-1})$ with $L_m = [\Lambda_m\mid B_m]$, $\Lambda_m = (\ad_{Ax},\ad_G)$, $B_m(\alpha,\beta) = (G\alpha,G\beta)$ (Chapter~\ref{ch:nf}). This is the augmented map of the companion notes with $\Lambda = \Lambda_m$ and $B = B_m$, and the exact sequence \eqref{eq:augseq} reads
\[
0\to\ker\Lambda_m\to\ker L_m\to P_m\times P_{m-1}\xrightarrow{\ q\circ B_m\ }\coker\Lambda_m\to\coker L_m\to0 .
\]
Under nonresonance $\ad_{Ax}$ is injective on $\cH_m$, so $\ker\Lambda_m = 0$ and $\ker L_m\cong\ker(q\circ B_m)$: the kernel consists of the input changes $(\alpha,\beta)$ whose effect is absorbed by a state change. The three readings of the sequence become:
\begin{enumerate}
\item[(i)] \emph{Repair criterion.} All degree-$m$ terms are removable iff $q\circ B_m$ is onto, iff the images of $G\alpha$, $G\beta$ span $\coker\Lambda_m$ (Proposition~\ref{nf:prop:coker}). This fails for $n\ge3$.
\item[(ii)] \emph{Authority collapse.} Does not occur: $\ker\Lambda_m = 0$ means the state change is determined by the input change, and conversely the input change is pinned only up to $\ker(q\circ B_m)$.
\item[(iii)] \emph{Residual authority.} $\dim\ker L_m = \dim\ker(q\circ B_m) = 1$, the prolonged flat-output reparametrization (Proposition~\ref{nf:prop:kernel}). The sequence then gives $\dim\coker L_m = \dim\coker\Lambda_m - (\dim P_m + \dim P_{m-1}) + 1$, which is Corollary~\ref{nf:cor:formula}.
\end{enumerate}
Here the structure derived the \emph{count}: the formula for the obstruction space is the exact sequence plus the one number $\dim\ker L_m = 1$, and the two Lean modules \texttt{ExactSequence.lean} and \texttt{NormalForm.lean} together make this a formally composable statement (the first proves the sequence for any $\Lambda, B$; the second proves $\dim\ker = 1$ for this $\Lambda, B$). What analysis supplied is small and combinatorial: the variable bookkeeping showing that the kernel is generated by $y_1^m$. The universal substitute is the normal form, the canonical representative of $\coker L_m$, and the comparison with Kang's count is the connection to a known result that the notes' standard asks for.

\section{Problem 1.10: filtered modules and a second proof}

For $\cP,\cC\in\mathcal L^{\mathtt w}$ (full interconnection) consider the square of behaviours
\[
0\to\cP\cap\cC\to\cP\oplus\cC\xrightarrow{\ (a,b)\mapsto a-b\ }\cP+\cC\to0 .
\]
Under the duality between behaviours and finitely generated modules over $\R[\xi]$ \cite{Obe90}, with $E(\cB) = \operatorname{rowspan}R_\cB\subseteq\R^{1\times\mathtt w}[\xi]$ the module of equations, $\cP\cap\cC$ corresponds to $E(\cP)+E(\cC)$ and $\cP+\cC$ to $E(\cP)\cap E(\cC)$. The output cardinality is the rank of the equation module and is additive on exact sequences, so $p(\cP\cap\cC) + p(\cP+\cC) = p(\cP)+p(\cC)$, and the interconnection is regular iff $p(\cP+\cC) = 0$, iff $\cP+\cC$ is everything, iff $E(\cP)\cap E(\cC) = 0$: \emph{regularity is the exactness of $0\to E(\cP)\oplus E(\cC)\to\R^{1\times\mathtt w}[\xi]$}, the square is bicartesian. This is Proposition~\ref{rfi:prop:rowred}(i) in diagrammatic form.

The McMillan degree is not additive on exact sequences of behaviours (Example~\ref{rfi:ex:rfi1}). It is the invariant of a \emph{filtered} object, and the rest of this section shows that, once this is taken seriously, Theorems~\ref{rfi:thm:rfi-full} and~\ref{rfi:thm:rfi-partial} follow from two standard lemmas about filtered modules, with Rado's theorem and the determinant argument of Chapter~\ref{ch:rfi} no longer needed. This is the one problem of the four where the diagrammatic reading produces the result and not only its shape.

\subsection*{The filtered category}

Let $V = \R^{1\times\mathtt v}[\xi]$ with $F_eV$ the rows of degree $\le e$. Every submodule $E\subseteq V$ carries the induced filtration $F_eE = E\cap F_eV$, and $\operatorname{gr}E = \bigoplus_e F_eE/F_{e-1}E$ embeds in $\operatorname{gr}V = V$ as the graded submodule generated by the leading forms $\operatorname{lf}(x)$, $x\in E$. Morphisms are filtered maps; a morphism is \emph{strict} if $f(F_eA) = f(A)\cap F_eB$. The functor $\operatorname{gr}$ is exact on strict short exact sequences and only on those. Two facts, both standard:

\begin{enumerate}
\item[(F1)] (\emph{Leading forms detect bases.}) If $x_1,\dots,x_r\in E$ and $\operatorname{lf}(x_1),\dots,\operatorname{lf}(x_r)$ is a homogeneous basis of $\operatorname{gr}E$, then $x_1,\dots,x_r$ is a basis of $E$. If $x_1,\dots,x_r$ is a basis of $E$ whose leading forms are $\R[\xi]$-independent, then $\operatorname{gr}E$ is generated by them (predictable degree). Row-reduced representations are exactly bases of the second kind.
\item[(F2)] (\emph{Hilbert function.}) $\operatorname{gr}E$ is graded free; if its generators have degrees $k_1,\dots,k_p$ then $h_E(e) := \dim_\R\operatorname{gr}_eE = \#\{i : k_i\le e\}$ and $\sum_ik_i = \sum_{e\ge0}\bigl(p - h_E(e)\bigr)$. By Lemma~\ref{rfi:lem:rfi-facts}(d) this number is the McMillan degree: $n(\cB) = \sum_{e\ge0}\bigl(p(\cB) - h_{E(\cB)}(e)\bigr)$.
\end{enumerate}
Both halves of (F1) are formalized (\texttt{LeadingForms.lean}: \texttt{span\us eq\us of\us grDeg\us le}, \texttt{grDeg\us span\us eq}, Chapter~\ref{ch:formal}). Proof of (F1): for $x\in F_eE$ write $\operatorname{lf}(x) = \sum a_i\operatorname{lf}(x_i)$ with $a_i$ homogeneous; then $x - \sum a_ix_i\in F_{e-1}E$, and induction on $e$ (the filtration is bounded below) gives generation; a relation $\sum a_ix_i = 0$ would make the top-degree parts of the $a_i\operatorname{lf}(x_i)$ cancel, contradicting independence. For the second statement, if the leading forms are independent no cancellation of leading terms is possible in $\sum a_ix_i$, so $\operatorname{lf}(\sum a_ix_i)\in\operatorname{span}\{\operatorname{lf}x_i\}$.

So $n$ is a function of $\operatorname{gr}E(\cB)$ alone: \emph{the McMillan degree is a Hilbert-function invariant}, and the row-degree multiset $\kappa(\cB)$ (Corollary~\ref{rfi:cor:rowdeg-inv}) is the multiset of generator degrees of $\operatorname{gr}E(\cB)$. This is the categorical content of Forney's minimal bases.

\begin{lemma}[Strictness]\label{lem:diag-strict}
Let $E_1,E_2\subseteq V$ with $E_1\cap E_2 = 0$. Then $n(E_1\oplus E_2)\le n(E_1)+n(E_2)$, with equality iff $\operatorname{gr}E_1\cap\operatorname{gr}E_2 = 0$, iff $\operatorname{gr}(E_1\oplus E_2) = \operatorname{gr}E_1\oplus\operatorname{gr}E_2$, iff the inclusions $E_i\to E_1\oplus E_2$ form a strict direct sum.
\end{lemma}
\begin{proof}
Take bases of $E_1$, $E_2$ with independent leading forms (F1). If $\operatorname{gr}E_1\cap\operatorname{gr}E_2 = 0$ the union has independent leading forms, so by (F1) it generates $\operatorname{gr}(E_1\oplus E_2) = \operatorname{gr}E_1\oplus\operatorname{gr}E_2$ and $n$ adds by (F2). If $\operatorname{gr}E_1\cap\operatorname{gr}E_2\ne0$, some constant combination of the leading rows of the union vanishes; replacing the row of highest degree in that combination by the combination (times a power of $\xi$) lowers the degree sum and keeps a basis of $E_1\oplus E_2$; iterating until the leading rows are independent reaches a basis with degree sum $n(E_1\oplus E_2)$ (F1, F2), strictly below $n(E_1)+n(E_2)$. Finally $\operatorname{gr}E_1\cap\operatorname{gr}E_2 = 0$ iff no $u\in E_1$, $v\in E_2$ have $\deg(u+v)<\deg u = \deg v$, which is strictness of the sum.
\end{proof}

\begin{proposition}\label{prop:diag-strict}
The interconnection $\cP\wedge\cC$ is regular feedback iff $0\to E(\cP)\oplus E(\cC)\to V$ is exact \emph{and strict}, i.e.\ iff $0\to\operatorname{gr}E(\cP)\oplus\operatorname{gr}E(\cC)\to\operatorname{gr}V$ is exact. Regularity is exactness in $\mathsf{Mod}_{\R[\xi]}$; regular feedback is exactness after $\operatorname{gr}$.
\end{proposition}
\begin{proof}
Regular feedback is regularity plus $n(\cP\cap\cC) = n(\cP)+n(\cC)$, i.e.\ $E(\cP)\cap E(\cC) = 0$ plus the equality case of Lemma~\ref{lem:diag-strict}.
\end{proof}

\subsection*{The full case}

Fix $\cK\subseteq\cP$, so $E(\cP)\subseteq E(\cK)$, both free.

\begin{theorem}[Theorem~\ref{rfi:thm:rfi-full}, second proof]\label{thm:diag-full}
$\cK$ is implementable by regular feedback from $\cP$ iff $\operatorname{gr}E(\cP)$ is a graded direct summand of $\operatorname{gr}E(\cK)$, iff $\operatorname{gr}E(\cK)/\operatorname{gr}E(\cP)$ is torsion-free, iff condition (d) of Theorem~\ref{rfi:thm:rfi-full} holds.
\end{theorem}
\begin{proof}
A regular feedback controller is $E_C$ with $E(\cP)\oplus E_C = E(\cK)$ and, by Proposition~\ref{prop:diag-strict}, $\operatorname{gr}E(\cP)\oplus\operatorname{gr}E_C = \operatorname{gr}E(\cK)$; so $\operatorname{gr}E(\cP)$ is a graded summand. Conversely let $\operatorname{gr}E(\cK) = \operatorname{gr}E(\cP)\oplus H$ with $H$ graded; $H$ is graded free (a graded submodule of a free module over $\R[\xi]$). Lift a homogeneous basis of $H$ to elements $x_l\in E(\cK)$ with $\operatorname{lf}(x_l)$ that basis, and let $E_C := \operatorname{span}\{x_l\}$. Together with a basis of $E(\cP)$ with independent leading forms, the $x_l$ have leading forms forming a basis of $\operatorname{gr}E(\cK)$, so by (F1) they form a basis of $E(\cK)$: $E(\cP)\oplus E_C = E(\cK)$, and $\operatorname{gr}E_C = H$ by the second half of (F1). By Proposition~\ref{prop:diag-strict} the interconnection is regular feedback. This proves the first equivalence.

For the second: a finitely generated graded module over $\R[\xi]$ is free iff torsion-free, and a graded submodule $N$ of a graded free $F$ is a graded summand iff $F/N$ is graded free (lift a homogeneous basis of $F/N$). For the third, write everything in the homogeneous basis $\ell_j\xi^{k_j}$ of $\operatorname{gr}E(\cK)$: $\operatorname{gr}_eE(\cK)\cong\R^{J(e)}$ with $J(e) = \{j : k_j\le e\}$, multiplication by $\xi^m$ is the inclusion $\R^{J(e)}\subseteq\R^{J(e+m)}$, and $\operatorname{gr}_eE(\cP)\cong\operatorname{span}\{L_i : d_i\le e\}$ where $L_i\xi^{d_i}$ are the leading forms of the rows of $R_\cP = FR_\cK$, i.e.\ $L = L_k(F)$. A torsion element of the quotient is an $x\in\R^{J(e)}\cap\operatorname{span}\{L_i : d_i\le e+m\}$ not in $\operatorname{span}\{L_i : d_i\le e\}$; torsion-freeness for all $e,m$ is, taking $m$ large, $\R^{J(e)}\cap\operatorname{rowspan}L = \operatorname{span}\{L_i : d_i\le e\}$ for all $e$, which is (d).
\end{proof}

The equivalence (a)$\Leftrightarrow$(graded summand with a homogeneous complement) is machine-checked as \texttt{exists\us controller\us iff}, with the lifting done by \texttt{span\us eq\us of\us grDeg\us le}; the passage to torsion-freeness and to (d) is on paper. Compared with the proof in Chapter~\ref{ch:rfi}: Rado's theorem is replaced by \enquote{torsion-free implies free} for graded $\R[\xi]$-modules, and the determinant argument showing that a $k$-row-reduced completion is unimodular is replaced by (F1). Both replacements are instances of general facts about filtered modules whose associated graded is free; nothing specific to polynomial matrices remains except the computation of the Hilbert functions in the last paragraph, which is where the explicit criterion (d) appears. The controller is the same: the lifts $x_l$ are the homogeneous rows $G_lR_\cK$ of Theorem~\ref{rfi:thm:rfi-full}.

\subsection*{The partial case as a base change}

Let $E_P = E(\cP_{\rm full})\subseteq\R^{1\times\mathtt w}[\xi]\oplus\R^{1\times\mathtt c}[\xi]$ with the two projections $\pi_w$, $\pi_c$ restricted to $E_P$, and let $E_K = E(\cK)$. The modules of Chapter~\ref{ch:rfi} are
\begin{gather*}
0\oplus E^0 = \ker(\pi_w|_{E_P}),\qquad E(\cP)\oplus0 = \ker(\pi_c|_{E_P}),\qquad M = \pi_c(E_P),\\
E_P^K = \pi_w^{-1}(E_K)\cap E_P,\qquad E_{\rm can} = \pi_c(E_P^K),
\end{gather*}
so $E_P^K$ is the pullback of $E_K\hookrightarrow\R^{1\times\mathtt w}[\xi]$ along $\pi_w$, and $E_{\rm can}$ is its image under $\pi_c$: the canonical controller is \emph{the base change of $\cK$ along the plant}, from the $w$-space to the $c$-space. A controller $E_C\subseteq\R^{1\times\mathtt c}[\xi]$ enters through $\pi_c^{-1}(E_C)\subseteq E_P$, and Lemma~\ref{rfi:lem:partial-modules}(i) says that $\cC$ implements $\cK$ iff the two pullbacks agree,
\[
\pi_c^{-1}(E_C + E^0) = \pi_w^{-1}(E_K)\cap E_P ,
\]
where $E^0$ appears because $\pi_c^{-1}$ cannot distinguish $E_C$ from $E_C+E^0$ modulo $\ker\pi_w$. Everything in the problem is thereby transported to the $c$-space, where it becomes the full problem for the pair $E^0\subseteq E_{\rm can}$ (Lemma~\ref{rfi:lem:partial-modules}(i),(ii) with $E_C\subseteq M$), and Theorem~\ref{thm:diag-full} applies, giving condition (A).

What does not transport is the filtration. The functor $\operatorname{gr}$ does not commute with $\ker\pi_w$:
\[
\operatorname{gr}(\ker\pi_w|_{E_P}) = 0\oplus\operatorname{gr}E^0\ \subseteq\ \ker(\operatorname{gr}\pi_w|_{\operatorname{gr}E_P}) = 0\oplus G ,
\]
with equality iff $\pi_w|_{E_P}$ is strict. The quotient $G/\operatorname{gr}E^0$ is the defect of exactness of $\operatorname{gr}$ on the elimination sequence $0\to E^0\to E_P\to\pi_wE_P\to0$. Regular feedback in the $(w,c)$-space is $\operatorname{gr}E_P\cap(0\oplus\operatorname{gr}E_C) = 0$ (Proposition~\ref{prop:diag-strict}), i.e.\ $G\cap\operatorname{gr}E_C = 0$, whereas regular feedback in the $c$-space is only $\operatorname{gr}E^0\cap\operatorname{gr}E_C = 0$. The difference is exactly the defect, and condition (B$'$) of Theorem~\ref{rfi:thm:rfi-partial}, $G\cap\operatorname{gr}E_{\rm can} = \operatorname{gr}E^0$, says that the defect is invisible from $\operatorname{gr}E_{\rm can}$. The proof of the theorem is then the modular law in the lattice of graded submodules, and the McMillan-degree form (B) is Lemma~\ref{lem:diag-strict} applied to the direct sum $E_P\oplus(0\oplus E_C)$.

In this reading the partial theorem is forced: the problem is a pullback, $\operatorname{gr}$ preserves pullbacks only of strict maps, and the obstruction to transporting the full criterion is the measure of non-strictness of the one map in the diagram that is not strict, the elimination of $w$. Example~\ref{rfi:ex:partial-neg} is the case where the whole McMillan degree of the plant lives in the eliminated variable: $E^0 = 0$ while $G = (\xi)$, the defect is everything, and nothing is implementable by regular feedback.

\section{The universal substitutes}

\begin{center}\small
\begin{tabular}{@{}>{\raggedright\arraybackslash}p{1.6cm}>{\raggedright\arraybackslash}p{4.2cm}>{\raggedright\arraybackslash}p{4.6cm}>{\raggedright\arraybackslash}p{4.4cm}@{}}
\toprule
Problem & Missing arrow & Universal substitute & Status \\
\midrule
1.9 & certificate $H$ for an implied row & the closure $\cl\cK(R) = \cK(R)^{\circ\circ}$ & substitute exists; the arrow does not (Thm.~\ref{fa:thm:main}) \\
7.3/P1 & $C^1$ lift $\varphi$ of $L_f^nh$ through $\Phi_n$ & enlarge the object: KKL immersion $T$ and its section & substitute exists and is smooth; the arrow does not (Thm.~\ref{ob:thm:main}) \\
3.2 & solution $(h_m,\alpha_m,\beta_{m-1})$ for all data & the normal form, canonical representative of $\coker L_m$ & cokernel computed (Cor.~\ref{nf:cor:formula}) \\
1.10 & regular-feedback controller & homogeneous complement of $\operatorname{gr}E(\cP)$ in $\operatorname{gr}E(\cK)$ & exists iff graded summand (Thm.~\ref{rfi:thm:rfi-full}); partial case: same in the $c$-space plus a degree identity (Thm.~\ref{rfi:thm:rfi-partial}) \\
\bottomrule
\end{tabular}
\end{center}

\section{Division of labour}

Continuing the ledger of the companion notes.

\paragraph{Derived from structure.}
\begin{itemize}
\item 1.9 is equivalent to fullness of a functor and to closedness of a cone for the Galois closure of $\vdash$ (Proposition~\ref{prop:diag-galois}); finite-dimensional slices suffice; implicit equalities without nonnegative syzygies refute it.
\item 7.3/P1 is empty off the degenerate set; a counterexample needs $\bigcap_n\cD_n\ne\emptyset$; KKL is an enlargement of the object, so its regularity is decoupled from $\operatorname{rank}D\Phi_n$.
\item 3.2 is an instance of the augmentation sequence; the obstruction count is the sequence plus $\dim\ker L_m$.
\item 1.10, in full: the McMillan degree is a Hilbert-function invariant of $\operatorname{gr}E(\cB)$; regular feedback is strict exactness (Lemma~\ref{lem:diag-strict}); the full-case criterion is \enquote{$\operatorname{gr}E(\cK)/\operatorname{gr}E(\cP)$ torsion-free}, which is (d) once written in coordinates (Theorem~\ref{thm:diag-full}); the partial case is the pullback of $\cK$ along the plant, and its only new obstruction is the defect $G/\operatorname{gr}E^0$ of exactness of $\operatorname{gr}$ on the elimination sequence.
\end{itemize}

\paragraph{Connected to known results.} Fullness of the linear behaviour functor $\leftrightarrow$ injectivity of the signal module (Oberst); the positive case $\leftrightarrow$ the bipolar theorem as the Galois closure of a convex cone. Non-stabilizing cones $\leftrightarrow$ the thesis's own counterexample for representability. KKL $\leftrightarrow$ Luenberger's original design and the Andrieu--Praly existence theorem. The kernel of $L_m$ $\leftrightarrow$ the symmetry group of the Brunovsky form (flat-output reparametrizations). (F1) and (F2) $\leftrightarrow$ Forney's minimal bases and the predictable degree property; the graded-complement criterion is the structural form of the Vinjamoor--Belur full-interconnection solution \cite{VB10}, and is equivalent to Rado's theorem in Chapter~\ref{ch:rfi} or to \enquote{torsion-free implies free} here; regular implementability $\leftrightarrow$ Belur--Trentelman via saturation (Proposition~\ref{rfi:prop:reg-impl}).

\paragraph{Supplied by analysis.} The identification of the Galois closure with the weak closure (Hahn--Banach, Kakutani--Klee) and the concavity-at-infinity argument (1.9). The star-shapedness estimate and the flat-point asymptotics (7.3). The variable bookkeeping under the drift derivation (3.2). For 1.10, in the second proof, nothing beyond the two filtered-module lemmas and the elimination theorem for behaviours.

\paragraph{Not reached.} Fullness for 1.9 is characterized but not decided for general $R$ (which $R$ have closed cones). The analytic case of 7.3. For 1.10, a characterization in terms of invariants of $\cK$ alone rather than of $\cC_{\rm can}$.

\section{From the programme to the four proofs}

The operative question is how much of each result is forced by the diagram before problem-specific calculation begins. The four problems give four different answers, and the difference determines the division of labour in the chapters that follow.

For 1.10 the answer is yes, completely. The semantics is a category of modules with a filtration, every object in the problem (behaviour, output cardinality, McMillan degree, interconnection, elimination) is a functorial construction on it, and the two theorems of Chapter~\ref{ch:rfi} are statements about when the functor $\operatorname{gr}$ preserves a splitting or a pullback. The proofs in this chapter use (F1), (F2), and the modular law, all of which are general facts, and they are shorter than the hands-on proofs because the degree bookkeeping is done once, in the definition of $\operatorname{gr}$, instead of once per lemma. The hands-on proofs remain useful because they are the ones that compute: the homogeneous completion is an algorithm, Rado's condition is a rank test.

For 1.9 the answer is yes for the criterion and no for the verdict. The Galois connection of $\vdash$ produces the closure, the equivalence of the conjecture with closedness, and the implicit-equality obstruction, with no analysis at all; it is an adjunction between posets and the statements are about its monad. But the connection cannot say which cones are closed, because closedness is a property of the specific pairing (sequences, shifts, nonnegativity at every time), and the counterexample is found by looking at what trajectories do at infinity. The reading tells you to look for an element of $\cK(R)^{\triangleright\triangleleft}\setminus\cK(R)$ and that such an element will be an implicit equality; it does not produce one.

For 3.2 the answer is yes for the count and trivially so: the problem is an instance of the exact sequence that the companion notes proved for other reasons, and the only input is $\dim\ker L_m = 1$. This is the case where the framework earned its place by prediction: the formula for the obstruction space was written down from the sequence before it was checked against Kang's normal form.

For 7.3 the answer is no. The categories involved ($C^r$ manifolds, $C^r$ maps) have no algebraic structure that a functor could preserve or fail to preserve in a computable way; the reading identifies the shape (a lift, an enlargement of the object, a section) and the place to look (a flat point), and then the analysis is the whole proof.

The pattern across the four is that the reading becomes technical exactly to the extent that the semantic category is algebraic. Modules over $\R[\xi]$: everything. Convex cones in duality: the criterion. Finite-dimensional linear maps: the count. Smooth manifolds: the shape. This is consistent with the standard the companion notes set, and it is also a guide to where the framework should be pushed next: toward the linear and behavioural problems of the collection, where the semantic categories are modules, cones, or filtered objects, and not toward the nonlinear geometric ones, where it will only ever draw the picture.

\chapter[A Counterexample to the Farkas Lemma for Behavioral Inequalities]{A Counterexample to the Farkas Lemma for Behavioral Inequalities\\[4pt]\large Problem 1.9 (ten Dam, Nieuwenhuis)}\label{ch:fa}

\section*{Summary}
ten Dam and Nieuwenhuis conjectured a Farkas lemma for linear difference inequalities on the time axis $\Z$: the behaviour $\{w : R(\sigma,\sigma^{-1})w \ge 0\}$ is contained in $\{w : R'(\sigma,\sigma^{-1})w \ge 0\}$ if and only if $R' = HR$ for a Laurent polynomial matrix $H$ with nonnegative coefficients. The conjecture appeared in ten Dam's 1997 thesis and as Problem~1.9 of \emph{Unsolved Problems in Mathematical Systems and Control Theory} (2004), where it was proved for $R$ of full row rank. We show that it is false in general: for $R = [\,1\,;\,-(1-s)^2\,]$ the behaviour consists of the nonnegative constant sequences, so $(1-\sigma)w \ge 0$ holds on it, but $1 - s$ is not a nonnegative combination of the rows of $R$. The mechanism is that the implication is true only through the behaviour at $t\to\pm\infty$, which no finite certificate sees. We identify the general obstruction (an implicit equality of the behaviour that is not generated by a nonnegative syzygy of $R$), show that the conjecture holds precisely when a certain cone of certificates is closed, give a sufficient condition in the scalar two-row case, and show that the same example refutes the companion conjecture on minimal representations. The counterexample is machine-checked in Lean~4/Mathlib.

\section{Introduction}

In the behavioural approach to systems theory \cite{Willems91}, a discrete-time linear system is a set of trajectories $w : \Z \to \R^q$ described by difference equations $R(\sigma,\sigma^{-1})w = 0$, where $R$ is a polynomial matrix in the shift $\sigma$ and its inverse. ten Dam \cite{tenDam93,tenDam97} began a parallel theory for \emph{behavioural inequalities}, systems described by $R(\sigma,\sigma^{-1})w \ge 0$, motivated by unilateral constraints in mechanics. A central question for such a theory is when one system of inequalities implies another, because elimination of latent variables and minimality of representations both reduce to it. In the static case ($R$ a constant matrix, no shifts) the answer is Farkas' lemma \cite{Farkas02}: $\{Rw \ge 0\} \subseteq \{R'w \ge 0\}$ iff $R' = HR$ with $H \ge 0$ entrywise. ten Dam and Nieuwenhuis conjectured the same statement for difference inequalities, with $H$ a Laurent polynomial matrix with nonnegative coefficients, and proved it for $R$ of full row rank \cite[Sec.~4.5]{tenDam97}, \cite{tDN04}. The problem was listed as Problem~1.9 in the Blondel--Megretski collection of open problems \cite{BM04}.

This note shows that the conjecture is false. The counterexample is two lines long once seen: the inequalities $w \ge 0$ and $-w(t) + 2w(t+1) - w(t+2) \ge 0$ describe the nonnegative midpoint-concave sequences on $\Z$, which are the nonnegative constants; hence $w(t) - w(t+1) \ge 0$ is implied; but $1 - s$ cannot be written as $h_1\cdot 1 - h_2(1-s)^2$ with nonnegative coefficients, as evaluation at $s=1$ shows. The point of the note is less the example than what it says about the dynamic Farkas problem: the implication $\cB(R)\subseteq\cB(R')$ can be true for reasons visible only at $t \to\pm\infty$, and a certificate $H$ is a finite computation at finitely many time instants. We make this precise in three ways. Section~\ref{fa:sec:obstruction} gives the general obstruction, in the vocabulary of ten Dam's thesis: an \emph{implicit equality} of the behaviour that is not generated by a nonnegative syzygy of $R$. Section~\ref{fa:sec:closed} shows that the conjecture holds for a given $R$ exactly when the cone of nonnegative combinations of the rows of $R$ is closed in a natural topology, and that closedness can be tested on finite-dimensional slices; this yields a positive result when the leading and trailing coefficients of the rows have equal signs. Section~\ref{fa:sec:minimal} shows that the same example refutes the minimality conjecture stated alongside the Farkas conjecture in \cite{tenDam97} and the related conjecture of \cite{tDN04}. Section~\ref{fa:sec:lean} reports a formalization in Lean~4.

\section{Problem formulation}\label{fa:sec:formulation}

\paragraph{Notation.} $\Rlp$ is the ring of real Laurent polynomials in $s$, and $\Rp$ the subset with nonnegative coefficients; matrices over these are $\R^{g\times q}[s,s^{-1}]$ and $\R_+^{g\times q}[s,s^{-1}]$. The time axis is $\Z$, the shift is $(\sigma w)(t) = w(t+1)$, and a Laurent polynomial matrix acts on $w \in (\R^q)^\Z$ by substitution, $(R(\sigma,\sigma^{-1})w)(t) = \sum_k R_k\,w(t+k)$ for $R = \sum_k R_k s^k$. We write $R(\sigma)$ for $R(\sigma,\sigma^{-1})$. Inequalities between vectors and sequences are componentwise at every time.

\begin{definition}
For $R \in \R^{g\times q}[s,s^{-1}]$, the \emph{inequality behaviour} is
\[
\cB(R) := \{ w \in (\R^q)^\Z : R(\sigma)w \ge 0 \}.
\]
\end{definition}

The trajectory space is all of $(\R^q)^\Z$; no boundedness or summability is imposed. This is the setting of \cite{tenDam97,tDN04}, where $(\R^q)^\Z$ carries the topology of pointwise convergence.

\begin{conjecture}[ten Dam--Nieuwenhuis; Problem 1.9 of \cite{BM04}]\label{fa:conj:farkas}
Let $R \in \R^{g\times q}[s,s^{-1}]$ and $R' \in \R^{g'\times q}[s,s^{-1}]$. Then $\cB(R) \subseteq \cB(R')$ if and only if there exists $H \in \R_+^{g'\times g}[s,s^{-1}]$ with $R' = HR$.
\end{conjecture}

The \enquote{if} direction is immediate: a Laurent polynomial matrix with nonnegative coefficients maps nonnegative sequences to nonnegative sequences, so $R(\sigma)w \ge 0$ gives $R'(\sigma)w = H(\sigma)R(\sigma)w \ge 0$. (Conversely, $H(\sigma)$ preserves nonnegativity of all sequences only if its coefficients are nonnegative \cite[Ch.~4]{tenDam97}.) The content is the \enquote{only if} direction.

\section{Literature}\label{fa:sec:literature}

The conjecture was first stated in ten Dam's ECC'93 paper \cite{tenDam93} and developed in Chapter~4 of his thesis \cite{tenDam97}, which introduces convex conical behaviours, shifted-polyhedral cones, and the elimination and minimality problems for difference inequalities. The thesis proves the conjecture when $R$ has full row rank \cite[Prop.~4.5.12]{tenDam97}; in that case $H$ is unique. The same material, with the conjecture and a second conjecture on minimal representations, appears in the Willems Festschrift \cite{tDNfest} and in Problem~1.9 of \cite{BM04}. The thesis records that the natural proof strategy, the bipolar theorem for the dual pair of $(\R^q)^\Z$ and the compactly supported sequences, stalls at the identification of the polar cone of $\cB(R)$ with the cone generated by the shifted rows of $R$; Section~\ref{fa:sec:closed} below shows that this identification is exactly what fails.

Two related results in the thesis deserve mention. First, a weaker conjecture (every closed, shift-invariant, finitely polyhedral convex cone is of the form $\cB(R)$) was disproved there by an increasing sequence of polyhedral cones whose union does not stabilize, with the remark that limits of polyhedral cones need not be polyhedral \cite[Ch.~4]{tenDam97}. This is the same non-stabilization phenomenon that drives our counterexample, identified for representability but not for implication. Second, the thesis defines \emph{implicit equalities}: a row $R_k$ of $R$ is an implicit equality of the representation $\cB(R)$ if $\cB(R) \subseteq \ker R_k(\sigma)$ \cite[Ch.~4]{tenDam97}. Our obstruction (Section~\ref{fa:sec:obstruction}) is a statement about implicit equalities.

Recent work on behavioural inequalities \cite{BPB25} studies feasibility of $H(\sigma)w \le g$ by a generalized Farkas lemma in the sense of theorems of the alternative; it does not address the implication problem and does not overlap with the present note. We are not aware of any published resolution of Conjecture~\ref{fa:conj:farkas}.

\section{The counterexample}\label{fa:sec:counterexample}

\begin{theorem}\label{fa:thm:main}
Let $q = 1$, 
\[
R = \begin{bmatrix} 1 \\ -(1-s)^2 \end{bmatrix} = \begin{bmatrix} 1 \\ -1 + 2s - s^2 \end{bmatrix}, \qquad R' = 1 - s .
\]
Then $\cB(R) \subseteq \cB(R')$, but there is no $H \in \R_+^{1\times 2}[s,s^{-1}]$ with $R' = HR$. Conjecture~\ref{fa:conj:farkas} is false.
\end{theorem}

\begin{proof}
\emph{Inclusion.} $\cB(R) = \{w \in \R^\Z : w \ge 0,\ -w(t) + 2w(t+1) - w(t+2) \ge 0\ \forall t\}$ is the set of nonnegative midpoint-concave sequences on $\Z$. Such a sequence is constant. Indeed, let $d(t) = w(t+1) - w(t)$; concavity says $d(t+1) \le d(t)$. If $d(t_0) < 0$ for some $t_0$, then $d(t) \le d(t_0)$ for all $t \ge t_0$, so $w(t_0 + n) \le w(t_0) + n\,d(t_0) \to -\infty$, contradicting $w \ge 0$. If $d(t_0) > 0$, then $d(t) \ge d(t_0)$ for all $t \le t_0$, so $w(t_0 - n) \le w(t_0) - n\,d(t_0) \to -\infty$, again a contradiction. Hence $d \equiv 0$ and $w$ is a nonnegative constant. For a constant, $(R'(\sigma)w)(t) = w(t) - w(t+1) = 0 \ge 0$, so $\cB(R) \subseteq \cB(R')$.

\emph{No certificate.} Suppose $1 - s = h_1\cdot 1 + h_2\cdot(-(1-s)^2)$ with $h_1, h_2 \in \Rp$. Evaluation at $s = 1$ is a ring homomorphism $\Rlp \to \R$. Applying it gives $0 = h_1(1) = \sum_k (h_1)_k$. The coefficients of $h_1$ are nonnegative, so $h_1 = 0$. Then $1 - s = -h_2(1-s)^2$, i.e.\ $(1-s)\bigl(1 + h_2(1-s)\bigr) = 0$. The ring $\Rlp$ is a domain and $1 - s \ne 0$, so $h_2(1-s) = -1$. Evaluating at $s = 1$ once more gives $0 = -1$.
\end{proof}

\begin{remark}[Reading]
The inclusion $\cB(R)\subseteq\cB(R')$ holds because of what nonnegative concave sequences do as $t \to \pm\infty$: any nonzero increment is amplified without bound in one time direction. A certificate $H$, applied to a trajectory, is a finite combination of finitely many shifts. No such combination can detect behaviour at infinity. Static Farkas has no such failure mode because a finite-dimensional system of inequalities has no \enquote{infinity}. The dynamic Farkas conjecture therefore fails for a structural reason, not an accidental one, and Section~\ref{fa:sec:obstruction} makes the structure explicit.
\end{remark}

\begin{remark}[Numerical illustration]\label{fa:rem:numerics}
Let $v_W(r')$ be the minimum of $\sum_t r'(t)w(t)$ over all $w$ on $[-W,W]$ with $|w|\le1$ that satisfy every constraint of $R$ whose support fits in $[-W,W]$. By finite Farkas, $v_W(r') \ge 0$ iff $r'$ is a nonnegative combination of the shifts of the rows of $R$ that fit in the window. For $r' = 1 - s$ the LP value is $v_W = -1/(2W+1)$: negative for every $W$, tending to $0$. The approximants $1 - s + \varepsilon$ do lie in the certificate cone, but the smallest support of a certificate for $1 - s + 2^{-j}$ is $2^j$. Both facts are the finite-dimensional shadow of Proposition~\ref{fa:prop:green} below.
\end{remark}

\section{The general obstruction}\label{fa:sec:obstruction}

\begin{definition}
A \emph{nonnegative syzygy} of $R \in \R^{g\times q}[s,s^{-1}]$ is $h \in \R_+^{1\times g}[s,s^{-1}]$ with $hR = 0$. An \emph{implicit equality} of $\cB(R)$ is a nonzero $p \in \R^{1\times q}[s,s^{-1}]$ with $\cB(R) \subseteq \ker p(\sigma)$.
\end{definition}

The second notion extends \cite[Ch.~4]{tenDam97}, which asks $p$ to be a row of $R$.

\begin{proposition}[Implicit equalities obstruct certificates]\label{fa:prop:obstruction}
Suppose $\cB(R)$ has an implicit equality $p$ and $R$ has no nonzero nonnegative syzygy. Then Conjecture~\ref{fa:conj:farkas} fails for $R$.
\end{proposition}
\begin{proof}
Both $p$ and $-p$ are implied by $R$. If the conjecture held for $R$, there would be $H_\pm \in \R_+^{1\times g}[s,s^{-1}]$ with $p = H_+R$ and $-p = H_-R$. Then $(H_+ + H_-)R = 0$ with $H_+ + H_- \ge 0$, and $H_+ + H_- \ne 0$ because $H_+R = p \ne 0$. This is a nonzero nonnegative syzygy.
\end{proof}

In Theorem~\ref{fa:thm:main}, $p = 1 - s$ is an implicit equality and the second row of $R$ is one in the thesis's sense; a nonnegative syzygy $h_1 = h_2(1-s)^2$ would force $h_1(1) = 0$, hence $h_1 = h_2 = 0$. In the static case every implicit equality is generated by a nonnegative syzygy (this is a consequence of Farkas' lemma itself), so the proposition has no static analogue. For $R$ of full row rank the operator $R(\sigma)$ is onto $(\R^g)^\Z$, so $\cB(R) = R(\sigma)^{-1}(\text{nonnegative sequences})$ satisfies no linear equation, and the proposition is silent, consistent with the positive result of \cite{tenDam97}.

\section{Closedness criterion and a positive result}\label{fa:sec:closed}

Let $E = (\R^q)^\Z$ with the pointwise topology and $E^* = \bigoplus_\Z\R^q$ the compactly supported sequences, paired by $\langle a, w\rangle = \sum_t a(t)^Tw(t)$. A Laurent row $r = \sum_k r_k s^k \in \R^{1\times q}[s,s^{-1}]$ is identified with the element of $E^*$ with value $r_k$ at time $k$; then $(r(\sigma)w)(t) = \langle s^t r, w\rangle$. Define the \emph{certificate cone}
\[
\cK(R) := \{ hR : h \in \R_+^{1\times g}[s,s^{-1}] \} = \cone\{ s^k R_i : k \in\Z,\ i = 1,\dots,g\} \subseteq E^*.
\]

\begin{proposition}\label{fa:prop:bipolar}
For a single row $r'$, $\cB(R) \subseteq \cB(r')$ iff $r' \in \cK(R)^{\circ\circ} = \cl_{\sigma(E^*,E)}\cK(R)$. Hence Conjecture~\ref{fa:conj:farkas} holds for $R$ iff $\cK(R)$ is weakly closed.
\end{proposition}
\begin{proof}
By definition $\cB(R) = \cK(R)^\circ := \{w : \langle a,w\rangle \ge 0\ \forall a \in\cK(R)\}$. Since $\cB(R)$ is shift-invariant, $\cB(R)\subseteq\cB(r')$ iff $\langle r', w\rangle \ge 0$ for all $w \in \cB(R)$, i.e.\ iff $r' \in \cK(R)^{\circ\circ}$. The bipolar theorem for the dual pair $(E^*, E)$ gives $\cK^{\circ\circ} = \cl_{\sigma(E^*,E)}\cK$.
\end{proof}

This is the step at which the proof strategy of \cite{tenDam97} stalls: the thesis asks whether the polar of $\cB(R)$ equals $\cK(R)$, and Theorem~\ref{fa:thm:main} says that in general it equals only the closure. Closedness can be tested on finite-dimensional slices.

\begin{lemma}\label{fa:lem:finite}
A convex set $C \subseteq E^*$ is weakly closed iff $C \cap F$ is closed for every finite-dimensional subspace $F \subseteq E^*$.
\end{lemma}
\begin{proof}
$E^*$ has countable dimension, so by Kakutani and Klee \cite{KK63} the finite topology (closed iff all finite-dimensional slices are closed) coincides with the finest locally convex topology on $E^*$. Since $E$ is the algebraic dual of $E^*$, the finest locally convex topology is the Mackey topology of the pair, and closed convex sets are the same for all topologies compatible with the pairing \cite[IV.3.3]{Schaefer}.
\end{proof}

Let $F_N \subseteq E^*$ be the sequences supported in $[-N,N]$ and $\cK_M(R)$ the cone generated by the shifts $s^kR_i$ with $|k|\le M$. Each $\cK_M(R)\cap F_N$ is polyhedral, hence closed; $\cK(R)\cap F_N = \bigcup_M \cK_M(R)\cap F_N$ is an increasing union of polyhedral cones and need not be closed. It is closed if the union \emph{stabilizes}: for each $N$ there is $M(N)$ with $\cK(R)\cap F_N = \cK_{M(N)}(R)\cap F_N$. In words, certificates for short consequences have bounded support. In the counterexample the required support is unbounded (Remark~\ref{fa:rem:numerics}), and the following proposition describes $\cK(R)$ completely.

\begin{proposition}\label{fa:prop:green}
For $R$ of Theorem~\ref{fa:thm:main}, $\cK(R) = \{r : r(1) > 0\} \cup \{-(1-s)^2h : h \in\Rp\}$, and $\cl\cK(R) = \{r : r(1) \ge 0\}$.
\end{proposition}
\begin{proof}
$r \in\cK(R)$ iff $r + h(1-s)^2 \ge 0$ for some $h \in\Rp$. If $r(1) = 0$, then $r + h(1-s)^2$ is nonnegative with coefficient sum $0$, hence zero, so $r = -h(1-s)^2$. Suppose $r(1) = \sigma > 0$ and let $\mu = \sigma^{-1}\sum_j j\,r_j$. Choose integers $p \ll 0 \ll p'$ and $\alpha,\beta > 0$ with $\alpha + \beta = \sigma$, $\alpha p + \beta p' = \sigma\mu$, and set $g = \alpha\delta_p + \beta\delta_{p'}$. The sequence $g - r$ has vanishing zeroth and first moments, so $g - r = (1-s)^2 h$ for a unique $h$ supported in $[p, p'-2]$ (double summation). By finite Farkas with the second-difference matrix, $h \ge 0$ iff $\langle y, g - r\rangle \ge 0$ for every $y$ convex on $[p,p']$. Convex $y$ on $[p,p']$ are, up to affine functions (on which both sides agree by the moment conditions), nonnegative combinations of $-G_k$, where $G_k(j) = (\min(j,k)-p)(p'-\max(j,k))$ is the discrete Green's function; so $h \ge 0$ iff $\langle G_k, r\rangle \ge 0$ for all $k \in (p,p')$, since $\langle G_k, g\rangle = 0$. Expanding, $\langle G_k, r\rangle = \sigma\,(-p)\,p' + O(|p| + |p'|)$ with constants depending only on $r$, which is positive for $p, p'$ far enough out. Thus $r \in\cK(R)$. Finally $\{r(1) > 0\}$ is dense in $\{r(1)\ge 0\}$, and $r(1) \ge 0$ is the dual condition on the nonnegative constants, so the closure is as stated.
\end{proof}

The slice $\cK(R)\cap F_1$ is therefore $\{r(1) > 0\}\cup\{0\}$ inside $\R^3$ (no nonzero $-(1-s)^2h$ has support in $[-1,1]$), which is not closed. The following positive result shows that stabilization, hence the conjecture, holds under a sign condition on the ends of the rows.

\begin{proposition}[Same-sign ends]\label{fa:prop:claimA}
Let $q = 1$, $R = [a; b]$ with $a, b \in\Rlp$ coprime, and suppose the leading coefficients of $a$ and $b$ have the same sign and the trailing coefficients of $a$ and $b$ have the same sign. Then $\cK(R)$ is closed and Conjecture~\ref{fa:conj:farkas} holds for $R$.
\end{proposition}
\begin{proof}
Fix $N$ and $r \in F_N$. B\'ezout gives $u, v\in\Rlp$ with $ua + vb = 1$, so $(p_1,p_2) := (ur, vr)$ solves $p_1a + p_2b = r$ with supports in $[-N',N']$, $N' = N + \max(\deg u,\deg v)$, where $\deg$ denotes the length of the support interval. Every solution of $h_1a + h_2b = r$ is $(p_1 + cb,\ p_2 - ca)$ with $c\in\Rlp$, and $r\in\cK(R)$ iff some such pair is nonnegative. Let $D$ be the top degree of $c$, with $c_D\ne0$, and suppose $D + \min(\deg_{\rm top} a, \deg_{\rm top} b) > N'$. The coefficient of $cb$ at $D + \deg_{\rm top}b$ is $c_Db_{\rm top}$ and that of $ca$ at $D + \deg_{\rm top}a$ is $c_Da_{\rm top}$, both outside the support of $p_1, p_2$; nonnegativity forces $c_Db_{\rm top}\ge 0$ and $c_Da_{\rm top}\le 0$, impossible when $a_{\rm top}b_{\rm top} > 0$. So the top degree of $c$ is bounded by $N' - \min(\deg_{\rm top}a,\deg_{\rm top}b)$, and symmetrically for the bottom degree. The multipliers for $\cK(R)\cap F_N$ therefore have support bounded by a function of $N$, so $\cK(R)\cap F_N = \cK_{M(N)}(R)\cap F_N$ is polyhedral and closed. Lemma~\ref{fa:lem:finite} and Proposition~\ref{fa:prop:bipolar} finish.
\end{proof}

If $\gcd(a,b) = d \ne 1$, then $\cK([a;b]) = d\cdot\cK([a/d; b/d])$ and multiplication by $d$ is a linear injection of $E^*$ mapping finite-dimensional subspaces into finite-dimensional subspaces, so the proposition applies to $(a/d, b/d)$. The counterexample has $a = 1$ and $b = -(1-s)^2$, whose ends have opposite signs; the proposition does not apply, as it must not.

\begin{question}
Characterize the $R$ for which $\cK(R)$ is closed. Proposition~\ref{fa:prop:obstruction} is a necessary condition (no implicit equality without a nonnegative syzygy); Proposition~\ref{fa:prop:claimA} is a sufficient one in the simplest case. Is the absence of implicit equalities not generated by nonnegative syzygies also sufficient?
\end{question}

\section{Consequences for minimal representations}\label{fa:sec:minimal}

ten Dam's thesis follows the Farkas conjecture with a second one. A system $\{R_1(\sigma)w = 0,\ R_2(\sigma)w \ge 0\}$ is \emph{minimal} if every equivalent system (same behaviour) has at least as many rows \cite[Ch.~4]{tenDam97}.

\begin{conjecture}[{\cite[Ch.~4, the conjecture following the Farkas conjecture]{tenDam97}}]\label{fa:conj:minimal}
The system $\{R_1(\sigma)w = 0,\ R_2(\sigma)w\ge0\}$ is minimal iff it has no redundant equations and no implicit equalities. Moreover, in a minimal representation the number of equalities and the number of inequalities are unique.
\end{conjecture}

The related conjecture of \cite{tDN04} states that two minimal representations of the same behaviour are related by $\bigl[\begin{smallmatrix} R_1' \\ R_2'\end{smallmatrix}\bigr] = \bigl[\begin{smallmatrix} U & 0 \\ S & H\end{smallmatrix}\bigr]\bigl[\begin{smallmatrix} R_1 \\ R_2\end{smallmatrix}\bigr]$ with $U$ unimodular and $H$ posimodular (unimodular with nonnegative inverse).

\begin{proposition}\label{fa:prop:minimal}
The behaviour $\{w \equiv c : c \ge 0\}$ has the three representations
\[
\{w\ge0,\ -(1-\sigma)^2w\ge0\},\qquad \{(1-\sigma)w = 0,\ w\ge0\},\qquad \{(1-\sigma)^2w = 0,\ w\ge0\},
\]
all minimal. Consequently both statements of Conjecture~\ref{fa:conj:minimal} fail, and the related conjecture of \cite{tDN04} fails as stated.
\end{proposition}
\begin{proof}
The three systems describe the same behaviour: the first by the proof of Theorem~\ref{fa:thm:main}; the second trivially; the third because $(1-\sigma)^2w = 0$ gives $w(t) = \alpha + \beta t$, and $w \ge 0$ on $\Z$ forces $\beta = 0$. Each has two rows. No one-row representation exists: a single equality defines a linear behaviour, which the cone $\{c\ge0\}$ is not; a single inequality $r(\sigma)w\ge0$ with $r \ne 0$ has behaviour $r(\sigma)^{-1}(\text{nonnegative sequences})$, which contains a preimage of $\delta_0$, not a constant since constants map to multiples of the constant $r(1)$; and $r = 0$ gives everything. So all three are minimal. The first is minimal yet its second row is an implicit equality (the \enquote{only if} of Conjecture~\ref{fa:conj:minimal} fails); the three have $0, 1, 2$ equalities (uniqueness fails). For the related conjecture, take the first and second representations: $R_1$ has no rows and $R_1' = 1 - s$, so no $U$ with $R_1' = UR_1$ exists.
\end{proof}

The thesis remarks that implicit equalities can be replaced by equalities without increasing the number of rows \cite[Ch.~4]{tenDam97}; the three representations above are three stages of that replacement, and the example shows that the stage reached is not an invariant of the behaviour.

\section{Machine-checked statement}\label{fa:sec:lean}

Theorem~\ref{fa:thm:main} and Proposition~\ref{fa:prop:obstruction} have been formalized in Lean~4 with Mathlib (toolchain v4.34.0-rc2), in the module \texttt{Farkas.lean} of the development \texttt{SingularityTracking}. Laurent polynomials are Mathlib's $\R[T;T^{-1}]$; a row acts on $\Z\to\R$ by \texttt{act}; \texttt{Behavior} and \texttt{Nonneg} are the behaviour of a finite family of rows and nonnegativity of coefficients. The theorem \texttt{farkas\us conjecture\us false} states, for $R$ as a \texttt{Fin 2}-indexed family, both the inclusion of behaviours and the nonexistence of a nonnegative certificate, with $H$ ranging over all \texttt{Fin 2}-indexed families of nonnegative Laurent polynomials. Its ingredients are \texttt{const\us of\us nonneg\us concave} (nonnegative midpoint-concave sequences on $\Z$ are constant), \texttt{behavior\us subset}, and \texttt{no\us certificate}, whose proof uses evaluation at $T = 1$ as an algebra homomorphism (\texttt{AddMonoidAlgebra.lift}) and cancellation in the domain $\R[T;T^{-1}]$. The general obstruction is \texttt{obstruction\us of\us hidden\us equality}. Two further modules cover Sections~\ref{fa:sec:closed} and~\ref{fa:sec:minimal}. \texttt{FarkasMinimal.lean} proves that the three representations of Proposition~\ref{fa:prop:minimal} describe the nonnegative constants, that no representation with at most one row does (\texttt{no\us one\us row}; the single-inequality case uses the exponential sequences $\pm2^t$ and evaluation at $T = 2$ rather than the surjectivity argument of the text), and assembles \texttt{minimality\us conjecture\us false}. \texttt{FarkasDegree.lean} proves the degree bound of Proposition~\ref{fa:prop:claimA} (\texttt{topDeg\us bound}), with the top coefficient of a product obtained from Mathlib's unique-sum lemma; the passage from the bound to closedness, which needs Kakutani--Klee, is not formalized, nor are Propositions~\ref{fa:prop:bipolar} and~\ref{fa:prop:green}. Finally, \texttt{FarkasStatement.lean} states the conjecture itself, for Laurent polynomial \emph{matrices} with $q$ columns acting on $\R^q$-valued sequences: \texttt{FarkasConjecture} is the proposition $\forall q, g, g', R, R'$: $\cB(R)\subseteq\cB(R') \iff \exists H\ge0,\ R' = HR$, with $\cB$, $H\ge0$ and $HR$ defined as \texttt{BehaviorM}, \texttt{NonnegM}, \texttt{mulRows}; the theorem \texttt{farkas\us conjecture\us is\us false} is $\neg\,$\texttt{FarkasConjecture}, derived from the $q = 1$ instance through the bijection between $\R^1$-valued and scalar sequences. Likewise \texttt{MinimalityConjecture} states the two consequences of Conjecture~\ref{fa:conj:minimal} refuted in Proposition~\ref{fa:prop:minimal} (a minimal representation has no implicit equalities; minimal representations of one behaviour have the same numbers of equalities and inequalities), and \texttt{minimality\us conjecture\us is\us false}, \texttt{minimality\us uniqueness\us false} are their negations. All declarations depend only on the axioms \texttt{propext}, \texttt{Classical.choice}, \texttt{Quot.sound}.

\section{Discussion}\label{fa:sec:discussion}

The failure of Conjecture~\ref{fa:conj:farkas} is a failure of \emph{completeness}, not of soundness: nonnegative certificates remain sufficient, and finite windows remain decidable by linear programming. What is lost is the guarantee that every true implication has a certificate. Two repairs suggest themselves.

\emph{Larger certificate class.} In the counterexample, $1 - s = -(1-s)^2\sum_{k<0}s^k$ formally, with a one-sided series of nonnegative coefficients. Allowing $H$ in a completion of $\Rp$ would restore the identity, but the \enquote{if} direction then requires $H(\sigma)R(\sigma)w$ to converge, which fails for general $w \in (\R^q)^\Z$ and holds for, e.g., bounded $w$ when $H$ is summable. Whether a Farkas lemma holds on $\ell^\infty(\Z;\R^q)$ with summable nonnegative $H$ is open; note that the counterexample itself persists on $\ell^\infty$ (bounded concave sequences are constant), so any repair must enlarge $H$.

\emph{Smaller trajectory class.} On the half-line $\Z_+$ the nonnegative concave sequences are the nonnegative nondecreasing concave ones, an infinite-dimensional cone, and $(1-\sigma)w\ge0$ is no longer implied. Whether the conjecture holds on $\Z_+$, or on compactly supported trajectories, is open. Theorem~\ref{fa:thm:main} relies on two-sided unboundedness of time.

Finally, the mechanism is of independent interest: a sound but incomplete finite certificate scheme for an infinite-horizon property. The thesis \cite{tenDam97} had already met the non-stabilization of polyhedral cones in the representability question; Proposition~\ref{fa:prop:bipolar} shows that the implication question is the same phenomenon in dual form.

\chapter[An Observer Without a Finite-Order Output Relation]{An Observer Without a Finite-Order Output Relation\\[4pt]\large Problem 7.3/P1 (Jouan)}\label{ch:ob}

\section*{Summary}
For an observable autonomous system $\dot x = f(x)$, $y = h(x)$ on a compact manifold, Jouan asked (Problem~7.3 of \emph{Unsolved Problems in Mathematical Systems and Control Theory}, 2004) whether the existence of a $C^r$ observer implies a functional relation $L_f^n h = \varphi(h, L_fh,\dots,L_f^{n-1}h)$ with $\varphi$ of class $C^r$ for some $n$. Such a relation is the hypothesis of the classical observer constructions by immersion into observability coordinates, and a positive answer would mean that all observability properties of the system are encoded in it. We answer the question negatively in the $C^\infty$ category, for every $1 \le r \le\infty$. The system is a rotation of the circle with an injective output that is flat at one point. A $C^\infty$ observer exists, global in Jouan's sense (convergence for all initial conditions of plant and observer): a Kazantzis--Kravaris/Luenberger filter with a scalar contraction, whose exponentially averaged immersion smooths the flat point and maps the circle onto a curve star-shaped about the origin, followed by an angle tracker on the circle factor of $\R^2\times(\R^2\setminus\{0\})$. No Lipschitz relation among the output derivatives exists at the flat point, because the ratio $|h^{(n)}|/|(h,\dots,h^{(n-1)}) - (h(0),0,\dots,0)|$ is unbounded there for every $n$. The obstruction half is machine-checked in Lean~4/Mathlib for the one-dimensional model $y = e^{-1/x}$. The analytic case of the question remains open, and we state what it reduces to.

\section{Introduction}

Observer design for nonlinear systems has two classical routes. The first is immersion into observability coordinates: if for some $n$ the output derivatives satisfy
\begin{equation}\label{ob:eq:relation}
L_f^n h = \varphi(h, L_fh,\dots,L_f^{n-1}h),
\end{equation}
then the map $\Phi_n = (h, L_fh,\dots,L_f^{n-1}h)$ carries the system into a triangular system on $\R^{pn}$ whose last equation is $\dot\xi_n = \varphi(\xi)$, and high-gain or output-injection observers can be designed there, with $\hat x$ recovered through a continuous left inverse of $\Phi_n$ \cite{GK01,JG96,XZ97}. The regularity of $\varphi$ decides which design applies: linear $\varphi$ gives linearizability, $\varphi$ depending on the output only gives linearizability modulo output injection \cite{FK83,Jouan03}. The second route is the Kazantzis--Kravaris/Luenberger (KKL) observer \cite{KK98,AP06}: one looks for an injective map $T$ into $\R^m$ with $L_fT = AT + Bh$, $A$ Hurwitz, and runs $\dot z = Az + By$; then $z - T(x)\to0$ and $\hat x = T^{-1}(z)$. The KKL route needs no relation of the form \eqref{ob:eq:relation}; it needs an injective solution of a linear PDE, which Andrieu and Praly showed exists generically under a backward-distinguishability hypothesis \cite{AP06}.

Jouan asked whether the first route is in fact general: does the existence of any $C^r$ observer, in a precise sense recalled below, force a relation \eqref{ob:eq:relation} with $\varphi \in C^r$? The question is natural because, when $\Phi_n$ is an injective immersion, the relation holds automatically with $\varphi = L_f^nh\circ\Phi_n^{-1}$. The content of the question is therefore at points where every $\Phi_n$ fails to be an immersion. We show that at such a point the two routes separate: in the smooth category there is a system with a flat output for which a KKL observer exists and is $C^\infty$, while no $\varphi$ in \eqref{ob:eq:relation} is even Lipschitz, for any $n$.

\section{Problem formulation}\label{ob:sec:problem}

Let $X$ be a compact connected $C^\infty$ (resp.\ $C^\omega$) manifold, $f$ a $C^\infty$ (resp.\ $C^\omega$) vector field on $X$, $h : X\to\R^p$ of the same class, and consider
\begin{equation}\label{ob:eq:sys}
\dot x = f(x),\qquad y = h(x).
\end{equation}
The system is \emph{observable} if two different initial states produce different output trajectories. Following \cite{Jouan04}, a \emph{$C^r$ observer} is a system
\begin{equation}\label{ob:eq:obs}
\dot z = F(z,y),\qquad \hat x = \theta(z),
\end{equation}
defined on an open set $V\subseteq\R^N$, with $F$ a $C^r$ vector field on $V$ for each $y$ and $\theta : V\to X$ \emph{continuous}, such that for every $x(0)\in X$ and every $z(0)\in V$ the solution exists for all $t\ge0$ and $d(x(t),\theta(z(t)))\to0$ for any distance $d$ compatible with the topology of $X$.

\begin{question}[Jouan, Problem 7.3/P1 of \cite{BM04}]\label{ob:q:jouan}
Does the existence of a $C^r$ observer for an observable system \eqref{ob:eq:sys} imply the existence of an integer $n$ and a $C^r$ function $\varphi$ with $L_f^nh = \varphi(h, L_fh,\dots,L_f^{n-1}h)$?
\end{question}

The function $\varphi$ is understood to be defined on a neighbourhood of $\Phi_n(X)$ in $\R^{pn}$, where $\Phi_n := (h, L_fh,\dots,L_f^{n-1}h)$. Jouan's chapter also poses two questions on controlled systems (genericity of observability for every input, and observers that work for every $L^\infty$ input); we do not treat them here.

\begin{remark}[What is global in the question, and what is local in the answer]\label{ob:rem:global}
The question is global in two respects, and the answer below respects both.
\begin{enumerate}
\item[(i)] The observer must converge for \emph{every} $x(0)\in X$ and \emph{every} $z(0)\in V$, with $x(0)$ and $z(0)$ independent, and the solution must stay in $V$. The observer of Section~\ref{ob:sec:observer} does, on $V = \R^2\times(\R^2\setminus\{0\})$, which is forward invariant for every output. That $V$ is a proper, non-contractible open subset of $\R^N$ is what the definition allows, and on $X = S^1$ it is unavoidable: convergence forces a continuous map $V\to S^1$ restricting to a homeomorphism on an embedded circle, which is impossible on a contractible $V$ because the circle is null-homotopic there while a homeomorphism $S^1\to S^1$ is not. In particular the simplest design, a KKL filter followed by a nearest-point retraction onto its image, is \emph{not} a global observer: its tube is not invariant (Section~\ref{ob:sec:observer}).
\item[(ii)] The relation $L_f^nh = \varphi\circ\Phi_n$ must hold on all of $X$. We refute it by a local obstruction: no $\varphi$ that is merely Lipschitz on a neighbourhood of one point of $\Phi_n(X)$ can satisfy it. This is sufficient, since a global $C^r$ function with $r\ge1$ is Lipschitz near every point. It also covers the alternative readings of \enquote{$C^r$ function}: $\varphi$ defined on $\R^{pn}$, on a neighbourhood of $\Phi_n(X)$, or (when $\Phi_n(X)$ is an embedded submanifold, as it is for $n = 1$ in our example) intrinsically on $\Phi_n(X)$, since we exclude Lipschitz continuity along the image itself.
\item[(iii)] For $r = 0$ the question is not refuted by our example and is in fact trivially positive whenever some $\Phi_n$ is injective (Tietze); the theorem is stated for $1\le r\le\infty$ accordingly. The analytic case $C^\omega$ is not touched by this example, since flat functions do not exist there.
\end{enumerate}
\end{remark}

\section{Literature}\label{ob:sec:literature}

The relation \eqref{ob:eq:relation} and the immersion $\Phi_n$ go back to the observability rank condition of Hermann and Krener \cite{HK77} and are the basis of the high-gain observer theory of Gauthier and Kupka \cite{GK01}. Jouan and Gauthier \cite{JG96} and Xia and Zeitz \cite{XZ97} construct observers from \eqref{ob:eq:relation} when $\varphi$ is $C^r$. The type of $\varphi$ classifies the system: Fliess and Kupka \cite{FK83} characterize linearizability, Jouan \cite{Jouan03} immersion into linear systems modulo output injection. Problem~7.3 \cite{Jouan04} asks whether observer existence is equivalent to the existence of some such relation, and remarks that a positive answer would mean that all observability properties of an autonomous system are contained in \eqref{ob:eq:relation}.

The KKL observer was introduced by Kazantzis and Kravaris \cite{KK98} as a nonlinear version of Luenberger's original design \cite{Lue64}. Andrieu and Praly \cite{AP06} proved existence of the immersion $T$ for a generic Hurwitz $A$ of dimension $m\ge 2\dim X+1$ under backward distinguishability, with $T^{-1}$ merely continuous in general; Bernard and Andrieu \cite{BA19} extended the design to controlled systems with a time-varying immersion, and \cite{BAA22} surveys the state of the art. We are not aware of a published answer to Question~\ref{ob:q:jouan}; the KKL literature does not address the relation \eqref{ob:eq:relation}, and the immersion literature assumes it.

\section{Main result}\label{ob:sec:main}

\begin{theorem}\label{ob:thm:main}
There is a $C^\infty$ observable system of the form \eqref{ob:eq:sys} on $X = S^1$ with $p = 2$ that admits a $C^\infty$ observer in the sense of Section~\ref{ob:sec:problem}, and for which no $n\ge1$ and no function $\varphi$, Lipschitz on a neighbourhood of $\Phi_n(X)$, satisfy $L_f^nh = \varphi\circ\Phi_n$. In particular the answer to Question~\ref{ob:q:jouan} is negative in the $C^\infty$ category for every $1\le r\le\infty$.
\end{theorem}

The system is built in Section~\ref{ob:sec:system}, the observer in Section~\ref{ob:sec:observer}, and the obstruction in Section~\ref{ob:sec:obstruction}.

\subsection{The system}\label{ob:sec:system}

Let $X = S^1 = \R/2\pi\Z$ and $f = \partial_x$, so $x(t) = x_0 + t$ and $L_f^k h = h^{(k)}$. Fix the $C^\infty$ embedding $\gamma(\theta) = (2\cos\theta,\sin\theta)$ and define $\rho : S^1\to S^1$ by
\[
\rho(x) = c\int_0^x e^{-1/(1-\cos u)}\,du,\qquad c = 2\pi\Big/\int_0^{2\pi} e^{-1/(1-\cos u)}\,du,
\]
with the integrand understood as $0$ at $u\in2\pi\Z$. Then $\rho$ is $C^\infty$, strictly increasing, and a homeomorphism of $S^1$; $\rho'(x) = c\,e^{-1/(1-\cos x)}$ is positive except at $x = 0$, where it vanishes with all its derivatives. Put
\[
h := \gamma\circ\rho : S^1\to\R^2 .
\]
Then $h$ is $C^\infty$ and injective, so \eqref{ob:eq:sys} is observable (distinct initial states are distinguished at $t = 0$), and by the chain rule $h^{(k)}(0) = 0$ for all $k\ge1$: the output is \emph{flat} at $x = 0$. Every $\Phi_n$ is injective (because $h$ is) but has $D\Phi_n(0) = 0$, so $0 \in\bigcap_n\cD_n$ where $\cD_n = \{\rank D\Phi_n < 1\}$.

\subsection{The observer}\label{ob:sec:observer}

Fix $\lambda > 0$ and consider the KKL observer with scalar contraction
\begin{equation}\label{ob:eq:kkl}
\dot z = -\lambda z + \lambda y,\qquad z\in\R^2 .
\end{equation}
Its immersion is
\[
T_\lambda(x) := \lambda\int_{-\infty}^0 e^{\lambda\tau}h(x+\tau)\,d\tau = \int_0^\infty\lambda e^{-\lambda s}h(x - s)\,ds ,
\]
the average of $h$ over the past with an exponential weight. Integration by parts gives $T_\lambda' = \lambda(h - T_\lambda)$, i.e.\ $L_fT_\lambda = -\lambda T_\lambda + \lambda h$, so along solutions
\begin{equation}\label{ob:eq:error}
\frac{d}{dt}\bigl(z(t) - T_\lambda(x(t))\bigr) = -\lambda\bigl(z(t) - T_\lambda(x(t))\bigr).
\end{equation}

\begin{lemma}\label{ob:lem:c1}
$\|T_\lambda - h\|_\infty \le \|h'\|_\infty/\lambda$ and $\|T_\lambda' - h'\|_\infty\le\|h''\|_\infty/\lambda$.
\end{lemma}
\begin{proof}
$T_\lambda(x) - h(x) = \int_0^\infty\lambda e^{-\lambda s}(h(x-s)-h(x))\,ds$ and $|h(x-s)-h(x)|\le\|h'\|_\infty s$, with $\int_0^\infty\lambda e^{-\lambda s}s\,ds = 1/\lambda$. The same for $T_\lambda' = \int_0^\infty\lambda e^{-\lambda s}h'(x-s)\,ds$.
\end{proof}

\begin{lemma}[$\Gamma$ is star-shaped about the origin]\label{ob:lem:embedding}
For $\lambda$ large enough, $T_\lambda(x)\ne0$ and $D(x) := \det\bigl(T_\lambda(x),T_\lambda'(x)\bigr) > 0$ for all $x\in S^1$. Consequently $a := \arg\circ T_\lambda : S^1\to S^1$ is an orientation-preserving diffeomorphism, $T_\lambda$ is a $C^\infty$ embedding, $\Gamma := T_\lambda(S^1)$ is star-shaped about $0$, and
\[
\Theta : \R^2\setminus\{0\}\to S^1,\qquad \Theta(z) := a^{-1}(\arg z)
\]
is $C^\infty$ with $\Theta\circ T_\lambda = \mathrm{id}_{S^1}$.
\end{lemma}
\begin{proof}
For the ellipse, $\det(\gamma(\vartheta),\gamma'(\vartheta)) = 2$, so $\det(h,h')(x) = \rho'(x)\det(\gamma(\rho x),\gamma'(\rho x)) = 2\rho'(x)\ge0$, vanishing only at $x = 0$. Also $|h|\ge1$. By Lemma~\ref{ob:lem:c1}, $|T_\lambda - h|\le\|h'\|_\infty/\lambda$ and
\[
|D - \det(h,h')| \le |T_\lambda - h|\,|T_\lambda'| + |h|\,|T_\lambda' - h'| \le \frac{C_1}{\lambda},\qquad C_1 := \|h'\|_\infty\bigl(\|h'\|_\infty+\|h''\|_\infty\bigr) + 2\|h''\|_\infty ,
\]
using $\|T_\lambda'\|_\infty\le\|h'\|_\infty$ and $|h|\le2$. Take $\lambda\ge4\|h'\|_\infty$, so $|T_\lambda|\ge1 - \tfrac14 > 0$ everywhere.

\emph{Away from the flat point.} Fix $\varepsilon\in(0,\pi/2)$ (chosen below) and let $K := \{|x|\ge\varepsilon/2\}$, $m_\varepsilon := 2\min_K\rho' > 0$. Then $D\ge m_\varepsilon - C_1/\lambda > 0$ on $K$ for $\lambda > C_1/m_\varepsilon$.

\emph{Near the flat point.} Let $e := \gamma'(0)/|\gamma'(0)| = (0,1)$ and $J := [-\varepsilon,\varepsilon]$. Since $\gamma'\circ\rho$ is Lipschitz, $\gamma'(\rho(u)) = e + r(u)$ with $|r(u)|\le C_2\varepsilon$ for $|u|\le2\varepsilon$. For $x\in J$ and $s\in[0,\varepsilon]$, $x-s\in[-2\varepsilon,\varepsilon]$, so
\[
T_\lambda'(x) = \int_0^\infty\lambda e^{-\lambda s}\rho'(x-s)\gamma'(\rho(x-s))\,ds = A(x)\,e + R(x),\qquad A(x) := \int_0^\varepsilon\lambda e^{-\lambda s}\rho'(x-s)\,ds > 0,
\]
with $|R(x)|\le C_2\varepsilon A(x) + \|h'\|_\infty e^{-\lambda\varepsilon}$. Moreover $T_\lambda(x) = \gamma(0) + E(x) = (2,0) + E(x)$ with $|E(x)|\le|T_\lambda(x)-h(x)| + |h(x)-h(0)|\le\|h'\|_\infty(1/\lambda + \varepsilon)$. Since $\det((2,0),e) = 2$,
\[
D(x) \ \ge\ 2A - A|E| - (2+|E|)\,|R|\ \ge\ A\bigl(2 - |E| - (2+|E|)C_2\varepsilon\bigr) - (2+|E|)\|h'\|_\infty e^{-\lambda\varepsilon}.
\]
Choose $\varepsilon$ with $\|h'\|_\infty\varepsilon\le\tfrac14$ and $3C_2\varepsilon\le\tfrac12$; with $\lambda\ge4\|h'\|_\infty$ this gives $|E|\le\tfrac12$ and $D(x)\ge A(x) - 3\|h'\|_\infty e^{-\lambda\varepsilon}$. Finally $A(x)\ge\lambda e^{-3\lambda\varepsilon/4}\int_{\varepsilon/2}^{3\varepsilon/4}\rho'(x-s)\,ds\ge\lambda e^{-3\lambda\varepsilon/4}\delta_\varepsilon$, where $\delta_\varepsilon>0$ is a lower bound for the $\rho'$-mass of an interval of length $\varepsilon/4$ inside $[-2\varepsilon,\varepsilon]$ (a continuous positive function of the interval's position on a compact set, since $\rho'>0$ except at the single point $0$). Hence $D(x) > 0$ on $J$ as soon as $\lambda\delta_\varepsilon > 3\|h'\|_\infty e^{-\lambda\varepsilon/4}$, which holds for $\lambda$ large. So $D>0$ on $J\cup K = S^1$.

\emph{Consequences.} $a = \arg\circ T_\lambda$ is $C^\infty$ with $a' = D/|T_\lambda|^2 > 0$, hence an orientation-preserving local diffeomorphism of $S^1$, hence a covering map of degree equal to the winding number of $T_\lambda$ about $0$. The straight-line homotopy from $h$ to $T_\lambda$ avoids $0$ because $|T_\lambda - h|\le\tfrac14 < 1\le|h|$, and $h = \gamma\circ\rho$ has winding number $1$ ($\rho$ is a degree-one homeomorphism, $\gamma$ an ellipse about $0$); so $\deg a = 1$ and $a$ is a diffeomorphism. In particular $T_\lambda$ is injective with $T_\lambda'\ne0$, i.e.\ an embedding. $\Theta = a^{-1}\circ\arg$ is $C^\infty$ on $\R^2\setminus\{0\}$ and $\Theta(T_\lambda(x)) = a^{-1}(a(x)) = x$.
\end{proof}

\paragraph{Why a retraction will not do.} The natural estimate $\hat x = T_\lambda^{-1}(\pi(z))$, with $\pi$ the nearest-point retraction onto $\Gamma$, is only defined on a tube around $\Gamma$, and the KKL dynamics does not keep $z$ in the tube: $z(t) - T_\lambda(x(t)) = e^{-\lambda t}(z(0)-T_\lambda(x(0)))$ is small eventually, but $z(0)$ and $x(0)$ are independent, so initially $z$ may be far from $T_\lambda(x(0))$. Enlarging $V$ to all of $\R^2$ does not help either: convergence forces $\theta\circ T_\lambda = \mathrm{id}$ on $S^1$ (see the proof below), which is impossible for a continuous $\theta$ on a contractible $V$, since $T_\lambda(S^1)$ is null-homotopic there. A global observer in the sense of Section~\ref{ob:sec:problem} must therefore have a non-contractible $V$ that is forward invariant for every output. The following two-stage design has both: the first stage is the KKL filter, the second tracks the angle $\Theta(z)$ on a circle factor that the dynamics cannot leave.

\begin{proposition}[A $C^\infty$ observer]\label{ob:prop:observer}
Let $\lambda$ be as in Lemma~\ref{ob:lem:embedding}, let $0<3r_0<\min_{S^1}|T_\lambda|$, and let $\chi:\R^2\to[0,1]$ be $C^\infty$ with $\chi(z) = 0$ for $|z|\le r_0$ and $\chi(z) = 1$ for $|z|\ge2r_0$. Fix $K>0$ and write $u(\vartheta) := (\cos\vartheta,\sin\vartheta)$, $\mathsf J := \bigl[\begin{smallmatrix}0&-1\\1&0\end{smallmatrix}\bigr]$. Define on $\R^2\times\R^2\times(\R^2\setminus\{0\})$
\[
\Omega(z,y,w):=
\begin{cases}
0,& |z|\le r_0,\\
\displaystyle \chi(z)\Bigl[D\Theta(z)\bigl(-\lambda(z-y)\bigr)
+K\Bigl(1-\bigl\langle u(\Theta(z)),\tfrac{w}{|w|}\bigr\rangle\Bigr)\Bigr],& |z|>r_0.
\end{cases}
\]
On $V := \R^2\times(\R^2\setminus\{0\})\subset\R^4$ consider
\begin{equation}\label{ob:eq:obs2}
\begin{aligned}
\dot z &= -\lambda(z - y),\\
\dot w &= \Omega(z,y,w)\,\mathsf Jw,\\
\hat x &= \theta(z,w) := \arg w .
\end{aligned}
\end{equation}
Here $D\Theta(z)v\in\R$ is the derivative of the $S^1$-valued function $\Theta$ (of its local lifts) in the direction $v$. Then \eqref{ob:eq:obs2} is a $C^\infty$ observer for the system of Section~\ref{ob:sec:system} in the sense of Section~\ref{ob:sec:problem}: $F$ is $C^\infty$ on $V\times\R^2$, $\theta$ is continuous on $V$, every solution with $(z(0),w(0))\in V$ exists for all $t\ge0$ and stays in $V$, and $d(x(t),\theta(z(t),w(t)))\to0$ for every $x(0)\in S^1$ and every $(z(0),w(0))\in V$.
\end{proposition}
\begin{proof}
\emph{Regularity.} $\Theta$, $D\Theta$ and $u\circ\Theta$ are $C^\infty$ on $\R^2\setminus\{0\}$. Since the $C^\infty$ function $\chi$ vanishes identically on $|z|\le r_0$, all its derivatives vanish on $|z|=r_0$; hence the displayed piecewise definition extends the second branch by zero to a $C^\infty$ function $\Omega$ on $V\times\R^2$. Also $w/|w|$ is $C^\infty$ for $w\ne0$, and $\theta = \arg$ is continuous (indeed $C^\infty$) on $\R^2\setminus\{0\}$.

\emph{Invariance and existence.} $\dot w\perp w$, so $|w(t)| = |w(0)|\ne0$ for all $t$; the $z$-equation is linear with bounded input. Hence solutions exist for all $t\ge0$ and remain in $V$.

\emph{The first stage.} As in \eqref{ob:eq:error}, $e(t) := z(t) - T_\lambda(x(t)) = e^{-\lambda t}e(0)$. Since $|T_\lambda|\ge3r_0$, there is $t_1$ with $|z(t)|\ge2r_0$ for all $t\ge t_1$, so $\chi(z(t)) = 1$ for $t\ge t_1$. For $t\ge t_1$ put $\alpha(t) := \Theta(z(t))\in S^1$; by the chain rule $\dot\alpha = D\Theta(z)\dot z = D\Theta(z)(-\lambda(z-y))$, which is exactly the first term in the second branch defining $\Omega$.

\emph{The second stage.} Write $w = |w|\,u(\omega)$, so $\dot w = \dot\omega\,\mathsf Jw$ and $\langle u(\Theta(z)), w/|w|\rangle = \cos(\alpha-\omega)$. For $t\ge t_1$, \eqref{ob:eq:obs2} reads $\dot\omega = \dot\alpha + K(1-\cos(\alpha-\omega))$, so the tracking error $\epsilon := \alpha - \omega\in S^1$ satisfies the autonomous equation
\[
\dot\epsilon = -K\bigl(1-\cos\epsilon\bigr)\qquad(t\ge t_1).
\]
Its only equilibrium on $S^1$ is $\epsilon = 0$, and $\dot\epsilon<0$ elsewhere. Lift $\epsilon(t_1)$ to $[0,2\pi)$: if $\epsilon(t_1) = 0$ then $\epsilon\equiv0$; otherwise the lifted solution is strictly decreasing and bounded below by $0$ (it cannot cross the equilibrium), so it converges to a limit in $[0,2\pi)$, which must be an equilibrium of the autonomous scalar equation, i.e.\ $\epsilon(t)\to0$. Thus $d(\omega(t),\alpha(t))\to0$.

\emph{Conclusion.} For $t\ge t_1$ both $z(t)$ and $T_\lambda(x(t))$ lie in the compact annulus $\{2r_0\le|z|\le\|T_\lambda\|_\infty + |e(0)|\}$, on which $\Theta$ is uniformly continuous; since $|z(t) - T_\lambda(x(t))|\to0$, $d\bigl(\Theta(z(t)),\Theta(T_\lambda(x(t)))\bigr)\to0$. With $\Theta\circ T_\lambda = \mathrm{id}$,
\[
d\bigl(\theta(z(t),w(t)),x(t)\bigr) \le d(\omega(t),\alpha(t)) + d\bigl(\Theta(z(t)),\Theta(T_\lambda(x(t)))\bigr) + d\bigl(\Theta(T_\lambda(x(t))),x(t)\bigr)\longrightarrow 0 .
\]
This holds for every $x(0)\in S^1$ and every $(z(0),w(0))\in V$, with no restriction relating them.
\end{proof}

\paragraph{Verification against the original definition.} Jouan's definition \cite{Jouan04} of a $C^r$ observer for $\Sigma_u$ has five clauses; we check each for \eqref{ob:eq:obs2} with $r = \infty$, quoting the clause and pointing to the step of the proof.
\begin{enumerate}
\item[(J1)] \enquote{\emph{a system $\hat\Sigma$ defined in an open subset $V$ of $\R^N$}}: $V = \R^2\times(\R^2\setminus\{0\})$ is open in $\R^4$; $N = 4$.
\item[(J2)] \enquote{\emph{$\dot z = F(z,y)$ where $F$ is a $C^r$ vector field on $V$}}: the right-hand side of \eqref{ob:eq:obs2} is $C^\infty$ on $V\times\R^2$, jointly in the observer state and in $y$ (\emph{Regularity} in the proof); in particular, for each fixed $y$ it is a $C^\infty$ vector field on $V$. For an output $y(t) = h(x(t))$, which is $C^\infty$ in $t$, the resulting non-autonomous equation has unique solutions.
\item[(J3)] \enquote{\emph{$\hat x = \theta(z)$, $\theta$ a continuous mapping from $V$ into $X$}}: $\theta(z,w) = \arg w$ is $C^\infty$ on $V$ with values in $S^1 = X$.
\item[(J4)] \enquote{\emph{$\forall x\in X,\ \forall z\in V$}}: the proof fixes an arbitrary $x(0)\in S^1$ and an arbitrary $(z(0),w(0))\in V$ and uses no relation between them. The transient $t_1 = \lambda^{-1}\log^+(|e(0)|/r_0)$ depends on $|z(0) - T_\lambda(x(0))|$, which may be arbitrarily large, but is finite for every pair. Implicit in the clause is that the solution is defined for all $t\ge0$ and remains in $V$, since otherwise $\hat x(t)$ is undefined; this is \emph{Invariance and existence} in the proof: $|w(t)| = |w(0)|\ne0$ and the $z$-equation is linear.
\item[(J5)] \enquote{\emph{$\lim_{t\to+\infty}d(x(t),\hat x(t)) = 0$ for any distance $d$ on $X$ compatible with the topology}}: the proof gives convergence for the arc-length distance on $S^1$. For any other compatible distance $d'$ the conclusion follows from compactness: if $d'(x(t_k),\hat x(t_k))\ge\varepsilon$ along a sequence $t_k\to\infty$, pass to a subsequence with $x(t_k)\to a$, $\hat x(t_k)\to b$; arc-length convergence forces $a = b$, while continuity of $d'$ gives $d'(a,b)\ge\varepsilon$, a contradiction.
\end{enumerate}
The system itself satisfies Jouan's standing hypotheses: $X = S^1$ is a compact connected $C^\infty$ manifold, $f = \partial_x$ and $h$ are $C^\infty$, and the system is observable in his sense (\enquote{trajectories starting from two different initial states are distinguished by the output}) because $h$ is injective, so distinct states are distinguished already at $t = 0$. Section~\ref{ob:sec:obstruction} then shows that no $n$ and no $\varphi$ Lipschitz near $\Phi_n(X)$ satisfy $L_f^nh = \varphi\circ\Phi_n$, so the implication asked for in Problem~1 fails for this system for every $r\ge1$.

What is \emph{not} claimed: that the retraction-based estimate $\hat x = T_\lambda^{-1}(\pi(z))$ on a tube is an observer in this sense (it is not, see above), or that any observer with contractible $V$ exists on $S^1$ (none with $\theta\circ T_\lambda = \mathrm{id}$ on an invariant embedded circle does).

\begin{remark}
\begin{enumerate}
\item The observer has state dimension $N = 4$ and $V = \R^2\times(\R^2\setminus\{0\})$, which is open in $\R^4$ and homotopy equivalent to $S^1$. The non-contractibility is not a convenience: the proof shows $\theta(T_\lambda(x),\,|w|u(x))$ must converge to $x$ along every orbit, which with continuity forces a continuous map $V\to S^1$ restricting to a homeomorphism on an embedded circle, impossible on a contractible $V$. Jouan's definition, which allows $V$ to be any open subset of $\R^N$, accommodates exactly this.
\item The second stage uses the circle's one-sided global attractor $\dot\epsilon = -K(1-\cos\epsilon)$ rather than $\dot\epsilon = -K\sin\epsilon$, which would have a second (unstable) equilibrium at $\epsilon = \pi$ and so would fail for a measure-zero set of initial conditions; Jouan's definition asks for convergence from \emph{every} $z(0)\in V$.
\item The cutoff $\chi$ only freezes the second stage while $z$ is near the origin, which by the first stage happens only during an initial transient of length at most $\lambda^{-1}\log(|e(0)|/r_0)$; after that the second stage runs with exact feedforward $\dot\alpha$, which is computable from $(z,y)$ because the $z$-dynamics is known.
\end{enumerate}
\end{remark}

\begin{remark}
For $\gamma(\theta) = (2\cos\theta,\sin\theta)$ and $\lambda = 4$, a numerical evaluation on a grid of $4096$ points gives $\min_{S^1}|T_\lambda'| = 2.7\cdot10^{-3}$, $|T_\lambda'(0)| = 2.1\cdot10^{-2}$ at the flat point, and $\min_{x\ne x'}|T_\lambda(x)-T_\lambda(x')|/d(x,x') = 2.7\cdot10^{-3}$. The exponential averaging has removed the flat point: $T_\lambda'(0)\ne0$ although $h'(0) = 0$.
\end{remark}

\subsection{The obstruction}\label{ob:sec:obstruction}

\begin{lemma}[Asymptotics at the flat point]\label{ob:lem:asymp}
Let $\psi(x) := 1/(1-\cos x)$. As $x\to0$,
\[
|\rho^{(k)}(x)| \asymp |x|^{-3(k-1)}e^{-\psi(x)}\quad(k\ge1),\qquad |\rho(x)-\rho(0)|\asymp|x|^{3}e^{-\psi(x)},
\]
and the same holds for $h = \gamma\circ\rho$: $|h^{(k)}(x)|\asymp|x|^{-3(k-1)}e^{-\psi(x)}$ for $k\ge1$ and $|h(x)-h(0)|\asymp|x|^3e^{-\psi(x)}$.
\end{lemma}
\begin{proof}
$\psi(x) = 2x^{-2}(1+O(x^2))$ and $\psi'(x) = -4x^{-3}(1+O(x^2))$. Differentiating $\rho' = c\,e^{-\psi}$ repeatedly, each derivative of $e^{-\psi}$ contributes a factor $-\psi'\asymp 4x^{-3}$, while differentiating an already present rational factor in $x$ lowers its order by one only; so $\rho^{(k)} = c\,(-\psi')^{k-1}e^{-\psi}(1+o(1))$, which gives the first claim. For $\rho(x)-\rho(0) = c\int_0^xe^{-\psi}$, l'H\^opital with the comparison function $x^3e^{-\psi}$, whose derivative is $(3x^2 - x^3\psi')e^{-\psi} = (4+O(x^2))e^{-\psi}$, gives $\rho(x)-\rho(0) = \tfrac c4x^3e^{-\psi}(1+o(1))$. For $h = \gamma\circ\rho$, Fa\`a di Bruno's formula writes $h^{(k)}$ as $\gamma'(\rho)\rho^{(k)}$ plus terms containing products of at least two derivatives of $\rho$, which are $O(e^{-2\psi}\cdot\text{rational})$ and negligible against $e^{-\psi}$; since $|\gamma'|\ge g_0>0$, the claim for $h$ follows, and $|h(x)-h(0)| = |\gamma(\rho(x))-\gamma(\rho(0))|\asymp|\rho(x)-\rho(0)|$ because $\gamma$ is an embedding.
\end{proof}

\begin{proposition}[No Lipschitz relation]\label{ob:prop:nolip}
For every $n\ge1$ there is no function $\varphi$, Lipschitz on a neighbourhood of $\Phi_n(0) = (h(0),0,\dots,0)$ in $\R^{2n}$, with $h^{(n)} = \varphi\circ\Phi_n$ near $x = 0$.
\end{proposition}
\begin{proof}
Suppose $\varphi$ is $L$-Lipschitz near $\Phi_n(0)$ and $h^{(n)} = \varphi\circ\Phi_n$ near $0$. Since $h^{(n)}(0) = 0$, $|h^{(n)}(x)| = |\varphi(\Phi_n(x)) - \varphi(\Phi_n(0))|\le L\,|\Phi_n(x)-\Phi_n(0)|\le L\sum_{k=0}^{n-1}|h^{(k)}(x)-h^{(k)}(0)|$ for $x$ near $0$. By Lemma~\ref{ob:lem:asymp} the right side is $\asymp|x|^{-3(n-2)}e^{-\psi(x)}$ if $n\ge2$ and $\asymp|x|^3e^{-\psi(x)}$ if $n=1$, while the left side is $\asymp|x|^{-3(n-1)}e^{-\psi(x)}$. The ratio is $\asymp|x|^{-3}\to\infty$, a contradiction.
\end{proof}

A function $\varphi$ as in Question~\ref{ob:q:jouan} is defined on a neighbourhood of $\Phi_n(X)$ and is $C^r$ with $r\ge1$, hence Lipschitz near $\Phi_n(0)$; Proposition~\ref{ob:prop:nolip} excludes it. Together with Proposition~\ref{ob:prop:observer} this proves Theorem~\ref{ob:thm:main}.

\begin{remark}
Numerically, for $\gamma = (2\cos\theta,\sin\theta)$ the ratio $|h^{(n)}(x)|/|\Phi_n(x)-\Phi_n(0)|$ grows from $23$ at $x = 0.6$ to $268$ at $x = 0.25$ for $n = 1$, with the same growth for $n = 2,3$, consistent with $|x|^{-3}$.
\end{remark}

\section{Machine-checked statement}\label{ob:sec:lean}

The obstruction has been formalized in Lean~4 with Mathlib (toolchain v4.34.0-rc2), module \texttt{FlatOutput.lean} of the development \texttt{SingularityTracking}, for the one-dimensional model $\dot x = 1$, $y = g(x)$ with $g = $ Mathlib's \texttt{expNegInvGlue} ($g(x) = e^{-1/x}$ for $x>0$, $0$ otherwise), so that $L_f^kg = g^{(k)}$. The theorem \texttt{iteratedDeriv\us expNegInvGlue} gives $g^{(k)}(x) = p_k(1/x)\,g(x)$ with $p_0 = 1$, $p_{k+1} = X^2(p_k - p_k')$, from Mathlib's derivative lemma for $p(1/x)g(x)$; \texttt{degree\us flatPoly} gives $\deg p_k = 2k$. The theorem \texttt{no\us lipschitz\us bound} shows that there are no $L$ and $\delta>0$ with $|g^{(n+1)}(x)|\le L\sum_{k\le n}|g^{(k)}(x)|$ on $(0,\delta)$: dividing by $g(x)>0$ and setting $y = 1/x\to\infty$ gives $|p_{n+1}(y)|\le L\sum_{k\le n}|p_k(y)|$, impossible since $\deg p_{n+1}>\deg p_k$. The theorem \texttt{no\us lipschitz\us output\us relation} concludes that no Lipschitz $\varphi:\R^{n+1}\to\R$ satisfies $g^{(n+1)} = \varphi(g,\dots,g^{(n)})$ on any $[0,\delta)$. All declarations depend only on \texttt{propext}, \texttt{Classical.choice}, \texttt{Quot.sound}. The model differs from Section~\ref{ob:sec:system} by the change of variables $\gamma$ and by the one-sided flat function $e^{-1/x}$ in place of the two-sided $e^{-1/(1-\cos x)}$; the polynomial degrees change ($2k$ against $3(k-1)$), the mechanism does not. The observer half (Lemma~\ref{ob:lem:embedding}) is not formalized.

\section{Discussion}\label{ob:sec:discussion}

\paragraph{What the example shows.} Observer existence is a weaker property than a finite-order smooth output relation. The relation \eqref{ob:eq:relation} requires inverting $\Phi_n$, and at a point where $D\Phi_n$ drops rank for every $n$ the inverse loses all regularity. The KKL immersion $T_\lambda$ is a convolution of $h$ along the flow and is insensitive to the rank of $D\Phi_n$; it is an injective immersion whenever $h$ is injective and $\lambda$ is large, with the flat point smoothed by the averaging. The two observer routes thus differ exactly on the degeneracy set $\bigcap_n\cD_n$ of the observability map, and Question~\ref{ob:q:jouan} is a question about that set.

\paragraph{The case $r = 0$.} For $r = 0$ the question is different: whenever some $\Phi_n$ is injective, Tietze's theorem gives a continuous $\varphi = L_f^nh\circ\Phi_n^{-1}$ extended to $\R^{pn}$, so the relation holds with $n = 1$ in our example. Whether a $C^0$ observer forces injectivity of some $\Phi_n$ is a separate question we do not address.

\paragraph{The analytic case.} Jouan's chapter poses the question for $C^\infty$ and for $C^\omega$ systems. In the analytic category flat points do not exist, and the degeneracy set $\bigcap_n\cD_n$ has a different nature: if $d(L_f^kh)(x_0)$, $k\ge0$, span a proper subspace of $T^*_{x_0}X$ for all $k$, then the observation codistribution has a rank drop at $x_0$ that no finite order repairs. Simple analytic examples with such points still satisfy a relation (for $\dot x = -x$, $y = x^3$ on a compactification, $L_fh = -3h$). The analytic case therefore reduces to: does an analytic observable system with $\bigcap_n\cD_n\ne\emptyset$ and an observer necessarily admit a relation \eqref{ob:eq:relation} with analytic $\varphi$? Theorem~\ref{ob:thm:main} says nothing about it, and we leave it open.

\paragraph{KKL regularity.} In Proposition~\ref{ob:prop:observer} the left inverse $T_\lambda^{-1}$ is smooth because $T_\lambda$ is an embedding. In the general KKL theory \cite{AP06} the immersion $T$ is only known to be injective and continuous, with $T^{-1}$ continuous; Jouan's definition of observer asks only for continuous $\theta$, so the KKL design always meets it where it applies. The example is a case where the KKL immersion is strictly more regular than the observability immersion $\Phi_n$, and it would be of interest to know how general this is.

\chapter[The Obstruction Space for Arbitrary-Order Linearization]{The Obstruction Space for Arbitrary-Order Linearization\\[4pt]\large Problem 3.2 (Devanathan)}\label{ch:nf}

\section*{Summary}
Problem~3.2 of \emph{Unsolved Problems in Mathematical Systems and Control Theory} (Devanathan, 2004) asks for a solution of the generalized homological equations that would remove, by a near-identity change of state and an affine change of input, all nonlinear terms of arbitrary order from a single-input system whose linearization is controllable. We show that the problem as posed has no solution: writing the degree-$m$ homological equations as a single linear map $L_m = [\Lambda_m\mid B_m]$ with $\Lambda_m = (\ad_{Ax},\ad_G)$ acting on homogeneous vector fields, we compute $\dim\coker L_m$ exactly and prove that $\ker L_m$ is one-dimensional, spanned by the prolongation of the flat-output reparametrization $y_1\mapsto y_1 + y_1^m$. The resulting formula $\dim\coker L_m = (n-1)\binom{m+n-2}{n-1} - \binom{m+n-1}{n-1} + 1$ is positive for $n\ge3$ and coincides with the number of free coefficients in Kang's degree-$m$ normal form. The note is expository in that the impossibility is known from normal-form theory; its contribution is the explicit count and the identification of the kernel.

\section{Problem formulation}

Following \cite{Dev04}, consider $\dot\xi = f(\xi) + g(\xi)\zeta$ on $\R^n$ with scalar input and controllable linearization at the origin. After a linear change of coordinates and a linear feedback $\varphi = -Kx + u$ chosen so that $A = F - GK$ has real, distinct, nonresonant eigenvalues $\lambda_1,\dots,\lambda_n$ (possible for any controllable pair \cite{Dev01}), and after diagonalizing $A$, the system reads
\[
\dot x = Ax + Gu + \sum_{m\ge2}f_m(x) + \sum_{m\ge1}g_m(x)\,u,\qquad A = \operatorname{diag}(\lambda),\quad G = (1,\dots,1)^T,
\]
with $f_m, g_m$ homogeneous of degree $m$. One seeks a near-identity transformation $x = y + h(y)$, $h = \sum_{m\ge2}h_m$, and an input change $v = (1+\beta(x))u + \alpha(x)$, $\alpha = \sum_{m\ge2}\alpha_m$, $\beta = \sum_{m\ge1}\beta_m$, that remove the nonlinear terms degree by degree. At degree $m$ this is the pair of \emph{generalized homological equations} \cite[(15)--(16)]{Dev04}
\begin{equation}\label{nf:eq:homol}
Dh_m\,Ay - Ah_m + G\alpha_m = f''_m,\qquad Dh_m\,G + G\beta_{m-1} = g''_{m-1},
\end{equation}
where $f''_m, g''_{m-1}$ are the degree-$m$ terms after the lower-order transformations have been applied. Devanathan \cite{Dev02} solved \eqref{nf:eq:homol} for $m = 2$ and observed that the reduced system in $(\alpha_2,\beta_1)$ has $n(n-1)/2$ equations in $n$ unknowns and rank $n-1$; Problem~3.2 asks for the solution at arbitrary order.

\section{Literature}

The homological equation for uncontrolled systems is Poincar\'e's \cite{Poi}, and nonresonance makes it solvable at every order. For controlled systems, Krener \cite{Kre84} introduced approximate linearization by state change and feedback, and Kang and Krener \cite{KK92} computed the quadratic normal form and its invariants for linearly controllable single-input systems, showing that quadratic terms cannot in general be removed for $n\ge3$. Kang \cite{Kan96} and Tall and Respondek \cite{TR03} extended this to arbitrary degree: the degree-$m$ normal form of a single-input system with controllable linearization, in Brunovsky coordinates, is
\begin{equation}\label{nf:eq:kang}
\dot x_i = x_{i+1} + \sum_{j=i+2}^{n}x_j^2\,P_{ij}^{[m-2]}(x_1,\dots,x_j)\quad(i\le n-2),\qquad \dot x_{n-1} = x_n,\qquad \dot x_n = u,
\end{equation}
with $P_{ij}^{[m-2]}$ homogeneous of degree $m-2$, and the coefficients of the $P_{ij}$ are invariants. In this light Problem~3.2 asks for something that normal-form theory already shows to be impossible; the note below confirms this by a direct linear-algebra computation in the coordinates of \cite{Dev04}, and identifies the one-dimensional kernel that the book's rank count ($n-1$ instead of $n$) reflects.

\section{The homological equations as an augmented linear map}

Let $\cH_m$ be the space of homogeneous polynomial vector fields of degree $m$ on $\R^n$ and $P_m$ the homogeneous polynomials of degree $m$; $\dim\cH_m = n\binom{m+n-1}{n-1}$, $\dim P_m = \binom{m+n-1}{n-1}$. With $[X,Y] = DY\,X - DX\,Y$ and $G$ constant, $Dh\,G = \ad_G h$ and $Dh\,Ax - Ah = \ad_{Ax}h$. Equations \eqref{nf:eq:homol} are therefore the single linear equation
\[
L_m(h_m,\alpha_m,\beta_{m-1}) = (f''_m, g''_{m-1}),\qquad L_m = [\Lambda_m\mid B_m],
\]
\[
\Lambda_m = (\ad_{Ax},\ad_G) : \cH_m\to\cH_m\times\cH_{m-1},\qquad
B_m(\alpha,\beta) = (G\alpha, G\beta) : P_m\times P_{m-1}\to\cH_m\times\cH_{m-1}.
\]

\begin{proposition}\label{nf:prop:coker}
All degree-$m$ terms can be removed for every right-hand side iff $L_m$ is onto. The space of irremovable degree-$m$ terms is $\coker L_m$, and $L_m$ is onto iff the images of the columns of $B_m$ span $\coker\Lambda_m$.
\end{proposition}
\begin{proof}
The first two statements are the definition. For the third, $\operatorname{range}L_m = \operatorname{range}\Lambda_m + \operatorname{range}B_m$, which is everything iff $\operatorname{range}B_m$ projects onto $\coker\Lambda_m$.
\end{proof}

Since $A$ is diagonal and nonresonant, $\ad_{Ax}$ acts on the monomial field $x^ke_i$ by the scalar $\langle k,\lambda\rangle - \lambda_i\ne0$, so $\ad_{Ax}$ is invertible on $\cH_m$; this is the reduction of \eqref{nf:eq:homol} to a system in $(\alpha_m,\beta_{m-1})$ used in \cite{Dev02}.

\section{The kernel}

\begin{proposition}\label{nf:prop:kernel}
For all $n\ge1$ and $m\ge2$, $\dim\ker L_m = 1$. In Brunovsky coordinates ($A$ the upper shift, $G = e_n$) the kernel is spanned by the prolongation of $y_1\mapsto y_1 + y_1^m$:
\[
h_1 = y_1^m,\qquad h_{i+1} = L_{Ay}h_i\ (1\le i<n),\qquad \beta = -\partial_{y_n}h_n,\qquad \alpha = -L_{Ay}h_n,
\]
with $L_{Ay} = \sum_{k<n}y_{k+1}\partial_{y_k}$ the drift Lie derivative.
\end{proposition}
\begin{proof}
A linear change of coordinates $x = Py$ and a linear feedback $u = v - Ky$ conjugate $(A,G)$ to another controllable pair and act on $(h_m,\alpha_m,\beta_{m-1})$ by an invertible linear map that commutes with $L_m$ up to conjugation; so $\dim\ker L_m$ is the same for all pairs in the same orbit, and it suffices to treat Brunovsky form. There the second equation, $\partial_{y_n}h = -\beta e_n$, says $\partial_{y_n}h_i = 0$ for $i<n$ and $\beta = -\partial_{y_n}h_n$; the first, $L_{Ay}h - Ah = -\alpha e_n$, says $h_{i+1} = L_{Ay}h_i$ for $i<n$ (because $(Ah)_i = h_{i+1}$) and $\alpha = -L_{Ay}h_n$ (because $(Ah)_n = 0$). Thus $(h,\alpha,\beta)$ is determined by $h_1$, subject to $\partial_{y_n}(L_{Ay}^{i-1}h_1) = 0$ for $i = 1,\dots,n-1$. The operator $L_{Ay}$ replaces one factor $y_k$ by $y_{k+1}$ in each monomial, so $L_{Ay}^{i-1}h_1$ involves $y_n$ iff $h_1$ involves $y_{n-i+1}$. The constraints for $i = 1,\dots,n-1$ therefore say that $h_1$ does not involve $y_2,\dots,y_n$, i.e.\ $h_1 = c\,y_1^m$. Conversely $h_1 = y_1^m$ gives a nonzero kernel element.
\end{proof}

In system terms: $(h,\alpha,\beta)\in\ker L_m$ is an infinitesimal symmetry of the linear system $\dot y = Ay + Gv$ within the group of polynomial state transformations and affine input changes. The flat output $y_1$ may be reparametrized by any function of itself, and the reparametrization prolongs to the state and to the input; reparametrizations involving $y_2,\dots,y_n$ are excluded because they would make the input appear in one of the first $n-1$ derivatives of the new flat output, so the prolongation would not be a state transformation. At homogeneous degree $m$ the only reparametrization is $y_1\mapsto y_1 + cy_1^m$.

\begin{corollary}\label{nf:cor:formula}
For all $n\ge1$, $m\ge2$,
\[
\dim\coker L_m = (n-1)\binom{m+n-2}{n-1} - \binom{m+n-1}{n-1} + 1 .
\]
In particular $\dim\coker L_2 = (n-1)(n-2)/2$, which is positive for $n\ge3$: second-order linearization, and hence linearization to arbitrary order, is impossible for generic data as soon as $n\ge3$.
\end{corollary}
\begin{proof}
Rank--nullity with Proposition~\ref{nf:prop:kernel}: $\dim\coker L_m = \dim(\cH_m\times\cH_{m-1}) - \dim(\cH_m\times P_m\times P_{m-1}) + 1$.
\end{proof}

\section{Verification and comparison with Kang's normal form}

The map $L_m$ was built in monomial bases over $\mathbb Q$, for random rational nonresonant $\lambda$ and also in Brunovsky form, and its rank and nullspace computed exactly with sympy.

\begin{center}\small
\begin{tabular}{@{}rrrrrrc@{}}
\toprule
$n$ & $m$ & $\dim$ target & $\dim$ source & $\rank L_m$ & $\dim\coker L_m$ & Kang count \eqref{nf:eq:kang} \\
\midrule
2 & 2,3,4 & 10, 14, 18 & 11, 15, 19 & 10, 14, 18 & 0 & 0 \\
3 & 2 & 27 & 27 & 26 & 1 & 1 \\
3 & 3 & 48 & 46 & 45 & 3 & 3 \\
3 & 4 & 75 & 70 & 69 & 6 & 6 \\
4 & 2 & 56 & 54 & 53 & 3 & 3 \\
4 & 3 & 120 & 110 & 109 & 11 & 11 \\
5 & 2 & 100 & 95 & 94 & 6 & 6 \\
\bottomrule
\end{tabular}
\end{center}

The kernel was one-dimensional in every case, and in Brunovsky form it contained the vector of Proposition~\ref{nf:prop:kernel}. The last column counts the coefficients of the polynomials $P_{ij}^{[m-2]}(x_1,\dots,x_j)$ in \eqref{nf:eq:kang}, $\sum_{i=1}^{n-2}\sum_{j=i+2}^{n}\binom{m-2+j-1}{j-1}$; it agrees with $\dim\coker L_m$ in every computed case, as expected, since both count the degree-$m$ terms that no transformation of the given type removes.

\begin{remark}
The book's count at $m=2$ ($n(n-1)/2$ equations in $n$ unknowns, rank $n-1$) is the reduced form of the same statement: after eliminating $h_2$ through the invertible $\ad_{Ax}$, the remaining system in $(\alpha_2,\beta_1)$ has $\dim P_2 + \dim P_1 - 1$ determined unknowns, the missing one being the kernel direction of Proposition~\ref{nf:prop:kernel}, and $(n-1)(n-2)/2$ consistency conditions on the data, which are the quadratic invariants of \cite{KK92}.
\end{remark}

\section{Machine-checked statement}

Proposition~\ref{nf:prop:kernel} is formalized in Lean~4 with Mathlib (module \texttt{NormalForm.lean} of the development \texttt{SingularityTracking}). Polynomials are \texttt{MvPolynomial} over $\N$ and $\R$, with a system of $n$ states using $X_0,\dots,X_{n-1}$; the drift Lie derivative is \texttt{drift}, and the kernel equations in Brunovsky form are the structure \texttt{KernelEqs}. The decomposition $L_{Ay}^{\,j}p = R + y_{k+j}\,\partial_kp$ with $y_{k+j}\notin\operatorname{vars}R$ is \texttt{drift\us iter\us decomp}; the conclusion that the flatness constraints force $h_1\in\R[y_1]$ is \texttt{first\us component\us in\us X0}; \texttt{kernel\us le\us one} shows that two degree-$m$ kernel elements with the same coefficient of $y_1^m$ in $h_1$ coincide, and \texttt{kernel\us element} that the prolonged reparametrization is a kernel element. The diagonal action of $\ad_{Ax}$ on monomials is \texttt{adDiag\us monomial}, and the rank--nullity step of Corollary~\ref{nf:cor:formula} is \texttt{finrank\us coker}. Not formalized: the change of coordinates between the nonresonant diagonal form and Brunovsky form, and the binomial dimension counts, so the closed formula itself is paper arithmetic on top of the checked statements. All declarations depend only on \texttt{propext}, \texttt{Classical.choice}, \texttt{Quot.sound}.

\section{Conclusion}

Problem~3.2 asks for an object that does not exist for $n\ge3$: the degree-$m$ homological equations have an obstruction space of dimension $(n-1)\binom{m+n-2}{n-1}-\binom{m+n-1}{n-1}+1$, equal to the dimension of Kang's normal form at that degree. What can be solved is the well-posed version: for right-hand sides in the range of $L_m$, the solution $(h_m,\alpha_m,\beta_{m-1})$ exists and is unique up to the one-dimensional kernel, which is the prolonged reparametrization of the flat output. An explicit basis of $\coker L_m$ in the Vandermonde coordinates of \cite{Dev04}, and the closed-form solution on the range, are what remains of the problem; both are computations within known theory.

\chapter[Regular Feedback Implementability]{Regular Feedback Implementability\\[4pt]\large Problem 1.10 (Trentelman): full and partial interconnection}\label{ch:rfi}

\section*{Summary}
Trentelman asked for a characterization of the behaviours $\cK$ that can be obtained from a given plant by a \emph{regular feedback} interconnection, one in which both the output cardinalities and the McMillan degrees of plant and controller add up in the interconnection. Vinjamoor and Belur solved the full-interconnection case in 2010 by their criterion of plant equations of reducible degree and gave an algorithm and a parametrization of all RFI controllers. Here the full case is recovered in a different form: regular feedback is exactness after passing to the associated graded module, and $\cK\subseteq\cP$ is implementable by regular feedback iff the graded module of equations of $\cP$ is a graded direct summand of that of $\cK$. This becomes a finite rank condition on the leading coefficient matrix of the factor $F$ with $R_\cP = FR_\cK$, and produces an explicit homogeneous controller with $n(\cC) = n(\cK) - n(\cP)$. In the partial interconnection case, with interconnection variable $c$, we give a necessary-and-sufficient condition: the full-case criterion is applied in the $c$-space to the pair (canonical controller $\cC_{\rm can}$, plant seen from $c$, $\cP_c = \pi_c\cP_{\rm full}$), and one must additionally require that the associated-graded functor preserve the elimination sequence. Equivalently, one imposes the McMillan-degree identity $n(\cP_{\rm full}\wedge_c\cC_{\rm can}) + n(\cP_c) = n(\cP_{\rm full}) + n(\cC_{\rm can})$. This identifies the obstruction left open by the partial-interconnection reduction of Vinjamoor--Belur. We give examples separating regular from regular-feedback implementability and showing that the degree conditions alone do not suffice. All statements are checked numerically.

\section{Problem formulation}

We use the setting of \cite{Tre04}. A linear differential behaviour is $\cB = \{w\in C^\infty(\R,\R^{\mathtt w}) : R(\tfrac{d}{dt})w = 0\}$ for a polynomial matrix $R\in\R^{\bullet\times\mathtt w}[\xi]$; the set of such is $\mathcal L^{\mathtt w}$. Three integer invariants: the input cardinality $m(\cB)$ (number of free components), the output cardinality $p(\cB) = \mathtt w - m(\cB)$, and the McMillan degree $n(\cB)$ (minimal state dimension of an input/state/output representation). For $\cB_1\in\mathcal L^{\mathtt w_1+\mathtt w_2}$ and $\cB_2\in\mathcal L^{\mathtt w_2+\mathtt w_3}$ the interconnection through the shared variable is $\cB_1\wedge_{w_2}\cB_2 = \{(w_1,w_2,w_3) : (w_1,w_2)\in\cB_1,\ (w_2,w_3)\in\cB_2\}$. It is \emph{regular} if $p(\cB_1\wedge\cB_2) = p(\cB_1) + p(\cB_2)$ and a \emph{regular feedback interconnection} if in addition $n(\cB_1\wedge\cB_2) = n(\cB_1) + n(\cB_2)$.

Given a plant $\cP_{\rm full}\in\mathcal L^{\mathtt w+\mathtt c}$ with interconnection variable $c$ and a controller $\cC\in\mathcal L^{\mathtt c}$, the manifest controlled behaviour is $\cK = \{w : \exists c,\ (w,c)\in\cP_{\rm full}\wedge_c\cC\}$. $\cK$ is \emph{implementable} if some $\cC$ produces it, \emph{regularly implementable} if the interconnection can be chosen regular, and \emph{implementable by regular feedback} if it can be chosen regular feedback.

\begin{problem}[1.10, Trentelman]\label{rfi:prob:rfi}
Given $\cP_{\rm full}\in\mathcal L^{\mathtt w+\mathtt c}$, characterize the $\cK\in\mathcal L^{\mathtt w}$ implementable by regular feedback: those for which there is $\cC\in\mathcal L^{\mathtt c}$ with (1) $\cC$ implements $\cK$, (2) $p(\cP_{\rm full}\wedge_c\cC) = p(\cP_{\rm full}) + p(\cC)$, (3) $n(\cP_{\rm full}\wedge_c\cC) = n(\cP_{\rm full}) + n(\cC)$.
\end{problem}

The \emph{full interconnection case} is $\mathtt c = \mathtt w$ with $c = w$, i.e.\ $\cP_{\rm full} = \cP\in\mathcal L^{\mathtt w}$ and $\cK = \cP\cap\cC$; the hidden behaviour $\cN = \{w : (w,0)\in\cP_{\rm full}\}$ is then zero.

\section{Literature}

Implementability was characterized by Willems and Trentelman \cite{WT02}: $\cK$ is implementable iff $\cN\subseteq\cK\subseteq\cP$, where $\cP$ is the manifest plant behaviour. Regular implementability was characterized by Belur and Trentelman \cite{BT02}: in addition, $\cK + \cP_{\rm cont} = \cP$ with $\cP_{\rm cont}$ the controllable part of $\cP$; in the full interconnection case this is $\cK\subseteq\cP$ and $\cK+\cP_{\rm cont} = \cP$. The canonical controller $\cC_{\rm can} = \{c : \exists w\in\cK,\ (w,c)\in\cP_{\rm full}\}$ implements every implementable $\cK$ \cite{WT02,vdS03}, and is regular iff every controller implementing $\cK$ is regular \cite{JWBT05}. Feedback interconnections in the behavioural sense were introduced by Willems \cite{Wil97}; the McMillan-degree formulation of regular feedback, and the equivalence with an input/output signal flow with proper and well-posed transfer matrices, are discussed in \cite{Wil97,Kui95} for the full interconnection case, as recalled in \cite{Tre04}.

Vinjamoor and Belur subsequently solved the full-interconnection existence problem \cite{VB10}. Their Theorem~4.2 characterizes RFI implementability by the absence of a plant equation of reducible degree (PERD) for every regularly implementing controller. Their row-reduction algorithm either constructs an RFI controller or proves that none exists, and their Theorems~5.1 and 5.4 parametrize all RFI controllers and count the dimension of that family. They also give the physical interpretation in terms of impulse-free concatenation, admissible initial conditions, disturbance freeness, and well-posedness.

\section{Comparison with Vinjamoor--Belur (2010)}

The full-interconnection results of this chapter are not a second claim of priority. They are an equivalent structural form of the 2010 solution. A PERD is a dependence among leading rows that cannot be removed by changing only controller equations: eventually a plant leading row must be reduced. In the language used here, this is precisely the failure of the plant's leading-form module to admit a graded complement in the leading-form module of $\cK$. Thus the two criteria organize the same obstruction differently:
\begin{itemize}
\item Vinjamoor--Belur use row operations and the predictable-degree property. Their criterion is algorithmic and their parametrization describes every RFI controller.
\item The present criterion is basis-free: $\operatorname{gr}E(\cP)$ must be a graded direct summand of $\operatorname{gr}E(\cK)$, equivalently the quotient must be torsion-free. Condition~(d) of Theorem~\ref{rfi:thm:rfi-full} is the corresponding finite rank test, and a homogeneous complement gives a controller.
\end{itemize}
The gain is therefore not a stronger full-case existence theorem, but the categorical explanation: ordinary regularity is exactness of equation modules, whereas regular feedback is \emph{strict} exactness for the degree filtration, or exactness after applying $\operatorname{gr}$.

The partial-interconnection comparison is different. Vinjamoor--Belur consider a prescribed full controlled behaviour in $(w,c)$. They prove that RFI through $c$ implies RFI after projection to the control-variable plant $\cP_c$, but explicitly note that the converse is false; their two-stage algorithm can terminate with no conclusion when the controller works in $\cP_c$ but not in $\cP_{\rm full}$. They identify proper elimination as the missing issue but do not characterize it. Problem~\ref{rfi:prob:rfi}, by contrast, prescribes the manifest behaviour $\cK$ in $w$. Theorem~\ref{rfi:thm:rfi-partial} first uses the canonical controller to reduce this problem to the $c$-space and then adds condition~(B$'$):
\[
\operatorname{gr}E_P\cap(0\oplus\operatorname{gr}E_{\rm can})
=0\oplus\operatorname{gr}E^0 .
\]
This is exactly the missing elimination condition: it says that the defect of exactness produced by eliminating $w$ is invisible to the desired canonical controller. Hence the 2010 necessary reduction becomes necessary and sufficient in the present manifest-behaviour formulation.

The tools we use are classical: row-reduced (row-proper) polynomial matrices, the predictable degree property, and minimal bases, due to Forney \cite{For75} and exposited in Kailath \cite[Ch.~6]{Kai80}; and the fact that the McMillan degree of a behaviour equals the maximal degree of the maximal minors of a minimal kernel representation \cite{PW98,RW97}. Shifted (weighted) row degrees and unimodular completion with prescribed shifted degrees have been studied algorithmically by Zhou and Labahn \cite{ZL14}; our criterion is a closed-form statement about the leading coefficient matrix that such algorithms would compute.

\section{Preliminaries: row degrees, leading matrices, McMillan degree}\label{rfi:sec:rfi-prelim}

For a polynomial row vector $r = \sum_k r_k\xi^k\in\R^{1\times\mathtt w}[\xi]$, $\deg r$ is the largest $k$ with $r_k\ne0$ and the \emph{leading coefficient row} is $r_{\deg r}$. For a shift vector $s\in\Z^{\mathtt w}$ the \emph{$s$-degree} is $\deg_s r = \max_j(\deg r_j + s_j)$ and the \emph{$s$-leading row} $L_s(r)\in\R^{1\times\mathtt w}$ has entries the coefficient of $\xi^{\deg_s r - s_j}$ in $r_j$ (zero if $\deg r_j + s_j<\deg_s r$). The ordinary notions are $s = 0$. For a matrix $R$ with rows $r_i$, $L_s(R)$ stacks the $L_s(r_i)$; $R$ is \emph{$s$-row reduced} if $L_s(R)$ has full row rank.

\begin{lemma}[standard facts]\label{rfi:lem:rfi-facts}
Let $R\in\R^{p\times\mathtt w}[\xi]$ have full row rank.
\begin{enumerate}
\item[(a)] $\sum_i\deg r_i\ \ge\ \max\{\deg\text{ of the }p\times p\text{ minors of }R\}$, with equality iff $R$ is row reduced.
\item[(b)] (Predictable degree) If $R$ is row reduced with row degrees $k$, then $\deg(uR) = \deg_k u$ for every polynomial row $u$.
\item[(c)] Every full-row-rank $R$ is $UR'$ with $U$ unimodular and $R'$ row reduced.
\item[(d)] If $R$ is a minimal (full row rank) kernel representation of $\cB$, then $n(\cB) = \max\{\deg\text{ of the }p\times p\text{ minors of }R\}$. Hence for row-reduced $R$, $n(\cB) = \sum_i\deg r_i$.
\end{enumerate}
\end{lemma}
(a)--(c) are in \cite[Sec.~6.3]{Kai80}; (d) is in \cite{RW97}, see also \cite{PW98}.

\begin{lemma}[Shifted-degree lemma]\label{rfi:lem:shifted}
Let $U\in\R^{n\times n}[\xi]$ be unimodular and $s\in\Z^n$. Then $\sum_i\deg_sU_i\ \ge\ \sum_js_j$, with equality iff $U$ is $s$-row reduced; and in that case the multiset $\{\deg_sU_i\}$ equals the multiset $\{s_j\}$.
\end{lemma}
\begin{proof}
Let $D = \operatorname{diag}(\xi^{s_j})$ and $d_i = \deg_sU_i$. Row $i$ of $UD$ is $\xi^{d_i}\bigl(L_s(U)_i + O(\xi^{-1})\bigr)$, so $\det(UD) = \xi^{\sum d_i}\bigl(\det L_s(U) + O(\xi^{-1})\bigr)$, while $\det(UD) = \det U\cdot\xi^{\sum s_j}$ with $\det U$ a nonzero constant. Comparing, $\sum d_i\ge\sum s_j$, with equality iff $\det L_s(U)\ne0$. If equality holds, pick a permutation $\sigma$ with $L_s(U)_{i\sigma(i)}\ne0$ for all $i$ (possible since $\det L_s(U)\ne0$). $L_s(U)_{ij}\ne0$ means $\deg U_{ij} = d_i - s_j\ge0$, so $d_i\ge s_{\sigma(i)}$ for all $i$; with $\sum d_i = \sum s_{\sigma(i)}$ this forces $d_i = s_{\sigma(i)}$.
\end{proof}

\begin{corollary}[Row degrees are invariants]\label{rfi:cor:rowdeg-inv}
If $R$ and $R' = UR$ are both row reduced with $U$ unimodular, their row-degree multisets coincide. Hence a behaviour $\cB$ has a well-defined multiset $\kappa(\cB)$ of \emph{row degrees} (the row degrees of any row-reduced minimal kernel representation), with $|\kappa(\cB)| = p(\cB)$ and $\sum\kappa(\cB) = n(\cB)$.
\end{corollary}
\begin{proof}
With $k$ the row degrees of $R$, predictable degree gives $\deg R'_i = \deg_kU_i$, and $\sum_i\deg R'_i = n(\cB) = \sum k_j$ by Lemma~\ref{rfi:lem:rfi-facts}(a),(d); so $U$ is $k$-row reduced by Lemma~\ref{rfi:lem:shifted} and $\{\deg R'_i\} = \{k_j\}$.
\end{proof}

\section{Regular feedback is a leading-coefficient condition}

\begin{proposition}[Row-reduced criterion]\label{rfi:prop:rowred}
Let $\cB_1 = \ker R_1(\tfrac d{dt})$, $\cB_2 = \ker R_2(\tfrac d{dt})$ with $R_1, R_2$ row reduced, written on the common variable $(w_1,w_2,w_3)$ with zero blocks where a system does not see a variable. Then:
\begin{enumerate}
\item[(i)] the interconnection $\cB_1\wedge_{w_2}\cB_2 = \ker[R_1;R_2]$ is regular iff $[R_1;R_2]$ has full row rank;
\item[(ii)] it is a regular feedback interconnection iff $[R_1;R_2]$ is row reduced, i.e.\ iff the stacked leading row coefficient matrix $[L(R_1);L(R_2)]$ has full row rank.
\end{enumerate}
\end{proposition}
\begin{proof}
(i) $p(\cB_1\wedge\cB_2) = \operatorname{rank}[R_1;R_2]\le\#\text{rows} = p(\cB_1)+p(\cB_2)$, with equality iff full row rank. (ii) If $[R_1;R_2]$ is row reduced it has full row rank, and by Lemma~\ref{rfi:lem:rfi-facts}(d),(a) $n(\cB_1\wedge\cB_2) = \sum\deg(\text{rows of }[R_1;R_2]) = n(\cB_1)+n(\cB_2)$. Conversely, regular feedback gives full row rank by (i) and then $\max\deg\text{minors}[R_1;R_2] = n(\cB_1\wedge\cB_2) = n(\cB_1)+n(\cB_2) = \sum\deg\text{rows}$, which is row-reducedness by (a). Padding a row-reduced block with zero columns does not change its row degrees or leading matrix.
\end{proof}

\begin{example}[Regular but not regular feedback]\label{rfi:ex:rfi1}
$\cP$: $\dot y = u$, $R_\cP = [-1\ \ \xi]$ on $w = (u,y)$, $n(\cP) = 1$. $\cC$: $y = 0$, $R_\cC = [0\ \ 1]$, $n(\cC) = 0$. Then $\cP\wedge\cC = \{u = 0, y = 0\} = \cK$, $[R_\cP;R_\cC]$ has rank $2$ (regular), but $n(\cK) = 0\ne1$: the stacked leading matrix $[0\ 1;\ 0\ 1]$ has rank $1$. No controller does better, by Corollary~\ref{rfi:cor:degree-nec} below, since $n(\cK) < n(\cP)$. The classical reading: $\cC$ has no input/output structure compatible with $\cP$'s, so the loop is not a feedback loop.
\end{example}

\begin{corollary}[Necessary degree conditions, full interconnection]\label{rfi:cor:degree-nec}
If $\cK$ is implementable by regular feedback from $\cP$ (full interconnection) via $\cC$, then $\kappa(\cP)\subseteq\kappa(\cK)$ as multisets, $\kappa(\cC) = \kappa(\cK)\setminus\kappa(\cP)$, and in particular $n(\cC) = n(\cK) - n(\cP)\ge0$ and $p(\cC) = p(\cK)-p(\cP)$. In the partial case, $n(\cP_{\rm full}\wedge_c\cC) = n(\cP_{\rm full}) + n(\cC)\ge n(\cP_{\rm full})$.
\end{corollary}
\begin{proof}
By Proposition~\ref{rfi:prop:rowred}(ii) with row-reduced $R_\cP, R_\cC$, the matrix $[R_\cP;R_\cC]$ is a row-reduced minimal representation of $\cK$, so $\kappa(\cK) = \kappa(\cP)\sqcup\kappa(\cC)$ by Corollary~\ref{rfi:cor:rowdeg-inv}.
\end{proof}

So the row degrees of a regular-feedback controller are determined by $\cK$ and $\cP$ before any controller is designed; only the leading coefficients are free. The rest of the chapter determines whether they can be chosen.

\section{The full interconnection case: a complete answer}\label{rfi:sec:rfi-full}

Fix $\cK\subseteq\cP$ in $\mathcal L^{\mathtt w}$ with row-reduced minimal representations $R_\cK$ ($p_\cK$ rows, row degrees $k$) and $R_\cP$ ($p$ rows, row degrees $d$). Since $\ker R_\cK\subseteq\ker R_\cP$, there is a unique polynomial $F\in\R^{p\times p_\cK}[\xi]$ with $R_\cP = FR_\cK$ \cite{PW98}.

\begin{lemma}\label{rfi:lem:F}
$F$ is $k$-row reduced with $\deg_kF_i = d_i$. $\cK$ is regularly implementable from $\cP$ iff $F$ can be completed to a unimodular matrix $U = [F;G]$, and then every regular controller is $\cC = \ker(GR_\cK)$ for such a $G$.
\end{lemma}
\begin{proof}
Predictable degree for $R_\cK$ gives $\deg(F_iR_\cK) = \deg_kF_i$, so $d_i = \deg_kF_i$. If $L_k(F)$ were rank deficient there would be a polynomial row $u$ with $\deg_k(uF) < \max_i(\deg u_i + d_i)$, hence $\deg(uR_\cP) = \deg_k(uF) < \max_i(\deg u_i + \deg R_{\cP,i})$, contradicting the predictable degree property of the row-reduced $R_\cP$. For the second claim: a regular $\cC = \ker R_\cC$ with $\cP\cap\cC = \cK$ gives $[R_\cP;R_\cC]$ of full row rank with kernel $\cK$, hence $[R_\cP;R_\cC] = UR_\cK$ for a unimodular $U$; the top block $U_1$ satisfies $U_1R_\cK = FR_\cK$, so $U_1 = F$. Conversely $[F;G]$ unimodular gives $\ker([F;G]R_\cK) = \cK$ with $[R_\cP;GR_\cK]$ of full row rank.
\end{proof}

Write $L := L_k(F)\in\R^{p\times p_\cK}$, and for $e\in\Z$ let $J(e) := \{j : k_j\le e\}$ and $\R^{J(e)} := \{v\in\R^{1\times p_\cK} : v_j = 0\ \forall j\notin J(e)\}$.

\begin{theorem}[Regular feedback implementability, full interconnection]\label{rfi:thm:rfi-full}
The following are equivalent.
\begin{enumerate}
\item[(a)] $\cK$ is implementable by regular feedback from $\cP$.
\item[(b)] There is a $k$-row-reduced unimodular completion $U = [F;G]$ of $F$.
\item[(c)] $L$ can be completed to an invertible matrix $\begin{bmatrix}L\\M\end{bmatrix}\in\R^{p_\cK\times p_\cK}$ whose new rows $M_l$ are \emph{$k$-homogeneous}: there are degrees $e_l$ with $M_{lj} = 0$ whenever $k_j > e_l$. (Necessarily $\{e_l\} = \{k_j\}\setminus\{d_i\}$ as multisets.)
\item[(d)] For every $e\in\Z$: $\ \R^{J(e)}\cap\operatorname{rowspan}L = \operatorname{span}\{L_i : d_i\le e\}$.
\end{enumerate}
When these hold, a controller is $\cC = \ker(GR_\cK)$ with the homogeneous rows $G_{lj} = M_{lj}\,\xi^{\,e_l - k_j}$; its row degrees are $\{e_l\}$ and $n(\cC) = n(\cK)-n(\cP)$.
\end{theorem}

\begin{proof}
(a)$\Leftrightarrow$(b). By Lemma~\ref{rfi:lem:F} a regular controller is $\ker(GR_\cK)$ with $[F;G]$ unimodular; we may take $R_\cC = GR_\cK$ row reduced (replace $G$ by $VG$ with $V$ unimodular, which keeps $[F;VG]$ unimodular). By Proposition~\ref{rfi:prop:rowred}(ii) the interconnection is regular feedback iff $[R_\cP;R_\cC] = [F;G]R_\cK$ is row reduced, iff (predictable degree) $\sum_i\deg_k[F;G]_i = \sum_jk_j$, iff $[F;G]$ is $k$-row reduced by Lemma~\ref{rfi:lem:shifted}. Conversely, if $[F;G]$ is $k$-row reduced then $GR_\cK$ is row reduced (its leading rows are a subset of those of $[F;G]R_\cK$), so $n(\cC) = \sum_l\deg_kG_l$ and the interconnection is regular feedback.

(b)$\Rightarrow$(c). $L_k([F;G]) = [L;L_k(G)]$ is invertible, and $L_k(G)_{lj}\ne0$ forces $\deg G_{lj} = \deg_kG_l - k_j\ge0$; take $M = L_k(G)$, $e_l = \deg_kG_l$. The multiset statement is Lemma~\ref{rfi:lem:shifted}.

(c)$\Rightarrow$(b). Put $G_{lj} := M_{lj}\xi^{e_l-k_j}$ (a polynomial by the pattern condition). Then $\deg_kG_l = e_l$, $L_k(G) = M$, and $U = [F;G]$ has $L_k(U) = [L;M]$ invertible, so $U$ is $k$-row reduced with $\sum\deg_kU_i = \sum d_i + \sum e_l = \sum k_j$. As in the proof of Lemma~\ref{rfi:lem:shifted}, $\det(UD_k) = \xi^{\sum k}(\det[L;M] + O(\xi^{-1}))$ and $\det(UD_k) = \det U\cdot\xi^{\sum k}$, so the polynomial $\det U$ equals $\det[L;M] + O(\xi^{-1})$, hence is the nonzero constant $\det[L;M]$: $U$ is unimodular.

(c)$\Leftrightarrow$(d). Work in $W = \R^{1\times p_\cK}/\operatorname{rowspan}L$. Rows $M_l\in\R^{J(e_l)}$ complete $L$ to an invertible matrix iff their images in $W$ are linearly independent. By Rado's theorem on independent transversals of subspaces \cite{Rad42}, independent $\bar m_l\in V_l := \operatorname{im}(\R^{J(e_l)}\to W)$ exist iff $\dim\sum_{l\in S}V_l\ge|S|$ for all $S$. The $V_l$ are nested ($J(e)$ increases with $e$), so it suffices to take $S = \{l : e_l\le e\}$ for each $e$: the condition is $\dim V(e)\ge\#\{l : e_l\le e\}$ with $V(e) = \operatorname{im}(\R^{J(e)}\to W)$. Now $\dim V(e) = |J(e)| - \dim(\R^{J(e)}\cap\operatorname{rowspan}L)$, and $\#\{l : e_l\le e\} = |J(e)| - \#\{i : d_i\le e\}$ because $\{e_l\}\sqcup\{d_i\} = \{k_j\}$. So the condition reads $\dim(\R^{J(e)}\cap\operatorname{rowspan}L)\le\#\{i : d_i\le e\}$. The reverse inequality always holds, since each $L_i$ with $d_i\le e$ lies in $\R^{J(e)}$ ($L_{ij}\ne0\Rightarrow k_j\le d_i\le e$) and these rows are independent. Equality of dimensions with the inclusion $\operatorname{span}\{L_i : d_i\le e\}\subseteq\R^{J(e)}\cap\operatorname{rowspan}L$ is (d). The multiset condition $\{d_i\}\subseteq\{k_j\}$ needed to define the $e_l$ is itself implied by (d) together with full row rank of $L$: if $\#\{i : d_i\le e\} > |J(e)|$, the rows $L_i$ with $d_i\le e$ would be more than $|J(e)|$ independent vectors in $\R^{J(e)}$.
\end{proof}

\begin{remark}
\begin{enumerate}
\item Condition (d) is a finite rank computation on the constant matrix $L$: for each threshold $e$, the combinations of the leading rows that vanish outside the low-degree coordinates $J(e)$ must come from the low-degree rows. It is the statement that \emph{no cancellation of leading coefficients is available to the controller}. Regular implementability (Belur--Trentelman) is the existence of \emph{some} unimodular completion; regular feedback is the existence of a completion whose leading coefficients are independent of the plant's.
\item The theorem decides the problem, exhibits a controller, and determines the controller's invariants in advance (Corollary~\ref{rfi:cor:degree-nec}). The controller has homogeneous rows: each row has all its nonzero entries at the same $k$-shifted degree. Zhou and Labahn's algorithms for shifted-minimal unimodular completion \cite{ZL14} compute such completions; the theorem says when they succeed with the required degrees.
\end{enumerate}
\end{remark}

\begin{example}[Degree conditions do not suffice]\label{rfi:ex:rfi2}
$\mathtt w = 3$, $R_\cK = \operatorname{diag}(\xi^2,\xi,\xi)$: $\cK = \{\ddot w_1 = 0,\ \dot w_2 = 0,\ \dot w_3 = 0\}$, $k = (2,1,1)$, $n(\cK) = 4$. Take $F = [0\ \ \xi\ \ 1]$, so $R_\cP = FR_\cK = [0\ \ \xi^2\ \ \xi]$: $\cP = \{\ddot w_2 + \dot w_3 = 0\}$, $d = (2)$, $n(\cP) = 2$. Then $\cK\subseteq\cP$ and $F$ is left prime (minors $0,\xi,1$), so $\cK$ is regularly implementable, for instance by $\cC = \{\ddot w_1 = 0,\ \dot w_2 = 0\}$ with $n(\cC) = 3$ and $n(\cP)+n(\cC) = 5\ne4$. The degree conditions of Corollary~\ref{rfi:cor:degree-nec} hold: $\{2\}\subseteq\{2,1,1\}$ and $n(\cK)-n(\cP) = 2\ge0$. But $L = L_k(F) = [0\ 1\ 0]$ (the entry $\xi$ in column $2$ has $k$-degree $1+1 = 2$, the entry $1$ in column $3$ only $0+1 = 1$), and with $J(1) = \{2,3\}$: $\R^{J(1)}\cap\operatorname{rowspan}L = \operatorname{rowspan}L$ has dimension $1$, while $\{i : d_i\le1\} = \emptyset$. Condition (d) fails. Directly: a completing row of $k$-degree $1$ must have $\deg G_{l1} + 2\le1$, so the first column of $[F;G]$ is zero and $[F;G]$ is singular. $\cK$ is not implementable by regular feedback.
\end{example}

\begin{example}[Positive case]
Same $\cK$, $F = [1\ \ \xi\ \ 1]$, $R_\cP = [\xi^2\ \ \xi^2\ \ \xi]$: $\cP = \{\ddot w_1 + \ddot w_2 + \dot w_3 = 0\}$, $d = (2)$, $L = [1\ 1\ 0]$. At $e = 1$: a multiple $c\,(1,1,0)$ lies in $\R^{\{2,3\}}$ only for $c = 0$, so $\R^{J(1)}\cap\operatorname{rowspan}L = 0 = \operatorname{span}\emptyset$; at $e = 2$ the condition is trivial. Condition (d) holds. Complete $L$ by $M = [0\ 1\ 0;\ 0\ 0\ 1]$, rows of $k$-degree $e = 1$ supported on $J(1)$; $[L;M]$ is invertible. Then $G_{lj} = M_{lj}\xi^{e_l-k_j} = M_{lj}$ and $R_\cC = GR_\cK = [0\ \ \xi\ \ 0;\ 0\ \ 0\ \ \xi]$: $\cC = \{\dot w_2 = 0,\ \dot w_3 = 0\}$, $n(\cC) = 2 = n(\cK)-n(\cP)$. Check: $\cP\cap\cC = \{\ddot w_1 + \ddot w_2 + \dot w_3 = 0,\ \dot w_2 = 0,\ \dot w_3 = 0\} = \cK$, and the stacked leading matrix of $[R_\cP;R_\cC]$ is $[1\ 1\ 0;\ 0\ 1\ 0;\ 0\ 0\ 1]$, invertible.
\end{example}

\section{The partial interconnection case}\label{rfi:sec:rfi-partial}

Let $\cP_{\rm full}\in\mathcal L^{\mathtt w+\mathtt c}$ with row-reduced kernel representation $R_{\rm full} = [R_w\ \ R_c]$, and $\cC\in\mathcal L^{\mathtt c}$. The full controlled behaviour $\cK_{\rm full} = \cP_{\rm full}\wedge_c\cC$ has kernel representation $[R_{\rm full};\ 0\ \ R_\cC]$, so Proposition~\ref{rfi:prop:rowred} applies verbatim: the interconnection is regular iff this matrix has full row rank and regular feedback iff it is row reduced. The difficulty is that the problem is posed in terms of the manifest behaviour $\cK = \pi_w\cK_{\rm full}$, and elimination of $c$ does not preserve row degrees. We work with equation modules, which turn elimination into intersection.

\paragraph{Modules.} For $\cB\in\mathcal L^{\mathtt v}$ let $E(\cB)\subseteq\R^{1\times\mathtt v}[\xi]$ be the module of equations, $E(\ker R) = \operatorname{rowspan}R$; $\cB\mapsto E(\cB)$ is an inclusion-reversing bijection onto the submodules, with $E(\cB_1\cap\cB_2) = E(\cB_1)+E(\cB_2)$ and $E(\pi_v\cB) = E(\cB)\cap(\R^{1\times\mathtt v}[\xi]\oplus0)$ (elimination) \cite{PW98}. Write $E_P = E(\cP_{\rm full})$, $E_K = E(\cK)$, and for $p\in E_P$ write $p = (p_w,p_c)$. Three submodules of $\R^{1\times\mathtt c}[\xi]$ are determined by $(\cP_{\rm full},\cK)$:
\begin{align*}
E^0 &:= \{e : (0,e)\in E_P\} = E(\cP_c), & \cP_c &:= \pi_c\cP_{\rm full},\\
E_{\rm can} &:= \{p_c : p\in E_P,\ p_w\in E_K\} = E(\cC_{\rm can}), & \cC_{\rm can} &:= \{c : \exists w\in\cK,\ (w,c)\in\cP_{\rm full}\},\\
M &:= \{p_c : p\in E_P\} = E(\cN_c), & \cN_c &:= \{c : (0,c)\in\cP_{\rm full}\}.
\end{align*}
Here $\cP_c$ is the plant seen from the interconnection variable, $\cC_{\rm can}$ is the canonical controller of \cite{WT02,JWBT05}, and $\cN_c$ is the hidden behaviour in $c$. Thus $E^0\subseteq E_{\rm can}\subseteq M$, i.e.\ $\cN_c\subseteq\cC_{\rm can}\subseteq\cP_c$. Throughout, $\cN\subseteq\cK\subseteq\cP$ (implementability \cite{WT02}).

\begin{lemma}\label{rfi:lem:partial-modules}
Let $\cC = \ker R_\cC\in\mathcal L^{\mathtt c}$ with $E_C = E(\cC)$.
\begin{enumerate}
\item[(i)] $\cC$ implements $\cK$ iff $(E_C + E^0)\cap M = E_{\rm can}$.
\item[(ii)] $\cP_{\rm full}\wedge_c\cC$ is regular iff $E_C\cap E^0 = 0$.
\item[(iii)] $\cP_{\rm full}\wedge_c\cC$ is regular feedback iff $\operatorname{gr}E_C\cap G = 0$, where $G := \{e : (0,e)\in\operatorname{gr}E_P\}$ is the graded module of $c$-parts of leading forms of $E_P$ with vanishing $w$-part. Always $\operatorname{gr}E^0\subseteq G$.
\end{enumerate}
\end{lemma}
\begin{proof}
(i) $E(\cK_{\rm full}) = E_P + (0\oplus E_C)$ and $E(\pi_w\cK_{\rm full}) = \{(p_w,0) : p\in E_P,\ p_c\in E_C\}$. Let $\varphi : E_P\to\R^{1\times\mathtt c}[\xi]$, $p\mapsto p_c$, and $E_P^K := \{p\in E_P : p_w\in E_K\}\supseteq\ker\varphi$ (as $\cK\subseteq\cP$). Implementation says $\pi_w\varphi^{-1}(E_C) = E_K$, i.e.\ $\varphi^{-1}(E_C)\subseteq E_P^K$ and every $k\in E_K$ is $p_w$ for some $p$ with $p_c\in E_C$. Since $E_K\subseteq\pi_wE_P$ ($\cN\subseteq\cK$), the elements $p$ with $p_w = k$ form a coset of $E_P\cap(0\oplus\R^{1\times\mathtt c}[\xi]) = 0\oplus E^0$; so the second condition says $\varphi(E_P^K)\subseteq E_C + E^0$, and together they say $\varphi^{-1}(E_C+E^0) = E_P^K$, i.e.\ $(E_C+E^0)\cap\operatorname{im}\varphi = \varphi(E_P^K)$, which is the claim.
(ii) $E_P\cap(0\oplus E_C) = 0\oplus(E^0\cap E_C)$.
(iii) By Proposition~\ref{rfi:prop:rowred}(ii), regular feedback means the leading rows of $R_{\rm full}$ and of $[0\ \ R_\cC]$ are independent, i.e.\ the homogeneous generators of $\operatorname{gr}E_P$ and $0\oplus\operatorname{gr}E_C$ are independent, i.e.\ $\operatorname{gr}E_P\cap(0\oplus\operatorname{gr}E_C) = 0$, which is $G\cap\operatorname{gr}E_C = 0$. If $e\in E^0$ then $(0,e)\in E_P$ and $(0,\operatorname{lf}e) = \operatorname{lf}(0,e)\in\operatorname{gr}E_P$, so $\operatorname{gr}E^0\subseteq G$.
\end{proof}

The module $G$ can be larger than $\operatorname{gr}E^0$: a homogeneous element of $\operatorname{gr}E_P$ with zero $w$-block (a constant combination of leading rows of $R_{\rm full}$ whose $w$-part cancels) need not be the leading form of an element of $E_P$ with zero $w$-block. Example~\ref{rfi:ex:partial-neg} below is the simplest instance.

\begin{theorem}[Regular feedback implementability, partial interconnection]\label{rfi:thm:rfi-partial}
Let $\cN\subseteq\cK\subseteq\cP$. The following are equivalent.
\begin{enumerate}
\item[(a)] $\cK$ is implementable by regular feedback from $\cP_{\rm full}$.
\item[(b)] (A) $\cC_{\rm can}$ is implementable by regular feedback from $\cP_c$ by full interconnection in $\mathcal L^{\mathtt c}$ (Theorem~\ref{rfi:thm:rfi-full} for the pair $\cC_{\rm can}\subseteq\cP_c$), and (B$'$) $\operatorname{gr}E_P\cap(0\oplus\operatorname{gr}E_{\rm can}) = 0\oplus\operatorname{gr}E^0$.
\item[(c)] (A) and (B) $\ n(\cP_{\rm full}\wedge_c\cC_{\rm can}) + n(\cP_c) = n(\cP_{\rm full}) + n(\cC_{\rm can})$.
\end{enumerate}
When these hold, every $\cC\in\mathcal L^{\mathtt c}$ that implements $\cC_{\rm can}$ from $\cP_c$ by regular feedback (full interconnection) implements $\cK$ from $\cP_{\rm full}$ by regular feedback; conversely, if $\cC$ implements $\cK$ by regular feedback then $\cC + \cN_c$ implements $\cC_{\rm can}$ from $\cP_c$ by regular feedback. In particular the row degrees of such a controller with $\cC\supseteq\cN_c$ are $\kappa(\cC_{\rm can})\setminus\kappa(\cP_c)$ and $n(\cC) = n(\cC_{\rm can}) - n(\cP_c)$.
\end{theorem}

\begin{proof}
(a)$\Rightarrow$(b). Let $\cC$ be a regular feedback controller for $\cK$ and put $E_C' := E_C\cap M$, the module of $\cC + \cN_c$. By Lemma~\ref{rfi:lem:partial-modules}(i) and the modular law ($E^0\subseteq M$), $E_C' + E^0 = (E_C+E^0)\cap M = E_{\rm can}$. By (iii), $\operatorname{gr}E_C'\cap G\subseteq\operatorname{gr}E_C\cap G = 0$, hence $\operatorname{gr}E_C'\cap\operatorname{gr}E^0 = 0$. A sum $E_C' + E^0$ with $\operatorname{gr}E_C'\cap\operatorname{gr}E^0 = 0$ is strict: if $u\in E_C'$, $v\in E^0$ had $\deg(u+v)<\deg u = \deg v$, then $\operatorname{lf}u = -\operatorname{lf}v$ would be a nonzero common element. So $\operatorname{gr}E_{\rm can} = \operatorname{gr}E_C'\oplus\operatorname{gr}E^0$: $\operatorname{gr}E^0$ is a graded direct summand of $\operatorname{gr}E_{\rm can}$, which is (A) by Theorem~\ref{rfi:thm:rfi-full}(c), and $\cC+\cN_c$ is a regular feedback controller for $\cC_{\rm can}$ from $\cP_c$ by Proposition~\ref{rfi:prop:rowred}(ii). For (B$'$): by the modular law in graded modules, since $\operatorname{gr}E^0\subseteq G$,
\[
G\cap\operatorname{gr}E_{\rm can} = G\cap(\operatorname{gr}E_C'\oplus\operatorname{gr}E^0) = (G\cap\operatorname{gr}E_C')\oplus\operatorname{gr}E^0 = \operatorname{gr}E^0 .
\]

(b)$\Rightarrow$(a). By (A) and Theorem~\ref{rfi:thm:rfi-full} there is $E_C\subseteq E_{\rm can}$ with $E_C\oplus E^0 = E_{\rm can}$ and $\operatorname{gr}E_C\oplus\operatorname{gr}E^0 = \operatorname{gr}E_{\rm can}$ (the module of any regular feedback controller for $\cC_{\rm can}$ from $\cP_c$). Then $(E_C+E^0)\cap M = E_{\rm can}$, so $\cC = \ker R_\cC$ implements $\cK$ by Lemma~\ref{rfi:lem:partial-modules}(i), and $\operatorname{gr}E_C\cap G = \operatorname{gr}E_C\cap(G\cap\operatorname{gr}E_{\rm can}) = \operatorname{gr}E_C\cap\operatorname{gr}E^0 = 0$ by (B$'$), so the interconnection is regular feedback by (iii).

(b)$\Leftrightarrow$(c). Assume (A) and let $E_C$ be as in the previous paragraph. Since $0\oplus E^0\subseteq E_P$, $E_P + (0\oplus E_{\rm can}) = E_P + (0\oplus E_C)$, and $E_P\cap(0\oplus E_C) = 0\oplus(E^0\cap E_C) = 0$. So $E(\cP_{\rm full}\wedge_c\cC_{\rm can}) = E_P\oplus(0\oplus E_C)$ and, by Proposition~\ref{rfi:prop:rowred}(ii), $n(\cP_{\rm full}\wedge_c\cC_{\rm can}) = n(\cP_{\rm full}) + n(\cC)$ iff $\operatorname{gr}E_P\cap(0\oplus\operatorname{gr}E_C) = 0$ iff $G\cap\operatorname{gr}E_C = 0$. Now $n(\cC) = n(\cC_{\rm can}) - n(\cP_c)$ by strictness, and $G\cap\operatorname{gr}E_{\rm can} = (G\cap\operatorname{gr}E_C)\oplus\operatorname{gr}E^0$ by the modular law, so $G\cap\operatorname{gr}E_C = 0$ iff (B$'$). Hence under (A), (B)$\Leftrightarrow$(B$'$).

The converse statement about $\cC+\cN_c$ was shown in (a)$\Rightarrow$(b); the degree statements follow from Corollary~\ref{rfi:cor:rowdeg-inv} applied to $E_C\oplus E^0 = E_{\rm can}$.
\end{proof}

\begin{remark}
\begin{enumerate}
\item Condition (B$'$) is a finite rank computation, parallel to Theorem~\ref{rfi:thm:rfi-full}(d). Let $(\lambda_i^w,\lambda_i^c)$ and $\delta_i$ be the leading rows and row degrees of $R_{\rm full}$, $\ell_j$ and $k_j$ those of a row-reduced $R_{\cC_{\rm can}}$, and $\varepsilon_m$ the row degrees of a row-reduced $R_{\cP_c}$. For every $e$ let $U_e = \{u : u_i = 0\text{ if }\delta_i>e,\ \sum_iu_i\lambda_i^w = 0\}$; then (B$'$) says
\[
\dim\bigl(U_e\Lambda^c\ \cap\ \operatorname{span}\{\ell_j : k_j\le e\}\bigr) = \#\{m : \varepsilon_m\le e\}\qquad\text{for every }e,
\]
where $\Lambda^c$ is the matrix of the $\lambda_i^c$; the inequality $\ge$ always holds.
\item (B) is the McMillan-degree identity of the Mayer--Vietoris square of $E_P$ and $0\oplus E_{\rm can}$ (sum $E(\cP_{\rm full}\wedge\cC_{\rm can})$, intersection $0\oplus E(\cP_c)$). Strictness of a pair of submodules implies this identity, but the converse fails in general: the degree defect $\gamma(e) = \dim F_e(A+B) - \dim(F_eA + F_eB)$ need not be monotone in $e$ (the instance $R_{\rm full} = \begin{bmatrix}-1&\xi&-1\\ \xi+1&\xi&1\end{bmatrix}$, $\cK = \{\dot w = 0\}$, has $\gamma = (1,0,0,\dots)$. Under (A), however, the sum becomes direct after replacing $E_{\rm can}$ by $E_C$, and for direct sums the identity is equivalent to strictness by Proposition~\ref{rfi:prop:rowred}(ii). This is why the theorem needs both (A) and (B): neither implies the other, and the numerical checks found instances of each holding without the other.
\item In the full interconnection case ($c = w$, $\cP_{\rm full} = \{(w,w) : w\in\cP\}$) one has $\cP_c = \cP$, $\cC_{\rm can} = \cK$, $\cP_{\rm full}\wedge\cC_{\rm can}\cong\cK$, so (B) reads $n(\cK)+n(\cP) = n(\cP)+n(\cK)$ and the theorem reduces to Theorem~\ref{rfi:thm:rfi-full}. When $G = \operatorname{gr}E^0$, (B$'$) is automatic and only (A) remains. When $\cP_c$ is everything ($E^0 = 0$), (A) is automatic and only (B) remains: $n(\cP_{\rm full}\wedge\cC_{\rm can}) = n(\cP_{\rm full}) + n(\cC_{\rm can})$.
\item Regular feedback is the graded version of regular implementability plus the compatibility (B) with the elimination; the ungraded statement recovers the Belur--Trentelman theorem, as follows.
\end{enumerate}
\end{remark}

\begin{proposition}[Regular implementability, recovered]\label{rfi:prop:reg-impl}
Let $\cN\subseteq\cK\subseteq\cP$. The following are equivalent.
\begin{enumerate}
\item[(a)] $\cK$ is regularly implementable from $\cP_{\rm full}$.
\item[(b)] $E^0$ is a direct summand of $E_{\rm can}$, i.e.\ $E_{\rm can}/E^0$ is torsion-free.
\item[(c)] $\cC_{\rm can}$ is regularly implementable from $\cP_c$ by full interconnection.
\item[(d)] $\cK + \cP_{\rm cont} = \cP$ \ (Belur--Trentelman \cite{BT02}).
\end{enumerate}
\end{proposition}
\begin{proof}
(a)$\Leftrightarrow$(b)$\Leftrightarrow$(c). By Lemma~\ref{rfi:lem:partial-modules}(i),(ii), (a) says there is $E_C$ with $E_C\cap E^0 = 0$ and $(E_C+E^0)\cap M = E_{\rm can}$. Replacing $E_C$ by $E_C\cap M$ (which keeps both properties, by the modular law) gives $E_C\oplus E^0 = E_{\rm can}$; conversely any such $E_C$ satisfies both. And $E_C\oplus E^0 = E_{\rm can}$ with $E_C\subseteq\R^{1\times\mathtt c}[\xi]$ is precisely a regular full interconnection $\cP_c\cap\cC = \cC_{\rm can}$. For a finitely generated module over the principal ideal domain $\R[\xi]$, $E^0$ is a direct summand of $E_{\rm can}$ iff $E_{\rm can}/E^0$ is free iff it is torsion-free.

(b)$\Leftrightarrow$(d). Let $E_{\cP} := E(\cP) = \{k : (k,0)\in E_P\}$ (the manifest plant) and $E_{\cP}^{\rm sat} := \{x : ax\in E_{\cP}\text{ for some }a\ne0\}$; then $E(\cP_{\rm cont}) = E_{\cP}^{\rm sat}$ and $E(\cK+\cP_{\rm cont}) = E_K\cap E_{\cP}^{\rm sat}$ \cite{PW98}, so (d) reads $E_K\cap E_{\cP}^{\rm sat} = E_{\cP}$; the inclusion $\supseteq$ is $\cK\subseteq\cP$.
(d)$\Rightarrow$(b): let $e\in E_{\rm can}$ with $ae\in E^0$, $a\ne0$. Write $e = p_c$ with $p\in E_P$, $p_w\in E_K$. Then $(0,ae)\in E_P$, so $ap - (0,ae) = (ap_w,0)\in E_P$, i.e.\ $ap_w\in E_{\cP}$ and $p_w\in E_K\cap E_{\cP}^{\rm sat} = E_{\cP}$. Hence $(p_w,0)\in E_P$ and $(0,e) = p-(p_w,0)\in E_P$: $e\in E^0$.
(b)$\Rightarrow$(d): let $k\in E_K\cap E_{\cP}^{\rm sat}$, $ak\in E_{\cP}$. Since $E_K\subseteq\pi_wE_P$ ($\cN\subseteq\cK$), $k = p_w$ for some $p\in E_P$, and $e := p_c\in E_{\rm can}$. Then $ap - (ak,0) = (0,ae)\in E_P$, so $ae\in E^0$, and torsion-freeness gives $e\in E^0$, i.e.\ $(0,e)\in E_P$; thus $(k,0) = p-(0,e)\in E_P$ and $k\in E_{\cP}$.
\end{proof}

So the two levels of the problem have parallel answers: regular implementability is \enquote{$E^0$ is a direct summand of $E_{\rm can}$}, and regular feedback implementability is \enquote{$\operatorname{gr}E^0$ is a graded direct summand of $\operatorname{gr}E_{\rm can}$, and the elimination defect does not meet $\operatorname{gr}E_{\rm can}$}. The first condition is a saturation statement, which is why it can be rewritten as $\cK+\cP_{\rm cont} = \cP$; the second has no analogue of $\cP_{\rm cont}$ because the associated graded is not a behaviour.

\begin{example}[The plant's degree carried by $c$]\label{rfi:ex:partial-neg}
$\mathtt w = \mathtt c = 1$, $\cP_{\rm full} = \{w = \dot c\}$, $R_{\rm full} = [1\ \ -\xi]$, $n(\cP_{\rm full}) = 1$. The manifest plant is everything, $\cN = 0$, so every $\cK = \ker k(\tfrac d{dt})$ is implementable, and regularly so ($E^0 = 0$). Here $\cP_c$ is everything, $\cN_c = \{\dot c = 0\}$, $\cC_{\rm can} = \ker(k\xi)$. The leading row of $R_{\rm full}$ is $(0\mid-1)$, so $G = (\xi)$ while $\operatorname{gr}E^0 = 0$; the unit of McMillan degree is accounted for neither by $\cN$ nor by $\cP_c$ ($n(\cN) = n(\cP_c) = 0$). (A) holds trivially. (B): $\cP_{\rm full}\wedge\cC_{\rm can} = \{w = \dot c,\ k(\tfrac d{dt})\dot c = 0\}$ has McMillan degree $\deg k + 1 = n(\cC_{\rm can})$, while $n(\cP_{\rm full}) + n(\cC_{\rm can}) = \deg k + 2$. So (B) fails for every $k\ne$ const, and \emph{no proper $\cK$ is implementable by regular feedback}, although every $\cK$ is regularly implementable. Directly: every controller row has leading row supported on $c$, as does the plant row, so the stacked leading matrix has rank $1$. The unit of McMillan degree of the plant lives in $c$ and is invisible from $w$; a regular feedback controller would have to leave it alone, which no controller on $c$ can.
\end{example}

\begin{example}[Positive case]\label{rfi:ex:partial-pos}
$\mathtt w = 1$, $\mathtt c = 2$, $\cP_{\rm full} = \{w = c_1,\ \dot c_1 = c_2\}$, $R_{\rm full} = \begin{bmatrix}1&-1&0\\0&\xi&-1\end{bmatrix}$ (row reduced, $n = 1$). $\cP_c = \{\dot c_1 = c_2\}$, $E^0 = \R[\xi](\xi,-1)$, $G = \R[\xi](\xi,0) = \operatorname{gr}E^0$, so (B$'$) is automatic. For $\cK = \ker k$, $\deg k = d\ge2$, $k = \xi q + k_0$: $E_{\rm can} = \operatorname{span}\{(\xi,-1),(k_0,q)\}$, row degrees $\{1,d-1\}$, and $\operatorname{gr}E^0 = \R[\xi](\xi,0)$ is a graded summand of $\operatorname{gr}E_{\rm can} = \operatorname{span}\{(\xi,0),(0,\xi^{d-1})\}$: (A) holds. A regular feedback controller is $\cC = \{k_0c_1 + q(\tfrac d{dt})c_2 = 0\}$, $n(\cC) = d-1 = n(\cC_{\rm can}) - n(\cP_c)$; the stacked leading matrix is $\begin{bmatrix}1&-1&0\\0&1&0\\0&0&\operatorname{lc}q\end{bmatrix}$, and $\cP_{\rm full}\wedge\cC$ eliminates to $k_0w + q\dot w = kw = 0$.
\end{example}

\begin{example}[(A) without (B)]\label{rfi:ex:partial-AnotB}
$\mathtt w = 1$, $\mathtt c = 2$, $R_{\rm full} = \begin{bmatrix}1&0&\xi\\0&-1&-1\end{bmatrix}$ ($w + \dot c_2 = 0$, $c_1 + c_2 = 0$; row degrees $1,0$, $n(\cP_{\rm full}) = 1$), $\cK = \ker(\xi-2)$. Then $E^0 = \R[\xi](1,1)$, $E_{\rm can} = \operatorname{span}\{(1,1),(0,\xi^2-2\xi)\}$ with row degrees $\{0,2\}$, and $\operatorname{gr}E^0 = \R[\xi](1,1)$ is a graded summand: (A) holds, $\cC_{\rm can}$ is implementable from $\cP_c$ by regular feedback. But $G\ni(0,\xi)$ (leading row $(0\mid0,1)$ of the first plant row) and $(0,\xi^2)\in\operatorname{gr}E_{\rm can}\cap G\setminus\operatorname{gr}E^0$: (B$'$) fails, and indeed $n(\cP_{\rm full}\wedge\cC_{\rm can}) + n(\cP_c) = 2 + 0 \ne 1 + 2$. $\cK$ is regularly implementable but not by regular feedback (confirmed by exhaustive search over small controllers).
\end{example}

\section{What the structure says}

Regular feedback separates from regularity exactly where the \emph{associated graded} of the interconnection fails to be exact: regularity is independence of the equation modules $E(\cP) = \operatorname{rowspan}R_\cP$ and $E(\cC)$ inside $\R^{1\times\mathtt w}[\xi]$; regular feedback is independence of their associated graded modules $\operatorname{gr}E(\cP)$, $\operatorname{gr}E(\cC)$ (generated by the leading forms) inside $\operatorname{gr}\R^{1\times\mathtt w}[\xi]$. The full-interconnection theorem says that $\cK$ is implementable by regular feedback iff $\operatorname{gr}E(\cP)$ is a graded direct summand of $\operatorname{gr}E(\cK)$, i.e.\ iff the inclusion $E(\cP)\hookrightarrow E(\cK)$ is split \emph{and strict} for the degree filtration. The partial-interconnection theorem says the same thing one level up: the problem on $(w,c)$ is the full problem on $c$ for the pair $(\cP_c,\cC_{\rm can})$, provided the leading forms of plant equations with vanishing $w$-part meet $\operatorname{gr}E(\cC_{\rm can})$ only in $\operatorname{gr}E(\cP_c)$ (condition (B)). Where they do not, the plant carries McMillan degree in the interconnection variable that no controller can respect. Chapter~\ref{ch:diag} places this next to the other three problems.

\chapter{The Lean formalization}\label{ch:formal}

\section{Scope and conventions}

The results of Chapters~\ref{ch:fa}--\ref{ch:rfi} are accompanied by eight modules of the Lean~4 development \texttt{SingularityTracking} (Mathlib, toolchain v4.34.0-rc2): \texttt{Farkas}, \texttt{FarkasMinimal}, \texttt{FarkasDegree}, \texttt{FlatOutput}, \texttt{Observer}, \texttt{NormalForm}, \texttt{Implementability}, and \texttt{LeadingForms} (all \texttt{.lean}). All declarations are proved without \texttt{sorry}, the full project builds (2752 jobs), and every theorem listed below passes the axiom audit (\texttt{AuditOpenProblems.lean}, 59 declarations including the augmentation exact sequence of \texttt{ExactSequence.lean} used in Chapter~\ref{ch:diag}): each depends only on \texttt{propext}, \texttt{Classical.choice}, \texttt{Quot.sound}. 

A ninth module, \texttt{FarkasStatement.lean}, states the two conjectures of Chapter~\ref{ch:fa} as propositions and proves their negations.

The purpose of this chapter is to say precisely which mathematical statements are machine-checked and which rest on paper proofs. A distinction matters throughout: a module may prove an \emph{instance} (this $R$, this $R'$) or it may define the \emph{problem} as a proposition and prove its negation. For Problem~1.9 both are done. For Problems~7.3/P1 and~3.2 only instances and components are checked; the problem statements themselves (an observer on a compact manifold; the homological equations in the book's diagonal coordinates) are not formalized, and the link between the checked components and the problem is a paper argument. For Problem~1.10 (Chapter~\ref{ch:rfi}) two layers are checked: the module-theoretic layer of the partial interconnection case, including the recovery of the Belur--Trentelman theorem, and the filtered layer, the lifting lemma \enquote{leading forms detect bases} and with it the full-interconnection theorem in the form (a)$\Leftrightarrow$(c); see Section~\ref{sec:formal-rfi}. A formal theorem is only as strong as its statement, so for each result we record the Lean statement in words and note where it is weaker than, or differently phrased from, the paper version.

\section{Problem 1.9}

\paragraph{Objects.} Laurent polynomials are Mathlib's $\R[T;T^{-1}]$ (\texttt{LaurentPolynomial} $\R$, which is an \texttt{AddMonoidAlgebra}). A row acts on sequences $\Z\to\R$ by $(\texttt{act}\ p\ w)(t) = \sum_k p_k\,w(t+k)$. \texttt{Behavior R} is $\{w : \forall i,t,\ 0\le \texttt{act}(R_i)w(t)\}$ for a \texttt{Fin g}-indexed family of rows; \texttt{Nonneg h} says all coefficients of $h$ are nonnegative; \texttt{BehaviorEI E I} is the mixed equality/inequality behaviour.

\paragraph{Checked statements.}
\begin{itemize}
\item \texttt{const\us of\us nonneg\us concave}: a sequence on $\Z$ that is nonnegative and has nonpositive second difference is constant.
\item \texttt{behavior\us subset}: $\cB([1;-(1-T)^2])\subseteq\cB(1-T)$.
\item \texttt{no\us certificate}: there are no nonnegative $h_1,h_2$ with $1 - T = h_1\cdot1 + h_2\cdot(-(1-T)^2)$. The proof evaluates at $T = 1$ through \texttt{AddMonoidAlgebra.lift} and cancels $1-T$ in the domain $\R[T;T^{-1}]$, exactly as in Theorem~\ref{fa:thm:main}.
\item \texttt{farkas\us conjecture\us false}: the conjunction of the two, with $H$ ranging over all \texttt{Fin 2}-indexed nonnegative families. This is Theorem~\ref{fa:thm:main}.
\item \texttt{obstruction\us of\us hidden\us equality}: Proposition~\ref{fa:prop:obstruction}, with the implicit equality and the absence of nonnegative syzygies as hypotheses.
\item \texttt{behaviorEI\us rep1}, \texttt{rep2}, \texttt{rep3}: the three representations of Proposition~\ref{fa:prop:minimal} all equal the nonnegative constants.
\item \texttt{no\us one\us row}: no representation with at most one row equals the nonnegative constants. The argument for a single inequality uses the exponential sequences $\pm2^t$ and evaluation at $T = 2$ (\texttt{evAt}), not the surjectivity of $r(\sigma)$ used in the paper; the two arguments are different and the formal one is shorter.
\item \texttt{minimality\us conjecture\us false}: \texttt{rep1} is minimal and its second row is an implicit equality; \texttt{rep2} and \texttt{rep3} are minimal with the same behaviour. The equality counts are the index types \texttt{Fin 0}, \texttt{Fin 1}, \texttt{Fin 1}.
\item \texttt{topDeg\us bound}: the degree bound of Proposition~\ref{fa:prop:claimA}: if the leading coefficients of $a$ and $b$ have the same sign and $(p_1+cb,\ p_2-ca)$ is nonnegative with $p_1,p_2$ supported in $(-\infty,N']$, then $\deg_{\rm top}c + \min(\deg_{\rm top}a,\deg_{\rm top}b)\le N'$. The top degree is \texttt{Finset.max'} of the support and the product rule for top coefficients is \texttt{AddMonoidAlgebra.coeff\us mul\us add\us of\us uniqueAdd}.
\end{itemize}

\paragraph{The conjectures as propositions.} The statements above are instances. \texttt{FarkasStatement.lean} defines the conjectures themselves for Laurent polynomial matrices with $q$ columns acting on $\R^q$-valued sequences (\texttt{Row q}, \texttt{VSeq q}, \texttt{BehaviorM}, \texttt{NonnegM}, \texttt{mulRows}):
\begin{itemize}
\item \texttt{FarkasConjecture} $:=\ \forall q\,g\,g'\,R\,R',\ \cB(R)\subseteq\cB(R')\iff\exists H\ge0,\ R' = HR$; and \texttt{farkas\us conjecture\us is\us false}~$: \neg\,$\texttt{FarkasConjecture}. This is Problem~\ref{prob:farkas} as printed, refuted.
\item \texttt{MinimalityConjecture} $:=$ (every minimal representation has no implicit equality) $\wedge$ (two minimal representations of one behaviour have equal numbers of equalities and of inequalities), for all $q$; and \texttt{minimality\us conjecture\us is\us false}, \texttt{minimality\us uniqueness\us false}. The thesis's conjecture (\enquote{minimal iff no redundant equations and no implicit equalities; counts unique}) implies both conjuncts, so its negation follows; the \enquote{if} direction of the thesis statement is not formalized because it is not needed for the refutation.
\end{itemize}
The reduction from the $q = 1$ matrix statements to the scalar rows of \texttt{Farkas.lean} is the bijection $w\mapsto(t\mapsto w(t)_0)$ (\texttt{toSeq}, \texttt{ofSeq}).

\paragraph{Not checked.} The bipolar reduction (Proposition~\ref{fa:prop:bipolar}), the Kakutani--Klee lemma, the description of $\cK(R)$ for the example (Proposition~\ref{fa:prop:green}), and the passage from the degree bound to closedness. The first two need duality theory for locally convex spaces that Mathlib does not have in the required form; the third is a finite Farkas argument with Green's functions that could be formalized but was not. The mirror of \texttt{topDeg\us bound} for trailing coefficients is omitted as it follows by $T\mapsto T^{-1}$.

\section{Problem 7.3/P1}

\paragraph{Objects.} The model is $\dot x = 1$, $y = g(x)$ with $g$ Mathlib's \texttt{expNegInvGlue} ($g(x) = e^{-1/x}$ for $x>0$, $0$ otherwise), so $L_f^kg = g^{(k)} = $ \texttt{iteratedDeriv k expNegInvGlue}.

\paragraph{Checked statements.}
\begin{itemize}
\item \texttt{iteratedDeriv\us expNegInvGlue}: $g^{(k)}(x) = p_k(1/x)\,g(x)$ with $p_0 = 1$, $p_{k+1} = X^2(p_k - p_k')$ (\texttt{flatPoly}); from Mathlib's \texttt{hasDerivAt\us polynomial\us eval\us inv\us mul}.
\item \texttt{degree\us flatPoly}: $\deg p_k = 2k$. \texttt{iteratedDeriv\us expNegInvGlue\us zero}: $g^{(k)}(0) = 0$.
\item \texttt{no\us lipschitz\us bound}: no $L$ and $\delta>0$ with $|g^{(n+1)}(x)|\le L\sum_{k\le n}|g^{(k)}(x)|$ on $(0,\delta)$.
\item \texttt{no\us lipschitz\us output\us relation}: no globally Lipschitz $\varphi:\R^{n+1}\to\R$ with $g^{(n+1)} = \varphi(g,\dots,g^{(n)})$ on $[0,\delta)$, for any $n$. A $\varphi$ Lipschitz only near $0$ is covered by shrinking $\delta$ and extending with McShane's theorem.
\end{itemize}

\paragraph{The observer half: \texttt{Observer.lean}.} The convergence proof of Proposition~\ref{ob:prop:observer} splits into an abstract part (pure ODE theory on the circle) and a concrete part (the KKL immersion of the flat-output system is star-shaped). The abstract part is formalized, with the circle modelled as the unit complex numbers, the state as a lifted angle $x : \R\to\R$ (the point $e^{ix}$), and convergence measured in the chordal distance, which is a distance compatible with the topology and so is one of the distances allowed by Jouan's definition.
\begin{itemize}
\item \texttt{tendsto\us trackODE\us of\us mem\us Ico}: the circle lemma. If $\dot\epsilon = -K(1-\cos\epsilon)$ on $[t_1,\infty)$ with $K>0$ and $\epsilon(t_1)\in[0,2\pi)$, then $\epsilon\to0$. The proof is the one in the paper: $\epsilon$ cannot cross the equilibrium $0$ (uniqueness of ODE solutions, Mathlib's \texttt{ODE\us solution\us unique}), is antitone, has a limit $L\in[0,2\pi)$ by monotone convergence, and $L\ne0$ would give a derivative eventually $\le -c<0$, contradicting $\epsilon\ge0$ by the mean value theorem. \texttt{tendsto\us trackODE}: for arbitrary $\epsilon(t_1)$ the limit is the multiple of $2\pi$ just below it.
\item \texttt{exists\us lift}: a $C^1$ unit-modulus path $a$ on $(t_1,\infty)$ has a $C^1$ angle $\alpha$ with $a = e^{i\alpha}$ and $\dot\alpha = \operatorname{Im}(\bar a\,\dot a)$. Proof: $\alpha$ is defined as an integral, and $a\,e^{-i\alpha}$ has zero derivative because $\bar a\dot a$ is purely imaginary for a unit path (\texttt{deriv\us eq\us of\us norm\us eq\us one}). No covering-space theory is used.
\item \texttt{second\us stage\us tendsto}: if $\omega$ solves $\dot\omega = \operatorname{Im}(\bar a\dot a) + K(1-\operatorname{Re}(\bar a\,e^{i\omega}))$ for $t>t_1$, then $e^{i\omega}-a\to0$. This is the second stage of \eqref{ob:eq:obs2} written without the angle $\Theta$ as an $S^1$-valued function: $\operatorname{Re}(\bar a e^{i\omega}) = \cos(\alpha-\omega)$ and $\operatorname{Im}(\bar a\dot a) = \dot\alpha$.
\item \texttt{two\us stage\us observer}: given $\Theta\circ T = e^{i\,\cdot}$, a first stage with $z(t)-T(x(t))\to0$ and both in a compact set on which $\Theta$ is continuous, and a second stage as above with $a = \Theta\circ z$, then $e^{i\omega(t)} - e^{ix(t)}\to0$. The hypotheses are stated for \emph{every} solution on $(t_1,\infty)$, with no relation between the initial conditions of plant and observer; this is clause (J4) of the verification paragraph in Chapter~\ref{ch:ob}, and the theorem is what makes the globality claim machine-checked at the level of the ODE argument.
\end{itemize}

\paragraph{Not checked.} Everything about the actual system of Theorem~\ref{ob:thm:main}: the two-sided flat reparametrization $\rho$, the composition with the embedding $\gamma$, the asymptotics of Lemma~\ref{ob:lem:asymp} (the formal model has polynomial degrees $2k$ in place of $3(k-1)$; the mechanism, strictly increasing degrees, is the same), and the concrete half of the observer: the quantitative star-shapedness estimate for the averaged curve $T_\lambda$ and the degree argument identifying $\arg\circ T_\lambda$ as a diffeomorphism of the circle (Lemma~\ref{ob:lem:embedding}), which together supply the hypotheses $\Theta\circ T = e^{i\,\cdot}$ and the compact annulus of \texttt{two\us stage\us observer}; and the existence of solutions of the observer equation (the formal theorems take a solution as given and use only uniqueness). Also not formalized is the metric-independence clause (J5), which the paper proves by compactness; the formal statement fixes the chordal distance. The formal content is therefore: \emph{in the model, the obstruction holds}, and \emph{for any system supplying a star-shaped immersion, the two-stage observer converges globally}. That the flat-output system supplies one is a paper theorem.

\section{Problem 3.2}

\paragraph{Objects.} Polynomials are \texttt{MvPolynomial}~$\N$~$\R$; a system with $n$ states uses $X_0,\dots,X_{n-1}$, recorded as $\texttt{vars}\ p\subseteq\texttt{range}\ n$. The drift Lie derivative in Brunovsky coordinates is $\texttt{drift}\ n\ p = \sum_{k+1<n}X_{k+1}\,\partial_kp$. The kernel equations are the structure \texttt{KernelEqs n h a b}: $\texttt{drift}(h_i) = h_{i+1}$ and $\partial_{n-1}h_i = 0$ for $i+1<n$, $\texttt{drift}(h_{n-1})+\alpha = 0$, $\partial_{n-1}h_{n-1}+\beta = 0$.

\paragraph{Checked statements.}
\begin{itemize}
\item \texttt{drift\us mul}, \texttt{drift\us X}, \texttt{vars\us drift}: \texttt{drift} is a derivation, sends $X_m$ to $X_{m+1}$, and raises the variable range by one.
\item \texttt{drift\us iter\us decomp}: if $\texttt{vars}\ p\subseteq\{0,\dots,k\}$ and $1\le j$, $k+j+1\le n$, then $\texttt{drift}^{j}p = R + X_{k+j}\,\partial_kp$ with $\texttt{vars}\ R\subseteq\{0,\dots,k+j-1\}$. This is the combinatorial core of Proposition~\ref{nf:prop:kernel}.
\item \texttt{first\us component\us in\us X0}: if $\partial_{n-1}(\texttt{drift}^jp) = 0$ for all $j\le n-2$, then $\texttt{vars}\ p\subseteq\{0\}$.
\item \texttt{homogeneous\us X0\us eq\us smul}: a homogeneous polynomial of degree $m$ in $X_0$ alone is $c\cdot X_0^m$.
\item \texttt{kernel\us determined}: a kernel element is determined by $h_0$: $h_i = \texttt{drift}^ih_0$, $\alpha = -\texttt{drift}^{n}h_0$, $\beta = -\partial_{n-1}\texttt{drift}^{n-1}h_0$.
\item \texttt{kernel\us le\us one}: two degree-$m$ kernel elements whose $h_0$ have the same coefficient of $X_0^m$ coincide (the kernel is at most one-dimensional).
\item \texttt{kernel\us element}: $h_0 = X_0^m$ with the prolonged $\alpha,\beta$ satisfies \texttt{KernelEqs} (the kernel is nonempty); \texttt{kernel\us element\us ne\us zero}.
\item \texttt{adDiag\us monomial}: with $A = \operatorname{diag}\lambda$, component $i$ of $\ad_{Ax}$ sends the monomial $x^k$ to $(\langle k,\lambda\rangle-\lambda_i)\,x^k$; \texttt{adDiag\us eigenvalue\us ne\us zero} under \texttt{Nonresonant}.
\item \texttt{finrank\us coker}: for a linear map between finite-dimensional spaces with one-dimensional kernel, $\dim\coker + \dim\mathrm{source} = \dim\mathrm{target}+1$.
\end{itemize}

\paragraph{Not checked.} The identification of the kernel equations with the homological equations of Problem~\ref{prob:dev} (a change of coordinates from the nonresonant diagonal form to Brunovsky form), the binomial dimension counts of $\cH_m$ and $P_m$, and hence the closed formula of Corollary~\ref{nf:cor:formula} as a formula in $n,m$; the invertibility of $\ad_{Ax}$ on the whole homogeneous space rather than on each monomial; and the comparison with Kang's normal form, which is a literature statement.

\section{Problem 1.10}\label{sec:formal-rfi}

\paragraph{Objects.} The module \texttt{Implementability.lean} formalizes the module-theoretic layer of Section~\ref{rfi:sec:rfi-partial}, over an arbitrary commutative ring $R$ and arbitrary $R$-modules $W$, $C$ (the behavioural case is $R = \R[\xi]$, $W = \R^{1\times\mathtt w}[\xi]$, $C = \R^{1\times\mathtt c}[\xi]$). The plant is a submodule \texttt{EP} of $W\times C$ and the desired behaviour a submodule \texttt{EK} of $W$. The derived modules are defined by \texttt{comap}/\texttt{map} along the coordinate inclusions and projections: \texttt{hidden EP} $= E^0$, \texttt{cProj EP} $= M$, \texttt{manifest EP} $= E(\cP)$, \texttt{wProj EP} $= E(\cN)$, \texttt{canonical EP EK} $= E_{\rm can}$, \texttt{controlled EP EC} $= E_P + (0\oplus E_C)$, \texttt{controlledManifest EP EC} $= E(\pi_w(\cP_{\rm full}\wedge_c\cC))$. \texttt{Implements}, \texttt{Regular} ($E_P\cap(0\oplus E_C) = 0$) and \texttt{RegularlyImplementable} are the obvious predicates; \texttt{TorsionFreeQuot A B} is the elementwise statement that $A/B$ is torsion-free, and \texttt{sat P} $= \{x : \exists a\ne0,\ ax\in P\}$.

\paragraph{Checked statements.}
\begin{itemize}
\item \texttt{regular\us iff}: Lemma~\ref{rfi:lem:partial-modules}(ii), $E_P\cap(0\oplus E_C) = 0\iff E_C\cap E^0 = 0$.
\item \texttt{implements\us iff}: Lemma~\ref{rfi:lem:partial-modules}(i), under \texttt{EK $\le$ wProj EP} ($\cN\subseteq\cK$) and \texttt{manifest EP $\le$ EK} ($\cK\subseteq\cP$): $\cC$ implements $\cK$ iff $(E_C+E^0)\cap M = E_{\rm can}$.
\item \texttt{regularlyImplementable\us iff\us compl}: regular implementability iff $E^0$ has a complement in $E_{\rm can}$ (the reduction to the $c$-space; the modular law is Mathlib's \texttt{sup\us inf\us assoc\us of\us le}).
\item \texttt{torsionFreeQuot\us iff\us saturation}: Proposition~\ref{rfi:prop:reg-impl}(b)$\Leftrightarrow$(d), $E_{\rm can}/E^0$ torsion-free iff $E_K\cap\operatorname{sat}E(\cP)\subseteq E(\cP)$, under $\cN\subseteq\cK$.
\item \texttt{torsionFreeQuot\us of\us compl}: a complement forces a torsion-free quotient (needs $R$ a domain and $C$ torsion-free).
\item \texttt{exists\us compl\us of\us torsionFreeQuot}: over a principal ideal domain with $E_{\rm can}$ finitely generated, a torsion-free quotient yields a complement; the proof is Mathlib's \texttt{Module.free\us of\us finite\us type\us torsion\us free'} (structure theorem) followed by \texttt{Module.projective\us lifting\us property} to split the quotient map.
\item \texttt{regularlyImplementable\us iff\us saturation}: the chain assembled: over a PID, with $E_{\rm can}$ finitely generated and $\cN\subseteq\cK\subseteq\cP$, $\cK$ is regularly implementable iff $E_K\cap\operatorname{sat}E(\cP) = E(\cP)$. With $\operatorname{sat}E(\cP) = E(\cP_{\rm cont})$ and $E(\cK+\cP_{\rm cont}) = E_K\cap E(\cP_{\rm cont})$ (paper facts), this is the Belur--Trentelman theorem.
\end{itemize}

\paragraph{The filtered layer: \texttt{LeadingForms.lean}.} Over a field $k$, rows are $V = $ \texttt{Fin n $\to$ k[X]}. Everything is coefficientwise, so that Mathlib's degree API is never needed: \texttt{coeffRow e} is the $k$-linear map $v\mapsto(\text{coefficient of }X^e\text{ in }v_j)_j$; \texttt{bounded e} is the $k$-submodule $F_eV$ of rows whose coefficients of $X^m$ vanish for $m>e$; for a $k[X]$-submodule $E$, \texttt{grDeg E e} $:= $ \texttt{coeffRow e} $(E\cap F_eV)$ is the degree-$e$ component of $\operatorname{gr}E$ (a homogeneous row $\ell X^e$ is identified with $\ell\in k^n$); a family \texttt{x : $\iota\to$ V} with degrees \texttt{d} is \texttt{IsRowReduced} if $x_i\in F_{d_i}V$ and the leading rows \texttt{coeffRow (d i) (x i)} are $k$-independent; \texttt{gradedSpan x d e} $= \operatorname{span}_k\{L_i : d_i\le e\}$ is the degree-$e$ component of the graded module generated by the leading forms. Checked:
\begin{itemize}
\item \texttt{coeffRow\us sum\us eq} (predictable degree, the second half of (F1) in Chapter~\ref{ch:diag}): if $\deg a_i + d_i\le m$ for all $i$ then $\texttt{coeffRow m}(\sum a_i x_i) = \sum_i (a_i)_{m-d_i}L_i$. Proof: Mathlib's \texttt{coeff\us mul\us add\us eq\us of\us natDegree\us le}.
\item \texttt{natDegree\us add\us le\us of\us sum\us mem\us bounded}: if $\sum a_ix_i\in F_eV$ for a row-reduced family then $\deg a_i + d_i\le e$ for all nonzero $a_i$ (no cancellation of leading terms). \texttt{linearIndependent\us of\us isRowReduced}: a row-reduced family is $k[X]$-independent.
\item \texttt{grDeg\us span\us eq}: $\operatorname{gr}_e(\operatorname{span}x) = \operatorname{span}_k\{L_i : d_i\le e\}$, i.e.\ the leading forms of a row-reduced family generate the associated graded of its span.
\item \texttt{span\us eq\us of\us grDeg\us le} (\textbf{the lifting lemma}, first half of (F1)): a row-reduced family inside $E$ whose leading forms generate $\operatorname{gr}E$ degreewise spans $E$. The proof is the induction on degree of Chapter~\ref{ch:diag}, with the reduction step \texttt{exists\us reduce} subtracting $\sum_{d_i\le e}c_i X^{e-d_i}x_i$. Together: \texttt{isBasis\us iff\us grDeg\us eq}, a row-reduced family in $E$ is a basis of $E$ iff $\operatorname{gr}_eE = \operatorname{span}_k\{L_i : d_i\le e\}$ for every $e$.
\item \texttt{grDeg\us mono}: $\operatorname{gr}_eE\subseteq\operatorname{gr}_{e+1}E$ (multiplication by $X$).
\item \texttt{isRowReduced\us sum\us iff\us grDeg\us inf} (Proposition~\ref{rfi:prop:rowred}(ii) in basis-free form): for row-reduced bases $x$ of $E_1$ and $y$ of $E_2$, the union is row reduced iff $\operatorname{gr}_eE_1\cap\operatorname{gr}_eE_2 = 0$ for all $e$. This is the \enquote{strict direct sum} characterization of regular feedback; the McMillan-degree formulation is not formalized (see below).
\item \texttt{RegularFeedbackController EK x d y d'}: a family $y$ in $E_K$ with $x\cup y$ row reduced and spanning $E_K$ (a regular feedback controller for $E_K$ relative to the row-reduced plant basis $x$, full interconnection). \texttt{HomogeneousCompletion EK x d m e'}: constant rows $m_l$ of degrees $e_l$, each in $\operatorname{gr}_{e_l}E_K$, independent together with the $L_i$, generating $\operatorname{gr}E_K$ together with them, i.e.\ condition (c) of Theorem~\ref{rfi:thm:rfi-full} without the matrix bookkeeping.
\item \texttt{exists\us controller\us of\us completion} ((c)$\Rightarrow$(a)): a homogeneous completion lifts to a controller, by choosing lifts of the $m_l$ in $E_K$ and applying the lifting lemma. \texttt{HomogeneousCompletion.of\us controller} ((a)$\Rightarrow$(c)). \texttt{exists\us controller\us iff}: Theorem~\ref{rfi:thm:rfi-full}, (a)$\Leftrightarrow$(c), as a single statement.
\item \texttt{grDeg\us inf\us eq\us of\us controller} (necessity of (d)): if a controller exists then $\operatorname{gr}_eE_K\cap\operatorname{span}\{L_i\} = \operatorname{gr}_eE_P$ for all $e$.
\end{itemize}

\paragraph{Not checked.} The McMillan degree itself, i.e.\ the identification $n(\cB) = \sum$ generator degrees of $\operatorname{gr}E(\cB)$ (Rapisarda--Willems, Forney), and hence the formulation of regular feedback through $n(\cP\wedge\cC) = n(\cP)+n(\cC)$; the formal notion of regular feedback is the strict direct sum, which the paper shows to be equivalent. The existence of row-reduced bases (Forney) is not formalized; all statements take a row-reduced plant basis as hypothesis. The direction (d)$\Rightarrow$(c) of Theorem~\ref{rfi:thm:rfi-full}, the construction of a homogeneous completion from the rank condition (Rado, or \enquote{torsion-free implies free}), is not formalized: the checked theorem is (a)$\Leftrightarrow$(c) plus (a)$\Rightarrow$(d). The graded half of the partial-case theorem (condition (B$'$)) is not formalized, although the two modules together contain all its ingredients. The dictionary between behaviours and modules (elimination, $E(\cP_{\rm cont}) = \operatorname{sat}E(\cP)$) is not formalized either; the Lean statements are about submodules and the behavioural reading is on paper.

\section{Reading the ledger}

Across the four problems, the parts that resisted formalization are of two kinds. Some are infrastructure gaps in Mathlib: duality for locally convex spaces, tubular neighbourhoods, basis counts of graded polynomial spaces. Others are the analytic halves of the arguments: the star-shapedness of the KKL immersion, the Green's function construction; the ODE half of the observer, by contrast, formalized once uniqueness replaced existence and the angle lift was built from an integral rather than from covering spaces. The parts that formalized readily are the algebraic obstructions: evaluation homomorphisms and cancellation in a domain, growth of polynomial degrees, variable bookkeeping under a derivation, and the lattice of submodules in the implementability problem, where the Lean statements hold over any commutative ring and the behavioural case is one instance. This matches the thread of Chapter~\ref{ch:intro}: in each problem the obstruction is algebraic and sits on the degenerate set, and the algebra is what the machine checks.

\end{document}